\documentclass{article}
\usepackage{amsmath,amsthm,amssymb,mathtools}
\usepackage[margin=3cm]{geometry}
\usepackage{charter}
\usepackage[affil-sl]{authblk}

\usepackage[utf8]{inputenc}
\usepackage[T1]{fontenc}
\usepackage[english]{babel}
\usepackage{amsfonts}
\usepackage{nicefrac}
\usepackage{graphicx}
\usepackage{xcolor}
\usepackage{dsfont}

\definecolor{darkred}{RGB}{150, 0, 0}
\definecolor{darkgreen}{RGB}{0, 150, 0}
\definecolor{darkblue}{RGB}{0, 0, 150}
\usepackage[breaklinks=true,colorlinks,citecolor=darkgreen,linkcolor=darkred,urlcolor=darkblue,bookmarks=false]{hyperref}

\usepackage{tikz}
\usetikzlibrary{shapes.geometric}
\usetikzlibrary{arrows.meta}
\usetikzlibrary{calc}
\usepackage[inline]{enumitem}
\usepackage{xspace}
\usepackage{tabularx}
\usepackage{booktabs}
\usepackage{capt-of}
\usepackage{supertabular}
\usepackage[binary-units=true]{siunitx}
\usepackage{enumitem}

\newcommand{\R}{\mathds{R}}
\newcommand{\Q}{\mathds{Q}}
\newcommand{\Z}{\mathds{Z}}
\newcommand{\B}[1]{\{0,1\}^{#1}}
\newcommand{\cC}{\mathcal{C}}
\newcommand{\SSTcone}{\mathcal{K}}

\DeclareMathOperator{\conv}{conv}

\newcommand{\group}{\Gamma}
\newcommand{\perm}{\gamma}
\newcommand{\pig}{\perm\in\group}
\newcommand{\sym}[1]{\mathcal{S}_{#1}}
\newcommand{\stab}[2]{\mathrm{stab}(#1,#2)}
\newcommand{\lex}[2]{#1 \preceq_{\mathrm{lex}}#2}
\newcommand{\orbit}[2]{\mathrm{orb}(#1,#2)}
\newcommand{\inv}[1]{{#1}^{-1}}

\newcommand{\st}{\,:\,}
\newcommand{\T}{^\top}

\newcommand{\sprod}[2]{{#1}\T{#2}}
\newcommand{\define}{\coloneqq}
\newcommand{\card}[1]{\lvert{#1}\rvert}
\DeclareMathOperator{\argmax}{argmax}
\newcommand{\tpgraph}{TP-graph\xspace}
\newcommand{\tpgraphs}{TP-graphs\xspace}

\theoremstyle{plain}
\newtheorem{theorem}{Theorem}[section]
\newtheorem{lemma}[theorem]{Lemma}
\newtheorem{proposition}[theorem]{Proposition}
\newtheorem{corollary}[theorem]{Corollary}
\newtheorem*{claim*}{Claim}
\newtheorem{claim}[theorem]{Claim}
\newtheorem{observation}[theorem]{Observation}
\theoremstyle{definition}
\newtheorem{remark}[theorem]{Remark}
\newtheorem{definition}[theorem]{Definition}
\newtheorem{example}[theorem]{Example}

\tikzset{
    triangleA/.style={
        anchor=corner 3,
        draw,
        shape border rotate=180,
        regular polygon,
        regular polygon sides=3,
        fill=white,
        node distance=2cm,
        minimum height=4em
    }
}

\tikzset{
    triangleB/.style={
        anchor= corner 1,
        draw,
        shape border rotate=80,
        regular polygon,
        regular polygon sides=3,
        fill=white,
        node distance=2cm,
        minimum height=4em
    }
}

\tikzset{
    triangleC/.style={
        anchor=corner 2,
        draw,
        shape border rotate=160,
        regular polygon,
        regular polygon sides=3,
        fill=white,
        node distance=2cm,
        minimum height=4em
    }
}

\tikzset{
    triangleD/.style={
        anchor=corner 3,
        draw,
        shape border rotate=-60,
        regular polygon,
        regular polygon sides=3,
        fill=white,
        node distance=2cm,
        minimum height=4em
    }
}

\tikzset{
    pentagon/.style={
        anchor=corner 3,
        shape border rotate=0,
        regular polygon,
        regular polygon sides=5,
        fill=white,
        node distance=2.5cm,
        minimum height=6em
    }
}

\title{The Impact of Symmetry Handling for the\\ Stable Set
Problem via Schreier-Sims Cuts\thanks{This article was partially funded by CMM FB210005 and by Fondecyt Regular Nr. 1221460.}}

\author[1]{Christopher Hojny}
\author[2]{Annika J\"ager}
\author[2]{Marc E. Pfetsch}
\author[3]{Jos\'e Verschae}
\affil[1]{%
  Eindhoven University of Technology,
  \emph{email} c.hojny@tue.nl
}
\affil[2]{%
  TU Darmstadt, Germany,
  \emph{email} \{jaeger|pfetsch\}@opt.tu-darmstadt.de
}
\affil[3]{%
  Pontificia Universidad Católica de Chile,
  Santiago de Chile, Chile,
  \emph{email} jverschae@uc.cl
}

\begin{document}

\maketitle

\begin{abstract}
  \noindent
  Symmetry handling inequalities (SHIs) are an appealing and popular tool for handling symmetries in integer programming.
  Despite their practical application, little is known about their interaction with optimization problems.
  This article focuses on Schreier-Sims (SST) cuts, a recently introduced family of SHIs, and investigates their impact on the computational and polyhedral complexity of optimization problems. Given that SST cuts are not unique, a crucial question is to understand how different constructions of SST cuts influence the solving process.

  First, we observe that SST cuts do not increase the computational complexity of solving a linear optimization problem over any polytope $P$. However, separating the integer hull of $P$ enriched by SST cuts can be NP-hard, even if $P$ is integral and has a compact formulation.
  We study this phenomenon more in-depth for the stable set problem, particularly for subclasses of perfect graphs. In doing so, we reveal a relation to mixed graphs and exploit results regarding the generalized stable set problem on bidirected graphs. For bipartite graphs, we give a complete characterization of the integer hull after adding SST cuts based on odd-cycle inequalities. For trivially perfect graphs, we observe that the separation problem is still NP-hard after adding a generic set of SST cuts. Our main contribution is to identify a specific class of SST cuts, called stringent SST cuts, that keeps the separation problem polynomial and a complete set of inequalities, namely SST clique cuts, that yield a complete linear description.

  We complement these results by giving SST cuts based presolving techniques and provide a computational study to compare the different approaches. In particular, our newly identified stringent SST cuts dominate other approaches. 
 \smallskip
 
 \noindent
 \textbf{keywords:} symmetry handling, stable set, perfect graph, Schreier-Sims table, integer hull
\end{abstract}

\section{Introduction}

Handling symmetries in integer programs has the goal to speed up the
solution process by avoiding the consideration of symmetric solutions.
Although many techniques have been developed for this goal, only little is
known about the interaction of symmetry handling methods and structures of
the problem to be solved.
Understanding this interplay opens the door for deriving stronger symmetry
handling methods that allow to solve the underlying problem faster than
with generic symmetry handling methods.
We therefore initiate in this article a systematic study of a popular
symmetry handling method, so-called SST cuts, and its interplay with a
classical problem in integer programming and combinatorial optimization,
the stable set problem.

Given an integer program~$\max_{x \in \Z^n}\{\sprod{c}{x} : Ax \leq b\}$, a
symmetry is a permutation~$\perm$ of the coordinate space that transforms
feasible solutions~$x$ into feasible solutions~$\perm(x) \define
(x_{\inv{\perm}(1)},\dots,x_{\inv{\perm}(n)})$ while preserving their
objective value.
The set of all symmetries of an integer program forms a group~$\group$, which
can be handled, among others, by symmetry handling inequalities (SHIs).
SHIs are inequalities that cut off some solutions, while guaranteeing that
for every feasible solution~$x$ there exists a symmetric
solution~$\perm(x)$, $\pig$, that satisfies all SHIs.
A special class of SHIs are SST cuts, which have been proposed by Liberti
and Ostrowski~\cite{LibO14} and Salvagnin~\cite{Sal18} and are derived from the
so-called Schreier-Sims Table of~$\group$.
More concretely, SST cuts are inequalities of the form~$x_j \leq x_i$ for
some carefully selected set of pairs of variable indices~$i$ and~$j$.
They can be constructed in polynomial
time if we are given a set of generators of~$\group$, see
Section~\ref{sec:SSTdefinition} for details, but they are not unique
because the indices can be selected in different ways.
Different parameterizations of SST cuts can thus potentially yield
substantial performance differences when solving an integer program.
The leading question of this article is therefore
\begin{quote}
  \emph{What is the impact of different parameterizations of SST cuts on
    the complexity of the underlying problem to be solved, both in theory
    and practice?}
\end{quote}
Without knowing the structure of the underlying problem to be solved, however,
it is arguably difficult to find a good parameterization.
We therefore focus on the stable set problem and the interplay of SST cuts
with the stable set polytope, because many structural results are known
for this polytope.

The goal of the stable set problem is to find, for an undirected simple
graph~$G = (V,E)$ and weights~$c_v \geq 0$, $v \in V$, a set~$A \subseteq
V$ of maximum weight~$\sum_{v \in A} c_v$ such that~$A$ is stable, i.e.,
no pair of distinct nodes in~$A$ is adjacent in~$G$.
A common integer programming formulation of this problem is
\[
  \max_{x \in \B{V}}\, \Big\{\sum_{v \in V} c_v\, x_v \;:\; x_u + x_v \leq
  1 \text{ for all } \{u, v\} \in E\Big\}
\]
from which the stable set polytope~$P(G) \define \conv\{ x \in \B{V} : x_u + x_v
\leq  1,\; \{u, v\} \in E\}$ is derived; here~$\conv(\cdot)$ denotes the
convex hull operator.
Since the stable set problem is NP-hard, it is unlikely to find complete
linear descriptions of~$P(G)$ for general graphs~$G$.
For the special class of perfect graphs~$G$, i.e., graphs for which on every induced subgraph the clique number and the chromatic number coincide, however, the stable set problem is
polynomial time solvable and the facet structure of~$P(G)$ is completely known.
Perfect graphs therefore form an ideal family of graphs to investigate the
interplay of SST cuts with the stable set polytope.
More concretely, we aim to understand, for a set~$S$ of SST cuts, the
object \[P(G,S) \define \conv\{ x \in P(G) \cap \B{V} : x \text{ satisfies the SST
  cuts~$S$}\},\] for which we formulate the following research questions.
\begin{enumerate}[label=(Q\arabic*),leftmargin=1cm]
\item\label{Q1} What is the complexity of optimizing a linear function over~$P(G,S)$
  when~$G$ is perfect?
\item\label{Q2} Every valid inequality for~$P(G,S)$ can be interpreted as a symmetry
  handling inequality that is compatible with the stable set problem
  on~$G$; facet defining inequalities are thus the strongest such SHIs.
  \begin{enumerate}[label=(Q2\alph*),leftmargin=1cm]
  \item Does the facet structure of~$P(G,S)$ depend on the parameterization
    of~$S$?
  \item Are there parameterizations of~$S$ such that a complete linear
    description of~$P(G,S)$ can be derived?
  \end{enumerate}
\end{enumerate}
Clearly, one would hope that neither the computational nor polyhedral
complexity increases when adding SST cuts. But the answer to the above questions is not immediate in
general, since SST cuts might change the structure of the underlying
problem and thus its complexity.

\paragraph{Outline and Contributions.}
After presenting related literature and terminology, Section~\ref{sec:SSTdefinition}
formally introduces SST cuts.
Section~\ref{sec:complexity} is dedicated to Question~\ref{Q1}.
We show that there is a significant difference between finding an optimal
solution in~$P(G)$ that adheres to SST cuts~$S$ and optimizing a linear
function over~$P(G,S)$.
Namely the former is polynomial time equivalent to solving the stable set
problem (and thus polynomial for perfect graphs), whereas the latter is
NP-hard even if~$G$ is perfect.
As mentioned above, there is thus a rich interplay between SST cuts and the
problem to be solved.

Section~\ref{sec:cuts} aims to understand this interplay from a polyhedral
perspective and provides some answers to Question~\ref{Q2}.
Due to the known facet structure of~$P(G)$ for perfect graphs and the
NP-hardness of optimization problems over~$P(G,S)$ for this class of
graphs, Section~\ref{sec:cuts} puts its focus on two subclasses of perfect
graphs: bipartite graphs and trivially perfect graphs.
For bipartite graphs, we provide a complete linear description
of~$P(G,S)$ for arbitrary sets~$S$ of SST cuts.
For trivially perfect graphs, we show that the parameterization of SST cuts
has a strong impact on the facet structure of~$P(G,S)$ and derive a
specific parameterization for which we can prove a complete linear
description of~$P(G,S)$.
This description is based on combining SST cuts with well-known clique
inequalities for the stable set polytope, yielding so-called SST clique
cuts.
SST clique cuts are also valid for~$P(G,S)$ for general graphs~$G$, and
are thus potentially stronger SHIs than SST cuts alone when solving stable
set problems by branch-and-cut.
We investigate the latter by a thorough computational study in
Section~\ref{sec:ComputationalResults}, which also concludes the article.

\paragraph{Techniques}
To achieve the results of Section~\ref{sec:cuts}, we use two different
approaches.
First, we show that there exists an auxiliary graph~$G'$ such that~$P(G,S)$
is the projection of a face of~$P(G')$.
This face is thus an extended formulation of~$P(G,S)$ and this result holds
for arbitrary graphs~$G$.
If~$G$ is bipartite, however, we can show that~$G'$ is t-perfect, which
allows us to take the known complete linear description of~$P(G')$ and
explicitly compute the projection~$P(G,S)$.

Second, we interpret SST cuts as directed edges that have been added to
$G$, thereby generating a mixed graph.
This allows to connect~$P(G,S)$ to the degree-two inequality
polyhedra derived from mixed graphs investigated by Johnson and
Padberg~\cite{Johnson1982}.
We use Tamura's characterization that provides a linear description of
$P(G,S)$ if and only if a carefully defined auxiliary graph
$\underline{G}_S$ is perfect~\cite{Tamura1997}.
This way, the question of a complete linear description of~$P(G,S)$
translates into the question of how SST cuts must be chosen such
that~$\underline{G}_S$ is perfect.
For trivially perfect graphs~$G$, we show that there is a parameterization~$S$
of SST cuts such that~$\underline{G}_S$ is perfect.
The results of Johnson and Padberg~\cite{Johnson1982} for degree-two
inequality polyhedra then allow us to show the description of~$P(G,S)$ via
SST clique cuts.

The latter analysis is heavily based on structural assumptions
on~$G$ made by~\cite{Tamura1997}.
We show that these assumptions indeed hold when the original graph~$G$ is
modified by incorporating implications of SST cuts into~$G$.
This results in the deletion of certain nodes and addition of certain
edges, which can be interpreted as presolving reductions of the given
graph~$G$.
These reductions are thus of independent interest when solving stable set
problems, because they simplify the graph's structure, and we investigate
the reductions' impact in our computational study.

\paragraph{Literature Review.}

The literature mainly discusses two lines of research to handle
symmetries in binary programs: the addition of (static) symmetry handling
inequalities to the problem formulation, which restricts the search space
of the original problem, and dynamic symmetry handling techniques, which
modify the branch-and-bound algorithm based on symmetry information, see,
among others,
\cite{Mar02,Mar03a,Ostrowski2008,OstrowskiAnjosVannelli2015,OstrowskiEtAl2011,LinderothEtAl2021}.
In this article, we follow the former line of research.

Let~$\sym{n}$ be the symmetric group on~$[n] \define \{1,\dots,n\}$.
A standard technique for deriving SHIs is to
enforce that among a set of symmetric solutions only those should be computed
that are lexicographically maximal.
Here, we say that two solutions~$x$ and~$y$ of a binary program with
symmetry group~$\group \leq \sym{n}$ are \emph{symmetric} if there exists~$\perm
\in \group$ such that~$y = \perm(x)$.
The set of all symmetric solutions of~$x$ is called the \emph{orbit}
of~$x$, denoted~$\orbit{x}{\group} \define \{\perm(x) \st \perm \in
\group\}$.
By overloading notation, we denote the orbits of variable indices~$i \in [n]$
by $\orbit{i}{\group} = \{\perm(i) \st \perm \in
\group\}$.
Moreover, we say that a vector~$x$ is lexicographically smaller or equal to
vector~$y$, denoted~$\lex{x}{y}$, if either~$x = y$ or~$x_i < y_i$ for the
first position~$i$ in which~$x$ and~$y$ differ.
If~$\mathcal{X}$ is the feasible region of a binary program, a valid
symmetry handling approach is to restrict the feasible region
to~$\mathcal{X} \cap \{ x \in \B{n} \st \lex{\perm(x)}{x} \text{ for all } \perm\in\group\}$.

Friedman~\cite{Friedman2007} describes how the lexicographic restriction
can be modeled by linear inequalities, which potentially need exponentially large
coefficients.
Since this might cause numerical instabilities, an alternative set of
inequalities with~$\{0, \pm 1\}$-coefficients is suggested by~\cite{HojnyPfetsch2019}.
The alternative set of inequalities for~$\lex{\perm(x)}{x}$ is derived via
a knapsack polytope associated with~$\conv\{x \in \B{n} \st
\lex{\perm(x)}{x}\}$ for a fixed permutation~$\perm$, the so-called
\emph{symresack}.
Although the alternative set of inequalities is exponentially large,
it can be separated in linear time~\cite{vDoornmalen2024}.
A linear time propagation algorithm for the constraint~$\lex{\perm(x)}{x}$
also exists~\cite{vDoornmaalenHojny2022}.
Moreover, there exists a family of permutations such that in each integer
programming formulation of a symresack, the size of coefficients or the
number of inequalities needs to be exponentially large~\cite{Hojny2020a}.

A drawback of Friedman's approach is that one needs to add a constraint for each permutation in a group.
But for specific symmetry groups, stronger results can be achieved.
If the symmetry group~$\group$ is cyclic, efficient propagation algorithms
for~$\{ x \in \B{n} \st \lex{\perm(x)}{x} \text{ for all }
\perm\in\group\}$ are described in~\cite{vDoornmaalenHojny2022}. If the variables of a binary program can be arranged in a matrix~$X \in
\B{p \times q}$ with~${n = p\cdot q}$, and the symmetry group~$\group$ acts
on the variables by permuting the columns of~$X$, the lexicographic
restriction boils down to sort the columns in lexicographic order, i.e., to
enforce~${\lex{X^{i+1}}{X^i}}$ for all~$i \in [q-1]$, where~$X^i$ denotes
the~$i$-th column of~$X$.
Bendotti et
al.~\cite{BendottiEtAl2021} describe a linear time propagation algorithm
for the set~$\mathcal{X}^{p,q} \define \{ X \in \B{p \times q} \st \lex{X^{i+1}}{X^i} \text{ for all
} i \in [q-1]\}$.
Moreover, if additionally every row of~$X$ has at most/exactly one~1-entry,
a complete linear description of the convex hull of these matrices, the
so-called packing/partitioning orbitope, is known~\cite{KaibelPfetsch2008}
and can be propagated in linear time~\cite{KaibelEtAl2011}. In general, however, a complete linear description
of~$\conv(\mathcal{X}^{p, q})$ is unknown~\cite{KaibelLoos2011} except for the case~$q = 2$
and~$p = 1$~\cite{HojnyPfetsch2019}.
That is, the strongest symmetry handling inequalities are unknown in general.

For general groups, Liberti~\cite{Liberti2012a} suggests to select a single
variable~$x_i$ and to add the inequalities~$x_i \geq x_j$ for all~$j \in
\orbit{i}{\group}$.
These inequalities partially handle symmetries of a single orbit
of~$\group$ and can only be used for a single variable.
By considering a subgroup, this idea can be iterated for a different variable. We detail this approach in the next section as these inequalities form
the main object of interest of this article.

As an extension of the presolving methods of this paper,~\cite{JaegerPfetsch2025} investigates a general framework to derive structure-preserving presolving steps that eliminate all symmetry, including computational results similar as in Section~\ref{sec:ComputationalResults}.
An iterated application of SST presolving appears as the main part of the presented presolving routine.

\section{Schreier-Sims Table Inequalities and Stable Sets}
\label{sec:SSTdefinition}

In this section, we formally introduce SST cuts and how they can be applied
to the stable set problem.
To this end, recall the SHIs~$x_i \geq x_j$ for~$j \in \orbit{i}{\group}$
and a fixed~$i \in [n]$.
These SHIs only handle symmetries on a single variable orbit and are thus
rather weak.
A simple idea to handle more symmetries is therefore to add these
inequalities for several orbits.
In general, this is not possible though, because all optimal solutions can
be cut off.
Liberti and Ostrowski~\cite{LibO14} and Salvagnin~\cite{Sal18} observed,
however, that combining these SHIs is possible when the group~$\group$ is
replaced by suitable subgroups of~$\group$ for subsequent
choices of variables.
This modification requires the concepts of stabilizers.

Let $\group \leq \sym{n}$ be a symmetry group of a binary program. The pointwise
\emph{stabilizer} of a set~$I \subseteq [n]$ is~$\stab{\group}{I} \define \{ \pig \st \perm(i) = i 
\text{ for } i \in I\}$.
If~$\group$ is given by a set of generating permutations~$\Pi \subseteq
\group$, generators of the stabilizer can be computed in time polynomial in~$n$
and~$\card{\Pi}$ via the so-called Schreier-Sims table,
see, e.g., \cite{Seress2003}.
In particular, the number of generators of the stabilizer group is
polynomial.
Since variable orbits $\orbit{i}{\group}$ can also be computed in time
polynomial in $n$ and $\card{\Pi}$, see~\cite{Seress2003}, also orbits of
stabilizer groups can be computed in polynomial time.

Using these concepts, an extended family of SHIs
can be found by the following so-called \emph{SST algorithm}.
It initializes~$\group' \gets \group$ and a sequence~$L \gets \emptyset$.
Afterwards, the following steps are repeated:
\begin{enumerate}[label=(A\arabic*),leftmargin=1cm]
\item\label{alg:A} select~$\ell \in [n] \setminus L$ and compute~$O_\ell \gets
  \orbit{\ell}{\group'}$;
\item\label{alg:B} append $\ell$ to~$L$ and update~$\group' \gets
  \stab{\group'}{L}$;
\item repeat the previous steps until~$L = [n]$.
\end{enumerate}
Due to the above discussion, this algorithm runs in polynomial time.
Moreover, note that we can terminate the SST algorithm once the
group~$\group'$ becomes trivial.
For technical reasons in Section~\ref{sec:triviallyPerfect}, however, we do
not include this check in the SST algorithm.

At termination, $L$ denotes an ordered sequence in $[n]$.
Each element $\ell$ in~$L$ is called
a \emph{leader} and each $f\in O_\ell\setminus\{\ell\}$ is a \emph{follower} of
$\ell$.
In the following, we will also use set notation for $L$, e.g.,
$\ell \in L$ for some element $\ell$ in the sequence.
For $\ell \in L$ and~$f\in O_\ell \setminus \{\ell\}$, a \emph{Schreier-Sims Table
  (SST) cut} is
\begin{align*}
  -x_{\ell} + x_{f} \leq 0,
\end{align*}
which for binary variables says that selecting the
follower~$f$, i.e., $x_f = 1$, implies $x_\ell = 1$.
Due to the modification of the group in each step of the algorithm, all SST
cuts can be used simultaneously to handle symmetries.
We refer the reader to~\cite{LibO14,Sal18} for details on correctness.

We usually refer to a single SST cut by a pair
$(\ell,f)$ with $\ell \in L$ and $f \in O_\ell \setminus \{\ell\}$.
Moreover, $S(L) \define \{(\ell,f) \st \ell\in L,\; f\in O_\ell \setminus \{\ell\}\}$ denotes the set of all SST cuts.
A set~$S \subseteq S(L)$ of SST cuts defines the cone
\[
  \SSTcone(S) \define \{x\in \R^n \st -x_{\ell} + x_f\le 0 \text{ for all } (\ell,f)\in
  S\}.
\]

This cone has attracted significant attention recently. It defines a fundamental domain for the action defined by $\Gamma$~\cite{claus_danielson_fundamental_2021,verschae_geometry_2023}, that is, it is the closure of the set of lexicographically maximal vectors within their orbits. Hence, it provides the best polyhedral approximation of lexicographically maximal vectors. Moreover, although the number of SST cuts can be quadratic, the number of facets of $\SSTcone(S)$ is only $O(n)$~\cite{verschae_geometry_2023}.

To apply SST cuts for solving the stable set problem~$\max\{\sprod{c}{x} :
x \in P(G)\}$ of an undirected graph~$G = (V,E)$, we require a symmetry
group of~$\group \leq \sym{V}$ of the problem, where~$\sym{V}$ is the symmetric group on~$V$.
Recall that a symmetry needs to preserve feasibility of solutions as well as
their objective value.
If all~$c_v$, $v \in V$, are identical, a valid choice for~$\group$ is the
\emph{automorphism group} of~$G$,
i.e., the group consisting of all permutations~$\perm\in\sym{V}$ such
that~$\{\perm(u), \perm(v)\} \in E$ holds for all~$\{u,v\} \in E$.
If not all weights~$c_v$ are identical, then one can select~$\group$ to be
the subgroup of the automorphism group that preserves node weights, i.e.,
$\perm$ additionally needs to satisfy~$c_{\perm(v)} = c_v$ for all~$v \in
V$.
In the remainder of the article, whenever we refer to a symmetry
group~$\group$ of the stable set problem, we always assume that~$\Gamma$
consists of automorphisms of~$G$ that preserve the node weights.

\section{Complexity}
\label{sec:complexity}

The primary goal of SST cuts is to accelerate the solution process.
That is, SST cuts should not make it harder to solve an optimization
problem.
Since we will see below that SST cuts indeed do not increase the
complexity of integer programs (up to a polynomial time overhead), a
secondary goal is to excel the effect of SST cuts by strengthening them.
To find the strongest SHIs for the stable set problem that are based on SST
cuts, we therefore suggest to find valid inequalities for~$P(G,S)$
and to analyze their separation problem, thus providing an answer to
Question~\ref{Q1}.

To find an answer to our question concerning the primary goal, we consider
general integer programs~$\max\{\sprod{c}{x} : Ax \leq b,\; x \in \Z^n\}$
with symmetry groups~$\group$.
We assume that~$\group$ is given by a set~$\Pi$ of generators, and
by~\cite{Seress2003}, we can further assume that~$\card{\Pi} \in
O(\text{poly}(n))$.

\begin{theorem}\label{thm:opt_CS}
  Assume that the program~$\max\,\{\sprod{c}{x} \st A x\le b,\; x\in\Z^n\}$ has~$\group \leq \sym{n}$ as symmetry group and that it can be solved in $T$~time. For leaders $L$ derived from~$\group$ and SST cuts $S = S(L)$, we can solve $\max\,\{\sprod{c}{x} \st Ax\le b,\; x\in \SSTcone(S)\cap\Z^n\}$ in~$T + \text{poly}(n)$ time.
\end{theorem}

\begin{proof}
  Let~$\hat{x}$ be an optimal solution of~$\max\,\{\sprod{c}{x} \st A x\le b,\; x\in\Z^n\}$.
  We construct an optimal solution~$\hat{x}'$ for $\max\,\{\sprod{c}{x} \st Ax\le b,\; x\in \SSTcone(S)\cap\Z^n\}$ in polynomial
  time.
  Consider the first leader $\ell_1 \in L$ and let $i_1 \in \argmax\,\{\hat{x}_i \st i \in O_{\ell_1}\}$ and~$\pig$ be such
  that~$\perm(i_1) = \ell_1$.
  Then, $\perm(\hat{x})$ satisfies the SST cuts $-x_{\ell_1}+x_f\le 0$ for all~$f\in O_{\ell_1}$.
  By replacing $\group$ by the stabilizer of~$\ell_1$ and~$\hat{x}$ by~$\perm(\hat{x})$, we can iterate the procedure for the
  remaining orbits to find a point~$\hat{x}' \in \orbit{x}{\group}$ that
  satisfies all SST cuts.
  Since~$\hat{x}$ is optimal, $\hat{x}'$ is optimal too.
  Because pointwise stabilizers can be computed in polynomial
  time~\cite[Ch. 4]{Seress2003}, $\hat{x}'$  can be constructed in polynomial time.
\end{proof}

We now shift the focus on the secondary goal: the generation of strong SHIs
for the stable set problem.
To this end, we study the separation problem for~$P(G,S)$ for a given set~$S$
of SST cuts that has been derived from a symmetry group~$\group$ of the
stable set problem.
That is, for a given vector~$x^\star \in \R^n$, we study the complexity of
deciding whether~$x^\star \in P(G,S)$.
It turns out that the separation problem is NP-hard even if~$G$ is a
perfect graph.

To prove this statement, we consider the perfect graph~$G_{m,\ell}$ which
consists of the disjoint union of~$\ell$ cliques of size~$m$.
Recall that a set of nodes is a \emph{clique} if every two distinct nodes in the set are adjacent.
We identify the nodes of~$G_{m,\ell}$ with pairs~$(i,j) \in [m] \times
[\ell]$, i.e., $(i,j)$ is the~$i$-th node of the~$j$-th clique, and assign
node~$(i,j)$ the weight~$i$.
The automorphism group of this weighted graph is isomorphic
to~$\sym{\ell}$, where a permutation~$\perm \in \sym{\ell}$ acts on the
nodes by permuting the clique label, i.e., $\perm(i,j) = (i,
\inv{\perm}(j))$.
That is, by interpreting the nodes as entries of a matrix, $\group$ acts on
the nodes by permuting the columns of this matrix.

\begin{theorem}\label{thm:NPhardness}
  Let~$\group$ be the automorphism group of the weighted
  graph~$G_{m,\ell}$.
  There exists a sequence of leaders~$L$ such that the separation problem
  for~$P(G_{m,\ell}, S(L))$ is NP-hard.
\end{theorem}

In our proof, we make use of the equivalence of optimization and
separation~\cite{GroetschelLovaszSchrijver1984}, i.e., we prove that the
complexity of optimizing a linear function over~$P(G_{m,\ell},S(L))$ is
NP-hard.
It is important to note that this does not contradict
Theorem~\ref{thm:opt_CS}.
In Theorem~\ref{thm:opt_CS}, the symmetry group~$\group$ needs to preserve
the objective function, whereas for the general optimization complexity we
need to consider arbitrary objectives (which are not necessarily invariant
under~$\group$).

\begin{proof}
  Let $w\in \Q^{m \times \ell}$ be some weight function.
  We will show that, for a specific choice of leaders $L$, the problem of finding $x \in P(G_{m,\ell},S(L))$ that maximizes $\sprod{w}{x}=\sum_{i=1}^{m}\sum_{j=1}^{\ell} w_{ij}x_{ij}$ is NP-hard. This is enough to show the lemma due to the equivalence of optimization and separation~\cite{GrotschelLovaszSchrijver1993}.

  Let us consider the following instance of 3D-matching. Let $X$, $Y$, $Z$ be three pairwise disjoint sets with $k$ elements each. We are given a collection of sets $\mathcal{T}\subseteq 2^{X \cup Y \cup Z}$, where $T\in \mathcal{T}$ is of the form $\{T_X,T_Y,T_Z\}$ with $T_X\in X$, $T_Y\in Y$, and $T_Z\in Z$. We must decide whether a 3D-matching exists, that is, if there exists a collection $\mathcal{T}'\subseteq \mathcal{T}$ that partitions $X\cup Y \cup Z$, i.e., for any $a\in X\cup Y \cup Z$ there exists exactly one $T\in \mathcal{T}'$ such that $a\in T$. It is well known that this decision problem is NP-complete~\cite{Karp72}.

  \begin{figure}
    \begin{center}
      \includegraphics[width=\textwidth,trim={0 0 0 5},clip]{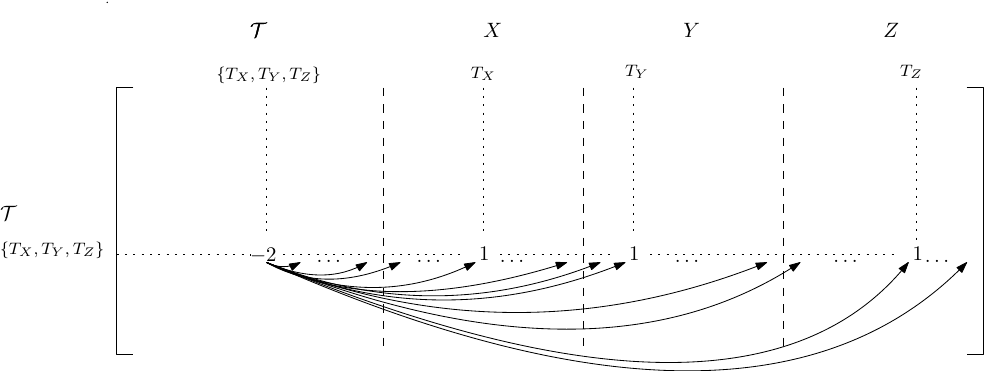}
      \caption{Scheme depicting the indices of a matrix $x\in P(G_{m,\ell})$ and weight matrix $w$. Rows are indexed by elements in $\mathcal{T}$ while columns by $\mathcal{T}\cup X \cup Y \cup Z$. The numbers displayed represent the coefficients of $w$ in the corresponding entry; entries not displayed on row $\{T_X,T_Y,T_Z\}$ are 0. Arrows represent SST cuts, where the tail is the leader and the head the follower.}
      \label{fig:reduction}
    \end{center}
  \end{figure}

  To construct our polytope, we define $m=\card{\mathcal{T}}$ and $\ell = \card{\mathcal{T}} + 3k$ and consider the polytope $P = P_{m,\ell}$. For a matrix $x\in P$ we can identify the indices of the rows with the set~$\mathcal{T}$. Similarly, due to the definition of $\ell$ we can identify the indices of columns of $x\in P$ with $\mathcal{T}\cup X \cup Y \cup Z$; see Figure~\ref{fig:reduction} for an schematic. Hence, a binary matrix $x$ belongs to $P_{m,\ell}$ if and only if at most one entry of the vector $(x_{T,b})_{T\in \mathcal{T}}$ is 1, for each $b\in \mathcal{T}\cup X \cup Y \cup Z$.
  As before, we define $\Gamma$ as the symmetric group $\mathcal{S}(\mathcal{T}\cup X \cup Y \cup Z)= \sym{\ell}$ that acts on $\R^{m\times \ell}$ by permuting the columns.

  We define a weight vector $w$ as
  \begin{equation*}
    w_{T,b} = \begin{cases}
                -2, & \text{ if } b=T,\\
                \phantom{-}1,  & \text{ if } b\in T,\\
                \phantom{-}0, & \text{ otherwise},
              \end{cases}
            \end{equation*}
            for $T \in \mathcal{T}$ and $b \in \mathcal{T}\cup X \cup Y \cup Z$. 

            Let us consider an arbitrary order of the collection $\mathcal{T}$ given by $T_1$, $T_2$, \dots, $T_m$. We consider the SST cuts where the $i$th leader is $(T_i,T_i)$. Hence, for a given leader $(T_i,T_i)$ its corresponding SST cuts are
            \begin{equation}
              \label{eq:SSTReduction}
              x_{T_i,b}  \le x_{T_i,T_i} \text{ where } b=T_j \text{ for some } j>i \text{, or } b \in X\cup Y \cup Z.
            \end{equation}
            This set of inequalities defines the set~$S(L)$ of SST cuts.
            The proof will be completed by showing the following result.

            \begin{claim} There exists $x\in P(G_{m,\ell},S(L))$ with $\sprod{w}{x} \ge k$ if an only if the initial instance admits a 3D-matching.
            \end{claim}

            To show the claim, let us first assume that $\mathcal{T}'\subseteq \mathcal{T}$ is a 3D-matching. Consider the matrix $x$ defined as follows. For any row $T=\{T_X,T_Y,T_Z\} \in \mathcal{T}'$ of the 3D-matching, we add four 1-entries, one per column indexed by $T$, $T_X$, $T_Y$ and $T_Z$. The rest of the entries are defined as zero. More precisely, we define
            \[
              x_{T,b} =
              \begin{cases}
                1, & \text{ if } T\in \mathcal{T}' \text{ and } b=T,\\
                1,  & \text{ if } T\in \mathcal{T}' \text{ and } b\in T,\\
                0, & \text{ otherwise.}
              \end{cases}
            \]

            First, observe that $\sprod{w}{x} = k$, as any row of $x$ indexed by $T\in \mathcal{T}'$ contributes exactly $-2+1+1+1=1$ to the total sum of $\sprod{w}{x}=\card{\mathcal{T}'}=k$. Clearly $x$ is integral. Let us now argue that $x\in P(G_{m,\ell}, S(L))$. Indeed, for any column indexed by~$T$, $x$ has exactly one 1-entry if $T\in \mathcal{T}'$, and no 1-entries otherwise. Consider now a column indexed by $b\in X \cup Y\cup Z$. Then, by construction of $x$, the entries of $x$ that equal 1 within column $b$ are $x_{T,b}$ where $b\in T$ and $T\in \mathcal{T}'$. As $\mathcal{T}'$ is a 3D-matching, column $b$ must have exactly one 1-entry. We hence conclude that $x\in P_{m,\ell}$. Finally, we note that the SST cuts \eqref{eq:SSTReduction} are satisfied. Indeed, we either have that row $T$ of~$x$ is identical to 0, if $T\not\in\mathcal{T}'$, or $x_{T,T}=1$ if $T\in\mathcal{T}'$. In either case the SST-cuts are satisfied by definition.

            Conversely, let $x$ be a binary matrix in $P(G_{m,\ell},S(L))$ with $\sprod{w}{x} \ge k$. Moreover, among all possible matrices $x$, choose one with a minimum number of 1-entries. Let us construct a 3D-matching. Indeed, we define the set $\mathcal{T}'\subseteq \mathcal{T}$ of all $T$ such that $x_{T,T}=1$. The SST-cuts imply that any row $T\not\in \mathcal{T}'$ does not contribute anything to the objective function, since $\sum_{b\in \mathcal{T}\cup X \cup Y \cup Z}w_{T,b}\,x_{T,b}=0$. Furthermore, each row of $w$ sums to $-2+1+1+1=1$, and hence, for each $T\in \mathcal{T}'$ it holds that $\sum_{b\in \mathcal{T}\cup X \cup Y \cup Z}w_{T,b}\,x_{T,b}\le 1$ since $x$ is binary and satisfies the SST cuts. Moreover, if  $\sum_{b\in \mathcal{T}\cup X \cup Y \cup Z}w_{T,b}\,x_{T,b}\le 0$, we could change the complete row $T$ to be zero. This would yield another vector $x\in P(G_{m,\ell}, S(L))$ with $\sprod{w}{x} \ge k$, contradicting the minimality of $x$. Hence we conclude that $\sum_{b\in \mathcal{T}\cup X \cup Y \cup Z}w_{T,b}x_{T,b}=1$ for each $T\in \mathcal{T}'$. In other words, for any $T\in \mathcal{T}$ and $b\in X\cup Y \cup Z $, we obtain that $x_{T,b}=1$ if and only if $T\in \mathcal{T}'$ and $b\in T$.
            In particular, this implies that $\card{\mathcal{T}'}= \sprod{w}{x} \ge k$.

            Let us use this to show that $\mathcal{T}'$ is a 3D-matching. Notice that if $T$, $T'\in \mathcal{T}'$, then $T$ and $T'$ are disjoint. Indeed, if $b\in T\cap T'$, then $x_{T,b}=x_{T',b}=1$, which would violate the inequalities of $P_{m,\ell}$. Hence, $\mathcal{T}$ contains at least $k$ sets, which are pair-wise disjoint. This already implies that $\mathcal{T}'$ is a 3D-matching, as $\card{\bigcup_{T\in \mathcal{T}'} T}\ge 3k$, and hence $\mathcal{T}'$ must cover all $X\cup Y\cup Z$ with pair-wise disjoint sets. The theorem follows.
\end{proof}

\section{Strengthened SST Cuts}
\label{sec:cuts}

The previous results show that there might be a substantial interplay between the polyhedral structure of an optimization problem and~$\SSTcone(S)$.
While this interplay might drastically change the
complexity of finding integer solutions, it also comes with the potential of identifying stronger SHIs.
In the remainder of the article, we investigate the latter aspect for the stable set problem.

Since the so-called \emph{edge formulation} used in the definition of $P(G)$ is
known to be rather weak, many additional cutting planes have been derived
in the literature
to strengthen the formulation, see, e.g.,
\cite{Borndoerfer1998,Padberg1973,NemhauserTrotter1975,Tro75}.
One particular type of a strong inequality is the \emph{clique inequality}.
Clearly, any stable set contains at most one node from any clique~$C$ which leads to the clique inequality $\sum_{v \in C} x_v \leq 1$. Moreover, if~$C$ is an inclusionwise maximal clique, the corresponding clique inequality defines a facet of~$P(G)$~\cite{Padberg1973}.
We denote by~$\cC = \cC(G)$ the set of maximal cliques of~$G$.
It is well known that if $G$ is a perfect graph, $P(G)$ is completely described by non-negativity inequalities $x_v\geq 0$ for all $v\in V$ as well as clique inequalities for all $C\in \cC$~\cite{GroetschelLovaszSchrijver1984}.
We leverage this insight in order to examine the polyhedral structure of $P(G,S)$ and therefore restrict ourselves to considering perfect graphs $G$.

Instead of using the edge formulation of $P(G)$ alone, we base our approach for the analysis of $P(G,S)$ on a different formulation.
This formulation was first proposed by Johnson and Padberg~\cite{Johnson1982} to generalize the stable set problem for so-called bidirected graphs representing degree-two inequalities.
These graphs potentially have three different types of edges which emerge from equipping every two end nodes of an edge with additional signs.
We formally introduce bidirected graphs following the notation of Tamura in~\cite{Tamura1997}:
\begin{definition}
A \emph{bidirected graph}~$G=(V,E)$ consists of a set of nodes~$V$ and a set of edges~$E$ that are described by their two incident nodes and additional signs at each end.
This means there can be~$(+,+)$-,~$(+,-)$- and~$(-,-)$-edges in~$E$. Every edge type corresponds uniquely to an inequality type in the following way:
\begin{equation*}
\begin{aligned}
x_i+x_j&\leq 1 & \text{ corresponds to } (+,+)\text{-edge } \{i,j\},\\
x_i&\leq x_j & \text{ corresponds to } (+,-)\text{-edge } \{i,j\},\\
x_i+x_j&\geq 1 & \text{ corresponds to } (-,-)\text{-edge } \{i,j\}.
\end{aligned}
\end{equation*}
The correspondence means that every inequality system consisting of inequalities of the above three forms is uniquely represented by a bidirected graph and every bidirected graph uniquely encodes such an inequality system.
\end{definition}
It is immediate to make the observation that while $(+,+)$-edges correspond to the usual edge inequalities in the stable set problem, the $(+,-)$-edges correspond to SST cuts which can be viewed as directed edges pointing from the follower node to the leader node.
This way, an undirected graph with additional SST cuts corresponds to a so-called \emph{mixed graph} or \emph{pure bidirected graph}, which are bidirected graphs without $(-,-)$-edges.
Therefore, we can interpret $P(G,S)$ as the polyhedron corresponding to the mixed graph we obtain by adding the SST cuts in $S$ to $G$ as directed edges.
Since we do not work in the full generality of bidirected graphs, in the following we adapt results on bidirected graphs for our mixed graphs emerging from $G$ and $S$.

The aim of this section is to investigate the impact of SST cuts on~$P(G)$ for
subclasses of perfect graphs: bipartite graphs and trivially perfect graphs.
Before deriving our formulations for these classes in Section~\ref{sec:LinearDescription}, we first start by defining some rules in Section~\ref{sec:presolving}
that allow us to fix variables based on SST cuts.
These rules will be useful for our study of bipartite and trivially perfect graphs and can also be used as presolving techniques in computational experiments (see Section~\ref{sec:ComputationalResults}).

\subsection{Presolving Reductions}
\label{sec:presolving}

A natural question is whether the graph~$G$ can be manipulated into a
graph~$G' = (V',E')$ in such a way that~$P(G,S)$ is affinely equivalent
to~$P(G')$. 
Tamura provides in~\cite{Tamura1997} a first hint towards an answer for this question in introducing a transformation of~$G$ whose maximum stable set value agrees with the optimal solution value of~$P(G,S)$.
To get there, Tamura showed how the generalized stable set problem for a bidirected graph can be reduced to a stable set problem in polynomial time with respect to the number of nodes of~$G$.
To obtain the corresponding undirected graph, one has to make sure that the bidirected graph is closed.
A mixed graph consisting of~$G$ and the SST cuts interpreted as directed edges is called \emph{closed} if (i) every edge implied by inequalities corresponding to edges in~$G$ and~$S$ has to exist in~$G$ or~$S$, respectively, and (ii) the graph does not have any self loops or parallel edges.
 In this section, we derive a graph~$G'$ that fulfills this task in incorporating some implications
of SST cuts by removing nodes and adding edges.
These operations can be interpreted as a presolving step.
We then show in Proposition~\ref{prop:additionOperation} that optimal solutions in~$P(G')$ correspond to optima in~$P(G,S)$.
In Lemma~\ref{lm:deletionOperation}, Lemma~\ref{lm:secondDelAndAddOperation}, and Proposition~\ref{prop:additionOperation} below,~$G = (V,E)$ denotes an undirected graph and~$c \in \Z^V$ is a weight vector. Additionally,~$\group$ is an automorphism group of~$G$ that respects the node weights~$c$,~$L$ is a sequence of leaders derived from~$\group$, and~$S = S(L)$ denotes the corresponding SST~cuts.

\begin{lemma}\label{lm:deletionOperation}
  Define $V' = V \setminus \{f \in V \st \{\ell, f\} \in
    E \text{ for some } (\ell, f) \in S\}$ and its induced subgraph $G' = (V', E[V'])$.
  Then, for all~$x \in P(G,S)$ and~$v \in V\setminus V'$, we have~$x_v = 0$.
\end{lemma}

\begin{proof}
  Let~$(\ell, f)$ be a leader-follower pair.  If~$x_{f} = 1$, the SST cuts
  imply~$x_{\ell} = 1$ as well. Since at most one of them is contained in a
  stable set if~$\{\ell, f\} \in E$, $x_{f}$ can be fixed to~0, which is
  captured by $G'$.
\end{proof}

This means that the stable set problem for~$G$ and~$G'$ are equivalent, i.e., one can remove followers that are adjacent to their leaders from~$G$.
We call this \emph{deletion operation} for a follower with an adjacent leader. The \emph{deletion presolving} deletes all nodes from $G$ that can be deleted by a deletion operation.
Note that this presolving step ensures that $G'$ together with the remaining SST cuts viewed as directed edges is a simple bidirected graph.

Since deletion presolving does not incorporate implications of SST
cuts~$(\ell, f)$ if~$\ell$ and~$f$ are not adjacent,
we modify the graph~$G$ further. The \emph{addition operation}
adds all edges~$\{v, f\}$ for every neighbor~$v$ of a leader~$\ell$ with follower~$f$ to~$E$. Then
setting~$x_f = 1$ forces~$x_v = 0$ for all neighbors~$v$ of~$\ell$. Analogously to the deletion presolving, the \emph{addition presolving} adds all edges that can be added to $G$ by applying an addition operation.

This presolving step exactly corresponds to adding all the necessary edges to $G'$ such that $G'$ with additional directed edges corresponding to the SST cuts is transitive.
This means that if there is a (directed) edge whose existence is implied by inequalities corresponding to edges in $G'$ and SST cuts, this (directed) edge is already contained in $G'$ (resp. corresponds to an SST cut).
This specifically entails that applying addition presolving several times has the same effect on $G'$ as applying the presolving step only once.
Moreover, one can observe that applying the addition presolving after the deletion presolving cannot cause the considered graph to loose its simpleness, 
i.e., the addition cannot add an edge $\{u,v\}$ to $G'$ while $(u,v)\in S$. Both observations are a consequence of the directed edges arising from symmetry and are generally not true for general mixed graphs.

\begin{lemma}
\label{lem:implied_edge}
Let $u\in O_{\ell_j}$ and $i<j$. Then $\{\ell_i,\ell_j\}\in E$ if and only if $\{\ell_i,u\}\in E$.
\end{lemma}

\begin{proof}
There exists $\gamma$ in the symmetry group of $G$ that stabilizes $\ell_i$ and maps $\ell_j$ to $u$.
Since $\gamma$ is a graph automorphism of $G$, it sends edges of $G$ to edges of $G$ and does the same for non-edges.
In particular, $\{\ell_i,\ell_j\}\in E$ if and only if $\gamma(\{\ell_i,\ell_j\})\define \{\gamma (\ell_i),\gamma (\ell_j)\}=\{\ell_i,u\}\in E$.
\end{proof}

\begin{lemma}
  \label{lm:secondDelAndAddOperation}
  Let $G'=(V',E')$ arise from $G$ by applying the deletion presolving first and afterwards the addition presolving for the SST cuts $S$. 
  Then again applying the deletion presolving on $G'$ for the same SST cuts $S$ does not cause the deletion of any node of $G'$.
  Further, again applying the addition presolving on $G'$ for the same SST cuts $S$ does not cause the addition of an edge to $G'$.
\end{lemma}

\begin{proof}
Assume that the second application of the deletion presolving deletes a node $u\in V'$.
This means that the addition presolving added an edge $\{\ell_j,u\}\in E[V']$ while $(\ell_j,u)\in S$.
Since the edge $\{\ell_j,u\}$ was added by an addition operation, we have $u\in O_{\ell_i}$ and $\{\ell_i,\ell_j\}\in E$ or $\ell_j\in O_{\ell_i}$ and $\{\ell_i,u\}\in E$.
Since $(\ell_i,\ell_j),(\ell_j,u)\in S$ implies that $(\ell_i,u)\in S$, the second case cannot occur after the first application of the deletion presolving.
Therefore, we can assume the first case.
If $j<i$, applying Lemma~\ref{lem:implied_edge} gives us that~$\{\ell_j,u\} \in E$ contradicting that the edge was added by an addition operation.
For $i<j$, again applying Lemma~\ref{lem:implied_edge}, we obtain $\{\ell_i,u\}\in E$ contradicting that $u\in V'$.

Now assume that the edge $\{u,v\}$ is added to $G'$ during a second application of the addition presolving.
This means that there exists an edge $\{\ell_i,v\}$ that was added during the first application of the addition presolving with $u\in O_{\ell_i}$.
Therefore, there exists $\{\ell_i,\ell_j\}\in E$ with $v\in O_{\ell_j}$ or $\{v,\ell_j\}\in E$ with $\ell_i\in O_{\ell_j}$.
In the first case, Lemma~\ref{lem:implied_edge} implies that $\{\ell_i,v\}\in E$ or $\{\ell_j,u\}\in E$ both implying that $\{u,v\}$ is added to $G'$ during the first application of the addition presolving because $u\in O_{\ell_i}$ as well as $v\in O_{\ell_j}$.
In the second case, because of the transitivity of SST cuts, we also have $(\ell_j,u)\in S$ and consequently, $\{u,v\}$ would have been added during the first application of the addition presolving.
\end{proof}
\begin{proposition}
  \label{prop:additionOperation}
  Let~$G' = (V', E')$ arise from~$G$ by applying the deletion presolving and then the addition
  presolving for a set of SST cuts.
  Suppose~$c_v \neq 0$ for all~$v \in V$.
  Then, every weight maximal stable set in~$G'$ is weight maximal in~$G$
  and satisfies all SST cuts.
\end{proposition}

\begin{proof}
  The implications of SST cuts are that setting~$x_f = 1$ for a follower~$f$
  implies that~$x_\ell = 1$ for the corresponding leader~$\ell$. We distinguish
  two cases: If $f$ and $\ell$ are adjacent, then the deletion operation
  covers the corresponding implication. Otherwise, the addition operation
  covers this case: if~$x_f = 1$, then the edges introduced by the addition
  operation cause~$x_v = 0$ for all neighbors~$v$ of~$\ell$.  Hence,
  if~$c_\ell > 0$, $x_\ell = 1$ in an optimal solution if~$x_f = 1$.
  Moreover, if~$c_\ell < 0$, then $x_f$ is not set to 1 in an optimal
  solution, since~$c_f = c_\ell < 0$ because~$\ell$ and~$f$ are symmetric.
  Finally, note that setting~$x_v = 1$ for some neighbor~$v$ of~$\ell$
  implies~$x_f = 0$ and~$x_\ell = 0$.  Thus, exactly the implications of SST
  cuts are incorporated by the deletion and addition operation, which keeps
  at least one optimal solution intact.
\end{proof}

The deletion and addition presolving steps can be used as a
symmetry-based presolving routine, which we call \emph{SST presolving}.
We investigate whether the graph manipulations of SST presolving indeed simplify the problem by conducting numerical experiments in Section~\ref{sec:ComputationalResults}.
We close this section by remarking that SST presolving does not incorporate
all implications of SST cuts into the underlying graph.
Indeed, Proposition~\ref{prop:additionOperation} only implies that an
optimal solution will adhere to the SST cuts; suboptimal solutions may
violate the inequalities as we illustrate by the following example.
Note that this is consistent with Theorems~\ref{thm:opt_CS} and~\ref{thm:NPhardness}.
Indeed, finding an optimal solution under SST cuts for a weight vector that is \emph{invariant} under $\group$ can be handled in polynomial time, as long as the original stable set problem is tractable.
For \emph{arbitrary} objectives, however, the optimization problem is NP-hard as the separation problem for~$P(G,S)$ is NP-hard in general (Theorem~\ref{thm:NPhardness} in Section~\ref{sec:complexity}).
That is, we cannot expect to find simple graph manipulations to express~$P(G,S)$ via the stable set polytope of a graph~$G'$.
\begin{example}
  \label{ex:presolving}
  Consider the graph~$G$ in Figure~\ref{fig:presolveA} whose symmetries are
  given by rotating the nodes along the cycle defined by the graph and
  ``reflecting'' the node labels along the shown lines.
  If we select the first leader~$\ell_1$ to be node~1, the deletion
  presolving removes nodes~$2$ and~$8$ from the graph.
  Figure~\ref{fig:presolveB} shows the resulting graph.
  Its only non-trivial symmetry is the reflection along the displayed line.
  If~$\ell_2 = 3$ is selected as second leader, its orbit is~$\{3,7\}$.
  As such, the deletion presolving does not remove any more nodes from the graph.
  The addition presolving, however, adds the edge~$\{4,7\}$, corresponding
  to the SST cut~$x_7 \leq x_3$.
  The resulting graph~$G'$ is shown in Figure~\ref{fig:presolveC}.
  Note that the single node~$7$ forms a stable set in~$G'$, but its
  characteristic vector does not adhere to the SST cut~$x_7 \leq x_3$.
  Hence, SST presolving does not model all implications of SST cuts.
\end{example}

\begin{remark}
  SST presolving can also introduce new symmetries.
  A symmetry of the graph in Figure~\ref{fig:presolveC} is to
  exchange nodes~5 and~7 while keeping all other nodes fixed.
  Obviously, this is not a symmetry of the original graph in
  Figure~\ref{fig:presolveA}. This phenomenon will be exploited in our computational experiments in Section~\ref{sec:ComputationalResults}.
\end{remark}

\begin{figure}[tb]
  \begin{minipage}[t]{0.31\linewidth}
    \centering
    \begin{tikzpicture}
      \tikzstyle{s} += [circle,draw=black,inner sep=0pt,minimum size=2mm];
      \node (1) at (0,0) [s,label=left:1] {};
      \node (2) at (1,-1) [s,label=below:2] {};
      \node (3) at (2,-1) [s,label=below:3] {};
      \node (4) at (3,0) [s,label=right:4] {};
      \node (5) at (3,1) [s,label=right:5] {};
      \node (6) at (2,2) [s,label=above:6] {};
      \node (7) at (1,2) [s,label=above:7] {};
      \node (8) at (0,1) [s,label=left:8] {};

      \draw[-] (1) -- (2) -- (3) -- (4) -- (5) -- (6) -- (7) -- (8) -- (1);
      \draw[-,red] (-.25,.5) -- (3.25,.5);
      \draw[-,red] (1.5,-1.25) -- (1.5,2.25);
      \draw[-,red] (0,-1) -- (3,2);
      \draw[-,red] (3,-1) -- (0,2);
      \draw[-,red] (1) -- (5);
      \draw[-,red] (2) -- (6);
      \draw[-,red] (3) -- (7);
      \draw[-,red] (4) -- (8);
    \end{tikzpicture}
      \captionof{figure}{Original graph and its reflection symmetries
        (along red lines).}
      \label{fig:presolveA}
  \end{minipage}
  \hfill
  \begin{minipage}[t]{0.31\linewidth}
    \centering
    \begin{tikzpicture}
      \tikzstyle{s} += [circle,draw=black,inner sep=0pt,minimum size=2mm];
      \node (1) at (0,0) [s,label=left:1] {};
      \node (3) at (2,-1) [s,label=below:3] {};
      \node (4) at (3,0) [s,label=right:4] {};
      \node (5) at (3,1) [s,label=right:5] {};
      \node (6) at (2,2) [s,label=above:6] {};
      \node (7) at (1,2) [s,label=above:7] {};
      \draw[-,red] (1) -- (5);

      \draw[-] (3) -- (4) -- (5) -- (6) -- (7);
    \end{tikzpicture}
      \captionof{figure}{Graph after applying addition/deletion for leader $\ell_1 = 1$.}
      \label{fig:presolveB}
  \end{minipage}
  \hfill
  \begin{minipage}[t]{0.31\linewidth}
    \centering
    \begin{tikzpicture}
      \tikzstyle{s} += [circle,draw=black,inner sep=0pt,minimum size=2mm];
      \node (1) at (0,0) [s,label=left:1] {};
      \node (3) at (2,-1) [s,label=below:3] {};
      \node (4) at (3,0) [s,label=right:4] {};
      \node (5) at (3,1) [s,label=right:5] {};
      \node (6) at (2,2) [s,label=above:6] {};
      \node (7) at (1,2) [s,label=above:7] {};

      \draw[-] (3) -- (4) -- (5) -- (6) -- (7);
      \draw[-] (4) -- (7);
    \end{tikzpicture}
      \captionof{figure}{Graph after applying addition/deletion for leader $\ell_1$ and then $\ell_2 = 3$.}
      \label{fig:presolveC}
  \end{minipage}
\end{figure}

\subsection{Complete Linear Descriptions for Special Perfect Graphs}\label{sec:LinearDescription}

In the rest of this section we study linear descriptions of $P(G,S)$ for bipartite graphs and trivially perfect graphs.
We focus on these two graph classes because of their well known facet descriptions which provide a strong starting point for the investigation of the impact of SST cuts on the polyhedral structure of the stable set problem.

We show for bipartite graphs~$G$ that the separation problem for~$P(G,S)$ is tractable and
give an explicit compact size extended formulation and linear outer description.
Unlike the case of trivially perfect graphs, for bipartite graphs this can be done regardless of the choice of leaders.
In the case of trivially perfect graphs~$G$, we show NP-hardness of the separation problem of $P(G,S)$ which motivates the derivation of leader sequences ensuring tractability.

\subsubsection{Bipartite Graphs}
\label{sec:bipartite}

A graph~$G = (V, E)$ is \emph{bipartite} if~$V$ admits a partition~$R \cup
B$ such that~$E \subseteq \{ \{u,v\} \st u \in R,\; v \in B\}$.
We refer to the two sets as the \emph{red} and \emph{blue} color class of the
bipartition, respectively.

To find a complete linear description of~$P(G,S)$, we introduce an auxiliary
graph $G^\star = (V^\star, E^\star)$.  The node set is~$V^\star = V \cup \bar{V}$, where
$\bar{V} \define \{v_1, \dots, v_n\}$ such that vertex $v_{\ell}$ corresponds to the leader~${\ell \in L}$;
 the edge set is~$E^\star = E \cup \bar{E}$, where
$\bar{E} \define \{ \{v_\ell, f\} \st f \in O_\ell,\; \ell \in L \}$. Note that $\ell \in
O_\ell$ and therefore~${\{v_\ell,\ell\} \in \bar{E}}$.
The
graph~$G^\star$ can be used to define an extended formulation of~$P(G,S)$, even
for non-bipartite graphs, as shown in the following proposition.

\begin{proposition}\label{prop:extendedFormulation}
  Given an undirected graph~$G = (V,E)$ it holds that
  \[
    P(G,S) = \{ x \in \R^V \st \exists y \in \R^{\bar{V}} \text{ with } (x,y) \in
    P(G^\star) \text{ and } x_\ell + y_{v_\ell} = 1 \text{ for } \ell
    \in L \}.
  \]
\end{proposition}

\begin{proof}
  Let  
  \[
    Q = \{ x \in \R^V \st \exists y \in \R^{\bar{V}} \text{ with } (x,y) \in
    P(G^\star) \text{ and } x_\ell + y_{v_\ell} = 1 \text{ for } \ell
    \in L \}.
  \]
  First, we show~$P(G,S)\subseteq Q$.
  Since~$P(G, S)$ is a polytope, it is sufficient to show that each vertex
  of~$P(G,S)$ is contained in~$Q$.
  Let~$x \in P(G, S) \cap \B{V}$ and define~$y \in \B{\bar{V}}$ via~$y_{v_\ell} = 1 -
  x_{\ell}$ for each~$\ell \in L$.
  It is enough to show~$(x,y) \in P(G^\star)$.
  To this end, let~$\ell \in [n]$.
  If~$y_{v_\ell} = 0$, then all edge inequalities of~$P(G^\star)$ hold trivially.
  Otherwise, if~$y_{v_\ell} = 1$, then~$x_\ell = 0$ and the SST cuts force~$x_v = 0$ for all~$v \in O_\ell$.
  Therefore, $x_v + y_{v_\ell} \leq 1$ holds for all neighbors~$v$ of~$v_\ell$ in~$G^\star$.
  Since~$G$ is an induced subgraph of~$G^\star$ and~$x \in P(G)$, also the remaining edge inequalities hold.
  We thus conclude~$(x,y) \in P(G^\star)$ as desired.
  
  Second, we prove~$Q \subseteq P(G,S)$.
  Since~$P(G^\star)$ is a 0/1-polytope, the assertion follows if we can show
  that every binary vector in~$Q$ is contained in~$P(G,S)$.
  Therefore, let~$(x,y) \in P(G^\star) \cap \B{V \cup \bar{V}}$ such that~$x_\ell + y_{v_\ell} = 1$ for
  all~$\ell \in L$.
  We claim that~$x$ satisfies all SST cuts~$S$.
  If this was not the case, there would exist a leader~$\ell \in L$ and
  follower~$f \in O_\ell$ such that~$x_\ell = 0$ and~$x_f = 1$.
  This implies~$y_{v_\ell} = 1$ and, since~$\{v_\ell, f\} \in E^\star$,
  $y_{v_\ell} + x_f \leq 1$ is a valid inequality for~$P(G^\star)$.
  Thus, $(x,y) \notin P(G^\star)$, a contradiction.
\end{proof}

Based on this extended formulation, we can find a complete linear
description of~$P(G,S)$ by projecting a face of~$P(G^\star)$ onto~$\R^V$.
While outer descriptions of such faces are unknown for general graphs, we can derive such an outer description for bipartite graphs. To do so, we say that an
(induced) path~$P$ in~$G$ is a sequence~$u_1,\dots,u_k$ of nodes such
that~$\{u_i,u_{i+1}\} \in E$ for all~$i \in [k-1]$ and there is no edge
connecting other nodes in~$P$.
To emphasize the start and end node of a path we also call it $(u_1-u_k)$-path.
We say that a path~$P$ is \emph{even} if~$\card{V(P)}$ is even, where $V(P)$ is the set of nodes of $P$.
If~$P$ contains the additional edge~$\{u_k, u_1\}$, we call~$P$ a
\emph{cycle}.
We say that the cycle is \emph{odd} if~$\card{V(P)}$ is odd.
Before we prove the main result on bipartite graphs, we provide an
auxiliary result first.
\begin{lemma}
  \label{lem:bipartiteTwoSides}
  Let~$G$ be a connected bipartite graph with bipartition~$R \cup
  B$ and $L$ be a sequence of leaders.
  Then, there exists at most one leader~$\ell \in L$ such that $O_\ell \cap R
  \neq \emptyset$ and~$O_\ell \cap B \neq \emptyset$.
\end{lemma}
\begin{proof}
  For the sake of contradiction, suppose there exist two distinct leaders~$\ell$ and~$\ell'$ whose orbits contain red and blue nodes.
  Let~$\Gamma_\ell$ and~$\Gamma_{\ell'}$ be the groups used to compute~$O_\ell$ and~$O_{\ell'}$, respectively, and let~$\ell$ appear before~$\ell'$ in~$L$.
  Then, every permutation in~$\Gamma_{\ell'}$ needs to stabilize~$\ell$ by~\ref{alg:B}.

  Since~$O_{\ell'}$ contains red and blue nodes, there is~$\perm \in \Gamma_{\ell'}$ that maps a red node~$v$ onto a blue node.
  Because~$\perm$ preserves bipartiteness of~$G$, all red nodes in the
  connected component of~$v$ need to be mapped onto a blue node and vice
  versa.
  Thus, since~$G$ is connected, $\perm$ exchanges the red/blue label of all
  nodes.
  This in particular means that the color label of~$\ell$ is changed, which contradicts that~$\gamma$ stabilizes~$\ell$.
  For this reason, there is at most one leader whose orbit contains red and
  blue nodes.
\end{proof}
\begin{theorem}
  \label{thm:descriptionBipartite}
  Let~$G = (V,E)$ be a connected bipartite graph. Then, $P(G,S)$ is completely described by box constraints~$0 \leq x_v \leq 1$ for~$v \in V$, edge
  inequalities $x_u+x_v\le 1$ for $\{u,v\}\in E$, SST cuts $S(L)$, and, for each leader~$\ell \in L$,
  \begin{align}
    \label{eq:SSTpath}
    -x_{\ell} + \sum_{v \in V(P)} x_v &\leq \frac{\card{V(P)}}{2} - 1,
    && \text{for all even $(\ell-f)$-paths~$P$,}\; f \in O_\ell.
  \end{align}
\end{theorem}
\begin{proof}
  Note that, by construction, the graph~$G^\star$ is bipartite if, for each~$\ell \in L$, the
  orbit~$O_\ell$ exclusively consists of either red or blue nodes.
  Moreover, if there exists an orbit with both colors, there exists exactly
  one~$\ell \in L$ such that~$O_\ell$ contains nodes with both colors by Lemma~\ref{lem:bipartiteTwoSides}.
  Hence, removing the node~$v_\ell$ and all its incident edges from~$G^\star$
  results in a bipartite graph.
  Such graphs are called almost bipartite. They are t-perfect,
  that is, their stable set polytope is completely described by box constraints, edge
  inequalities, and odd cycle inequalities
  \[
    \sum_{v \in V(C)} x_v \leq \frac{\card{V(C)} - 1}{2}
  \]
  for every odd cycle~$C$ in the graph, cf.~\cite{FonluptUhry1982}.
  To derive the desired outer description of~$P(G,S)$, we project~$P(G^\star)$
  onto~$P(G,S)$ using the equations~$x_{\ell} = 1 - y_{v_\ell}$ for all~$\ell \in
  L$.
  We discuss this projection for the three different types of inequalities
  for~$P(G^\star)$ in turn.

  \emph{Odd Cycle Inequalities:}
  Since~$G$ is bipartite, $C$ is an odd cycle in~$G^\star$ if and only if it
  contains a node~$v_\ell$, $\ell \in L$, such that~$O_\ell$ contains red
  and blue nodes.
  By construction of~$G^\star$, each such cycle contains~$\ell$ as well as a
  node~$f \in O_\ell \setminus \{\ell\}$, and no such cycle can
  contain~$v_{\ell'}$ for~$\ell' \in L \setminus \{\ell\}$.
  For this reason, each odd cycle can be decomposed into an even path~$P$ with
  endpoints~$\ell$ and~$f$ all of whose nodes are contained in~$V$, as well
  as the path~$\ell - v_\ell - f$.
  If we eliminate the~$y$-variables using the
  equation~$y_{v_\ell} = 1 - x_\ell$ for all~$\ell \in L$, the odd cycle
  inequalities~$y_{v_\ell} + \sum_{v \in V(P)} x_v \leq \frac{\card{V(P)}}{2}$
  are transformed into even path inequalities
  \[
    -x_\ell + \sum_{v \in V(P)} x_v \leq \frac{\card{V(P)}}{2} - 1.
  \]

  \emph{Edge Inequalities:}
  Let~$x_u + x_v \leq 1$ be an edge inequality for~$P(G^\star)$.
  If one of the nodes, say~$v$, is~$v_\ell$ for a leader~$\ell \in L$, the
  edge inequality transforms into~$x_u + 1 - x_\ell \leq 1$, which is the same
  as~$x_u - x_\ell \leq 0$.
  If~$u = \ell$, then the inequality is redundant.
  Otherwise, by definition of~$\bar{E}$, $u \in O_\ell \setminus \{\ell\}$,
  i.e., the transformed edge inequality is an SST cut for~$P(G,S)$.
  Finally, if neither~$u$ nor~$v$ is~$v_\ell$ for some~$\ell \in L$, then
  the edge inequality corresponds to an edge inequality for~$P(G,S)$.

  \emph{Box Constraints:}
  All box constraints transform to box constraints for~$P(G,S)$,
  except for~$0 \leq y_{v_\ell} \leq 1$ for~$\ell \in L$.
  Using~$y_{v_\ell} = 1 - x_\ell$ yields the box constraint for~$x_\ell$.
\end{proof}
\begin{corollary}\label{cor:PGSBipartiteWithoutCycles}
  Let~$G = (V,E)$ be a connected bipartite graph, $L$ be a sequence of leaders, and $S = S(L)$.
  Then, $P(G,S)$ is completely described by box constraints, edge
  inequalities, and SST cuts if and only if every orbit~$O_\ell$, $\ell \in
  L$, contains only nodes from one color class.
\end{corollary}
\begin{proof}
  On the one hand, if every orbit~$O_\ell$ contains only nodes from one
  color class, then every path from~$\ell$ to a follower is necessarily odd
  as both endpoints of the path are contained in the same color class.
  Hence, no inequality of type~\eqref{eq:SSTpath} exists.
  On the other hand, if there exists an even leader-follower path~$P$, the
  corresponding inequality is needed in an outer description of~$P(G,S)$:
  If we remove the inequality from the description of~$P$ given in
  Theorem~\ref{thm:descriptionBipartite}, consider the solution~$x$ that
  is~$\frac{1}{2}$ on~$P$ and~0 otherwise.
  Then, $x$ has a unique pre-image~$(x,y)$ in the extended formulation, whose
  support defines an odd cycle~$C$ in~$G^\star$ and such that each entry in the
  support is~$\frac{1}{2}$.
  For~$(x,y)$, the inequality for~$C$ evaluates to $\frac{\card{V(C)}}{2}$.
  It is well-known that the only odd cycle inequality that cuts off~$(x,y)$
  is the odd cycle inequality for~$C$.
  Moreover, no edge inequality or box constraint is violated by~$(x,y)$.
  Hence, the only inequality in the description of
  Theorem~\ref{thm:descriptionBipartite} that separates~$x$ is the path
  inequality for~$P$, which shows that it is needed in an outer description.
\end{proof}

\begin{remark}
  If the path of Inequality~\eqref{eq:SSTpath} consists of two nodes, i.e., it
  is~$\ell-f$ for some leader/follower pair, then the inequality reduces
  to~$x_f \leq 0$.
  That is, it fixes the follower to~0, which corresponds to the deletion operation.
\end{remark}
We also remark that \emph{SST path cuts}~\eqref{eq:SSTpath} can be derived for
non-bipartite graphs, too, because they correspond to odd cycle inequalities
for~$P(G^\star)$. In particular, as odd cycle inequalities can be separated in
polynomial time~\cite{GrotschelLovaszSchrijver1993}, SST path cuts are
polynomial time separable.

\subsubsection{Trivially Perfect Graphs}
\label{sec:triviallyPerfect}
Trivially perfect graphs (\tpgraphs) are perfect graphs that form a subclass of interval graphs.
An undirected graph~$G = (V,E)$ is an \emph{interval graph} if for each~$v
\in V$ there is an interval~$I_v\subseteq \R$ such that, for all
distinct~$u$, $v \in V$, we have~$\{u,v\} \in E$ if and only if~$I_u \cap I_v
\neq \emptyset$.
An interval graph is \emph{trivially perfect (TP)} if its
interval collection $(I_v)_{v\in V}$ can be chosen to be laminar. That is, for all $u$, $v\in V$ if $I_u\cap I_v\neq \emptyset$ then $I_v\subseteq I_u$ or $I_u\subseteq I_v$.

An important observation is that the polytope $P_{m,\ell}$ in Theorem~\ref{thm:NPhardness} corresponds to the clique formulation of a trivially
perfect graph.
\begin{corollary}\label{cor:NPhardPerfect}
  Let $G$ be a trivially perfect graph, $\group$ an automorphism group for~$G$, and $S=S(L)$ a set of SST cuts for some set of leaders $L$. Then the separation problem for $P(G,S)$ is NP-hard.
\end{corollary}
\begin{proof}
	For given natural numbers $m$ and $\ell$, let $G=(V,E)$ be a graph consisting of~$\ell$ node-disjoint cliques, with $m$ nodes each. That is, we have that $V=V_1\dot{\cup} \ldots \dot{\cup} V_\ell$ and $E= E_1\dot{\cup} \ldots \dot{\cup} E_\ell$, where~$(V_i,E_i)$ is a complete graph of $m$ nodes for each $i\in\{1,\ldots,\ell\}$. It is easy to see that $G$ is trivially perfect: for each node $v\in V_i$ we define an interval $I_v \define I_i$ for each $i$, where $I_i\cap I_j=\emptyset$ for $i\neq j$. Then we observe that the clique formulation of $P(G)$ is exactly $P_{m,\ell}$ as defined in Section~\ref{sec:complexity}. The result follows from Theorem~\ref{thm:NPhardness}.
\end{proof}
Since the separation problem of $P(G,S)$ is NP-hard, even if $G$ is
trivially perfect, there is a
complex interaction between~$S$ and the underlying graph's structure.
The goal of this section is to better understand this interaction and to derive a class
of SST cuts, given by a carefully selected sequence of leaders, such that we can fully describe~$P(G,S)$ for a \tpgraph~$G$. In particular, we give a specific construction of SST cuts that avoids the NP-hardness of the separation problem for $P(G,S)$ for any trivially perfect graph.
To this end, it will be useful to consider a strengthening of SST cuts to so-called SST clique cuts.
These are a special form of biclique inequalities which were introduced in~\cite{Johnson1982}.

\begin{definition}
Let~$G=(V,E)$ be a bidirected graph.
For a subset~$W\subseteq V$ we denote by~$W_+\subseteq W$ (resp.~$W_-\subseteq W$) the subset of nodes in~$W$ incident to edges that have a plus end (resp. minus end) at this node.
A completely connected subgraph~$C=(W,E_W)$ of~$G$ is called \emph{biclique} if~$W_+\cap W_-=\emptyset$ holds.
\end{definition}

Johnson and Padberg~\cite[Prop.~6]{Johnson1982} showed that all binary solutions of the inequality system corresponding to a bidirected graph satisfy the so-called \emph{biclique inequalities}
\begin{equation*}
\sum_{i\in W_-}(1-x_i)+\sum_{i\in W_+}x_i\leq 1.
\end{equation*}
For mixed graphs, one always has~$\card{S_-}\leq 1$ for a biclique. Moreover, in the particular setting where mixed graphs are derived from SST cuts, we obtain the following corollary.
\begin{corollary}
  Let~$G$ be an undirected graph.
  If~$\ell$ is a leader for the symmetry group of the stable set problem over $G$ with orbit $O_{\ell}$, the
  \emph{SST clique cut}
  \begin{equation}\label{eq:SSTclique}
    -x_\ell + \sum_{f \in C} x_f \leq 0
  \end{equation}
  is an SHI for each clique~$C \subseteq O_{\ell}$.
\end{corollary}

Note that in~\cite{Johnson1982} the definition of a biclique differs from our definition in also demanding maximality under the stated conditions.
Anyhow, clearly the validity of SST clique cuts remains if the associated clique is not maximal in~$O_{\ell}$.

Since a single follower defines a clique in~$G$, SST clique cuts can be seen as a generalization of SST cuts.
Consequently, extending SST cuts to SST clique cuts exploits the potential of SST cuts within the broader framework of biclique inequalities and therefore makes this type of cut as strong as possible.
Among these cuts, Johnson and Padberg characterized the ones inducing facets of the according polyhedron for closed bidirected graphs (see Theorem 10 in~\cite{Johnson1982}).
They call them \emph{strong biclique inequalities} which are biclique inequalities for bicliques~$C=(W,E_W)$ with the following additional property: There may not exist a node~$i\in V\setminus W$ such that for all~$j\in W_+$ there exists a~$(\cdot,+)$-edge~$\{i,j\}$ and for all~$j\in W_-$ there exists a~$(\cdot,-)$-edge~$\{i,j\}$.
We obtain the facet characterization in our special setting as a corollary.
\begin{corollary}
  \label{lem:SSTcliqueFacet}
  Let~$G = (V,E)$ be an undirected graph presolved by SST presolving for SST cuts $S$ 
  and let~$\ell \in V$ be a leader with orbit~$O_\ell$.
  Moreover, let $C\subseteq O_\ell\setminus \{\ell\}$ be a clique.
  Then the SST clique cut ~$\sum_{i \in C} x_i \leq x_\ell$ defines a facet of~$P(G,S)$ if and only if there exists no~$f\in O_\ell \setminus (\{\ell\}\cup C)$ such that for every $v\in C$ either $v\in O_f$ or $\{f,v\}\in E$ where $O_f$ is the orbit of $f$.
\end{corollary}
Since clique cuts as well as SST clique cuts do both fall under the category of biclique cuts, they can be seen as a natural generalization of the clique cuts in the linear description of the stable set polyhedron for perfect graphs.
Indeed, these cuts will enable us to formulate a complete linear description of~$P(G,S)$ for certain SST cuts~$S$.
After that, we present a mechanism for deriving such SST cuts~$S$ for which the mentioned formulation of~$P(G,S)$ is applicable.
As cited in~\cite{Tamura1997}, Ikebe and Tamura characterized when $P(G,S)$ can be completely described using biclique inequalities.
This description relies heavily on a modification of the underlying bidirected graph which is replacing all edge signs by~$(+,+)$-signs.
In our setting, this means that the SST cuts~$S$ are incorporated into~$G$.
The modified graph~$\underline{G}_S$ is accordingly the undirected graph we obtain by
starting with~$\underline{G}_S = G$ and adding for all SST cuts~$(\ell ,f)$ the edge~$\{\ell ,f\}$ to~$\underline{G}_S$.
In~\cite{Ikebe1996}\footnote{Unfortunately, we do not have access to this source. Anyway, the direction of the equivalence relevant for our description of the polyhedron~$P(G,S)$ is a direct consequence of the result that a graph is biperfect if and only if it is perfect, conjectured by Johnson and Padberg in~\cite{Johnson1982} and shown by  Sewell in~\cite{Sew96}.} it is shown that for a closed bidirected graph, all strong biclique inequalities describe the integer hull of solutions of the inequality system corresponding to the graph completely if and only if the modified graph is perfect.
\begin{theorem}[Corollary of \cite{Ikebe1996}]
  \label{thm:SSTclique}
Let~$G=(V,E)$ be an undirected graph with a set of SST cuts~$S$ such that SST presolving has already been applied on~$G$ regarding~$S$.
The polyhedron~$P(G,S)$ is completely described by box constraints, clique cuts, and SST clique cuts if and only if~$\underline{G}_S$ is a perfect graph.
\end{theorem}
To make use of this theorem, we show that a carefully constructed family of SST cuts, so-called stringent SST cuts as defined next, guarantee~$\underline{G}_S$ to be a perfect graph if $G$ is a trivially perfect graph.
To do so, we use the fact that an undirected graph is a \tpgraph if and only if it does not have paths or cycles on four nodes as induced subgraphs.
This means whenever there is a path $(\{u,v\},\{v,w\},\{w,y\})$ in a \tpgraph $G$, the edge $\{u,w\}$ or the edge $\{v,y\}$ also has to be contained in~$G$.
\begin{definition}
Let~$(\ell_1, \dots, \ell_n)$ be the leaders of a family of SST cuts~$S$.
For every~$i \in [n]$, let~$O_i$ be the orbit of~$\ell_i$ from Step~\ref{alg:A} of the SST algorithm, and let~$\mathfrak{O}'_i = \mathfrak{O}_i \setminus \{\ell_1,\dots,\ell_{i-1}\}$ with~$\mathfrak{O}_i = \bigcup_{j = 1}^{i-1} O_j$.
The family $S$ is called \emph{stringent} if for all $i\in [n]$ we have~$\ell_i\in [n] \setminus \{\ell_1,\dots,\ell_{i-1}\}$ and
\begin{equation}
N(u)\cap \mathfrak{O}'_i\subseteq \{v\in V \;:\; \exists \; j<i \;\text{s.t.}\; u,v\in O_j\} \ \text{for all}\ u\in O_i,
  \label{eq:stringent}
\end{equation}
where $N(u)$ is the neighborhood of $u$ in $G$.
\end{definition}
While general SST cuts allow to select the leaders in an arbitrary
order, stringent SST cuts enforce some hierarchy in the leader selection.
This hierarchy allows us to show perfectness for $\underline{G}_S$ and consequently provide a complete description of~$P(G,S)$ for \tpgraphs~$G$:
\begin{theorem}
  \label{thm:SSTperfect}
Let~$G = (V,E)$ be a trivially perfect graph.
Let~$L$ be a sequence of leaders such that the corresponding SST cuts~$S = S(L)$ are stringent.
Assume that before constructing $\underline{G}_S$, SST presolving is applied on $G$ for $S(L)$.
Then $\underline{G}_S$ is a perfect graph.
\end{theorem}
We defer the proof of Theorem~\ref{thm:SSTperfect} to the next section, and discuss the
result first.
The description of~$P(G,S)$ from the previous theorem can be separated
efficiently, since we can show that the formulation has at most~$O(n^2)$
constraints. Moreover, none of the assumptions in Theorem~\ref{thm:SSTperfect} can be
dropped since there exist counterexamples for the respective cases as we
discuss next.

Consider the \tpgraph~$G$ from Figure~\ref{fig:stringent} and suppose that we
select~$\ell_1 = 3$ as first leader with orbit~$O_1 = \{3,4,5,10\}$.
If we do not require to compute stringent SST cuts, we are allowed to
select~$\ell_2 = 2$ as second leader with orbit~$O_2 = \{2,7,12\}$.
Indeed, for stringent SST cuts we are not allowed to select~$\ell_2 = 2$
since node $4$ is in $N(2)\cap \{4,5,10\}$.
The dashed lines in Figure~\ref{fig:stringent} demonstrate that, caused by the addition presolving, the graph $\underline{G}_S$ contains an odd hole and is therefore not perfect.
Hence, the linear description in Theorem~\ref{thm:SSTclique} results in a polyhedron that has fractional vertices.

\begin{figure}[tb]
\centering
    \begin{tikzpicture}
      \tikzstyle{s} += [circle,draw=black,inner sep=0pt,minimum size=2mm];

      \node[pentagon] (E) at (0,0) {};
      \node[white] (E1) at (E.corner 1) {};
      \node[white] (E2) at (E.corner 2) {};
      \node[white] (E3) at (E.corner 3) {};
      \node[white] (E4) at (E.corner 4) {};
      \node[white] (E5) at (E.corner 5) {};

      \node[triangleA] (A) at (E3) {};
      \node[s,fill=green,label=below:{$12$}] (A1) at (A.corner 1) {};
      \node[s,fill=orange,label=right:{$11$}] (A2) at (A.corner 2) {};
      \node[s,fill=blue,label=left:{$10$}] (A3) at (A.corner 3) {};

      \node[triangleB] (B) at (E5) {};
      \node[s,fill=green,label=left:{$7$}] (B1) at (B.corner 1) {};
      \node[s,fill=orange,label=right:{$\ell_3=9$}] (B2) at (B.corner 2) {};
      \node[s,fill=blue,label=right:{$5$}] (B3) at (B.corner 3) {};

      \node[triangleC] (C) at (E2) {};
      \node[s,fill=green,label=left:{$8$}] (C1) at (C.corner 1) {};
      \node[s,fill=orange,label=right:{$6$}] (C2) at (C.corner 2) {};
      \node[s,fill=blue,label=left:{$\ell_1=3$}] (C3) at (C.corner 3) {};

      \node[triangleD] (D) at (E1) {};
      \node[s,fill=green,label=right:{$\ell_2=2$}] (D1) at (D.corner 1) {};
      \node[s,fill=orange,label=left:{$1$}] (D2) at (D.corner 2) {};
      \node[s,fill=blue,label=below:{$4$}] (D3) at (D.corner 3) {};

      \draw[dashed] (C2)--(A3);
      \draw[dashed] (C2)--(D3);
      \draw[dashed] (D3)--(B1);
      \draw[dashed] (B1)--(A2);

      \draw[dashed, gray, opacity=0.7] (C2)--(B3);
      \draw[dashed, gray, opacity=0.7] (C1)--(A3);
      \draw[dashed, gray, opacity=0.7] (C1) to [bend right=15] (B3);
      \draw[dashed, gray, opacity=0.7] (C1) to [bend left=20] (D3);
      \draw[dashed, gray, opacity=0.7] (D3)--(A1);
      \draw[dashed, gray, opacity=0.7] (D2)--(B1);
      \draw[dashed, gray, opacity=0.7] (D2)--(A1);
      \draw[dashed, gray, opacity=0.7] (B3)--(A2);
    \end{tikzpicture}
      \captionof{figure}{Example for non-stringent SST cuts where different colors indicate different node weights. Dashed lines are the edges added by the addition presolving\protect\footnotemark .}
      \label{fig:stringent}
\end{figure}
\footnotetext{Instead of giving the nodes different node weights to make sure that the nodes in a triangle are asymmetric, one can also attach asymmetric graphs to the nodes. I.e., we can replace the colors by adding a common neighbor to the former blue and green node and adding another neighbor to the former blue node.}
Moreover, since SST clique cuts dominate SST cuts, it is necessary to
replace SST cuts by SST clique cuts.
Also, the requirement of being trivially perfect and to apply the deletion presolving
are necessary for the validity of the theorem: If~$G$ is a path of length~$6$ and the first leader is an end node of the path,~$\underline{G}_S$ will be an odd hole of size 5 with an attached triangle and therefore not perfect.
Finally, consider the graph consisting of only two adjacent nodes. Then an SST cut between those nodes can be introduced and the polytope described by the SST cut, the edge constraint and the variable bounds, has the fractional vertex~$(\frac{1}{2},\frac{1}{2})$. Thus, we need to apply deletion presolving to avoid having fractional vertices.

We conclude this section by an observation that is crucial for the proof of Theorem~\ref{thm:SSTperfect}.
\begin{observation}
  \label{obs:stringent}
For a trivially perfect graph~$G=(V,E)$ and stringent SST cuts~$S$ for which SST presolving has been applied, one has $\bigcup_{u\in O_i}N(u)\cap \mathfrak{O}'_i=\emptyset$ for all~$i\in [n]$.
Indeed, assume by contradiction that for some~$v\in N(u)\cap \mathfrak{O}'_i$  there exists~$j<i$ with~$u$,~$v\in O_j$. Since~$G$ is a trivially perfect graph, the neighborhoods of~$u$ and~$v$ are nested.
Moreover, since~$u$,~$v\in O_{j}$, the nodes have the same degree in~$G$.
Consequently, we have $N(u)\cup \{u\}=N(v)\cup \{v\}$.
Hence, the permutation of~$u$ and~$v$ is an objective value preserving automorphism of~$G$ and therefore, since neither~$u$ nor~$v$ are stabilized, they are both contained in the same orbit.
Thus, one is a follower of the other such that the deletion operation removes one of the nodes.
\end{observation}

\subsection{Proof of Theorem~\ref{thm:SSTperfect}}
\label{sec:sstCliqueProof}

To prove Theorem~\ref{thm:SSTperfect}, we use the strong perfect graph theorem.
This means we will show that $\underline{G}_S$ does not contain any odd holes or antiholes as induced subgraphs.
Recall that an \emph{odd hole} is a cycle of odd length at least $5$ that has no chords, i.e., there exist no edges between two non-consecutive nodes in the cycle, and an \emph{odd antihole} is the complement graph of an odd hole.

Throughout this section, we assume that~$G = (V,E)$ is a weighted \tpgraph and that the family of SST cuts $S$ is stringent.
Further, we assume that SST presolving has been applied to $G$ regarding $S$ before $\underline{G}_S$ has been constructed.

To prove perfectness, we are going to assign orientations to edges in $\underline{G}_S$ emerging from the orientations of involved SST cuts in $S$.
Considering various subgraphs with different orientation configurations will show that on paths and cycles without chords of size four the orientations form alternating patterns.
This insight enables us to prove that $\underline{G}_S$ cannot contain any odd holes or odd antiholes and, thus, it is a perfect graph.

\begin{definition}
We say that an edge $\{u,v\}$ of~$\underline{G}_S$ is \emph{of type $\alpha$} if $\{u,v\}\in E$, \emph{of type~$\beta$} if $\{u,v\}$ was added by an addition operation, and \emph{of type~$\delta$} if $(u,v)\in S$.
Further, we say that an edge $\{u,v\}$ of type~$\beta$ or~$\delta$ is oriented towards $u$ if the involved SST cut has~$u$ as its follower.
To emphasize that an edge is oriented towards~$u$, we say that it is of type~$\beta_u$ (resp.\ $\delta _u$).
For consecutive edges $e_1=\{u_1,u_2\},e_2=\{u_2,u_3\},\dots ,e_k=\{u_k,u_{k+1}\}$ in $\underline{G}_S$ we call a sequence $(\epsilon_1,\dots,\epsilon_k)$ of edge types the \emph{type sequence} of $(e_1,\dots,e_k)$ if for each $i\in [k]$ the type of edge~$e_i$ is~$\epsilon_i$.
\end{definition}

So far it is unclear whether the graph~$\underline{G}_S$ is simple or not.
While an edge of type~$\alpha$ makes an edge of type~$\delta$ impossible between the same nodes because of the deletion presolving, the graph's definition allows edges of type~$\beta$ as well as~$\delta$ between the same two nodes.
Anyhow, in the case where~$\underline{G}_S$ has parallel edges of the same orientation, Lemma~\ref{lem:implied_edge} implies the existence of an edge that contradicts that the graph is presolved by deletion presolving.
The following lemma will show that opposite orientations are also impossible implying that the graph~$\underline{G}_S$ is simple.

\begin{lemma}
\label{lem:no_same_type}
There exists no edge~$\{u,v\}$ in $\underline{G}_S$ of type~$\beta$ or~$\delta$ that is oriented towards $u$ as well as~$v$.
\end{lemma}

\begin{proof}
First, assume that the edge is of type~$\beta _u$ as well as~$\beta _v$.
Then there exist edges $\{u,\ell_i\},\{v,\ell_j\}\in E$ with~$v\in O_{i}$ and~$u\in O_{j}$.
All involved nodes are distinct, because if~$\ell_i=\ell_j$, both~$u$ and~$v$ would have been deleted by deletion operations.
Without loss of generality we assume that~$i<j$.
Then by Lemma~\ref{lem:implied_edge}, $\{\ell_i,\ell_j\}\in E$.
Hence, there exists a walk $(\{u,\ell_i\},\{\ell_i,\ell_j\},\{\ell_j,v\})$ in~$G$ and therefore~$\{\ell_i,v\}\in E$ or~$\{u,\ell_j\}\in E$ because~$G$ is a \tpgraph.
Since~$v\in O_{i}$ and~$u\in O_{j}$ this leads in any case to a contradiction regarding the existence of the edge~$\{u,v\}$ in~$\underline{G}_S$ because of the deletion presolving.
Now, we assume that the edge is of type~$\beta _u$ as well as~$\delta _v$.
Then we get that~$v\in O_i$ with $\{v,\ell_i\}\in E$ for some~$i$ and therefore~$v$ would have been deleted by a deletion operation.
Finally, it clearly is not possible that an edge is of type~$\delta _u$ and~$\delta _v$ by the Schreier-Sims-Algorithm.
\end{proof}

\begin{lemma}
\label{lem:P3}
Let the induced subgraph of~$\underline{G}_S$ on the nodes~$u$,~$v$,~$w$ be a path of length three with edges~$\{u,v\}$ and~$\{v,w\}$.
If the corresponding type sequence does not contain~$\alpha$, the edges have opposite orientations.
This means, the type sequence is either $((\beta/\delta)_u,(\beta/\delta)_w)$ or $((\beta/\delta)_v,(\beta/\delta)_v)$.
\end{lemma}

\begin{proof}
By symmetry, we only need to consider one of the two type sequences. We prove that none of the remaining sequences are possible and distinguish all four cases displayed in the first two lines of Figure~\ref{fig:P34Cycles}.

Case 1: Assume the corresponding type sequence is $(\beta_v,\beta_w)$.
This means there exist edges $\{u,\ell_i\}$, $\{v,\ell_j\}\in E$ with~$v\in O_{i}$ and~$w\in O_{j}$.
All involved nodes are pairwise distinct due to~$G$ being presolved and the considered subgraph being a path.
On the one hand, if~$j<i$, we have~$\{\ell_i,\ell_j\}\in E$ by Lemma~\ref{lem:implied_edge}.
So we have the walk $(\{u,\ell_i\},\{\ell_i,\ell_j\},\{\ell_j,v\})$ in~$G$.
Since~$G$ is a \tpgraph, there must exist~$\{u,\ell_j\}\in E$ or~$\{\ell_i,v\}\in E$.
The second is impossible since~$v$ would have been deleted by a deletion operation.
If~$\{u,\ell_j\}\in E$, the edge~$\{u,w\}$ would have been added by an addition operation.
We thus obtain a contradiction to the assumption that~$\{u,\ell_i,\ell_j,v\}$ induces a path.
On the other hand, if~$i<j$ we get by~$S$ being stringent that $N(\ell_j)\cap \mathfrak{O}'_j=\emptyset$ using Observation~\ref{obs:stringent} and therefore~$v\notin \mathfrak{O}'_j$.
Consequently, $v\in \{\ell_1,\dots ,\ell_{j-1}\}$ and therefore,~$w\in O_j$ implies that~$\{v,w\}\in E$ by Lemma~\ref{lem:implied_edge} which contradicts this edge being of type~$\beta$.

Case 2: Assume the corresponding type sequence is $(\delta_v,\beta_w)$.
This means that~$u=\ell_i$ and there exists an edge~$\{v,\ell_j\}\in E$ with~$w\in O_j$.
Clearly, all involved nodes are distinct because if~$\ell_j=\ell_i$, the node~$v$ would have been deleted by a deletion operation.
If~$j<i$, the existence of the edge~$\{v,\ell_j\}\in E$ and the fact that~$v\in O_i$ implies that~$\{\ell_i,\ell_j\}\in E$ by Lemma~\ref{lem:implied_edge}.
So the edge~$\{u,w\}$ would have been added by an addition operation.
If~$i<j$, we can apply the exact same argument as in Case 1.

Case 3: Assume the corresponding type sequence is $(\beta_v,\delta_w)$.
Then for some~$i\in [n]$ we have~$\ell_i\in N(u)$ and~$v\in O_i$ as well as~$w\in O_j$ with~$\ell_j=v$.
Since~$v\in O_i$, we have~$j>i$ because otherwise~$v$ would have been stabilized.
Hence, we obtain~$w\in O_i$ as a consequence of~$v\in O_i$ and~$w\in O_j$.
Therefore, the edge~$\{u,w\}$ would have been added by an addition operation.

Case 4: Assume the corresponding type sequence is $(\delta_v,\delta_w)$.
Then clearly we have~$(u,w)\in S$ and therefore~$\{u,w\}$ is an edge in~$\underline{G}_S$.
\end{proof}

\begin{figure}[t!]
    \centering

\begin{tikzpicture}[thick]
      \tikzstyle{s} += [circle,draw=black,inner sep=0pt,minimum size=2mm];

      \begin{scope}[yshift=2.5cm, xshift=4.0cm]
      \node[s, label=below:{$u$}] (A0) at (0,0) {};
      \node[s, label=below:{$v$}] (A1) at (1.5,0) {};
      \node[s, label=below:{$w$}] (A2) at (3,0) {};
      \node[s, label=above:{$\ell_i$}] (A3) at (0.75,1) {};
      \node[s, label=above:{$\ell_j$}] (A4) at (2.25,1) {};

      \draw (A0)--(A3);
      \draw (A1)--(A4);
      \draw[-latex] (A0)--(A1) node[midway, above] {$\beta$};
      \draw[-latex] (A1)--(A2) node[midway, above] {$\beta$};

      \path let \p1 = (A1), \p2 = (A3), \n1 = {atan2(\y2-\y1,\x2-\x1)}, \n2 = {veclen(\x2-\x1,\y2-\y1)} in
      node[ellipse, at = ($(\p1)!.5!(\p2)$), rotate=\n1, minimum width=\n2+10pt, minimum height=0.5cm, fill=blue!20, fill opacity=0.4] {};

      \path let \p1 = (A2), \p2 = (A4), \n1 = {atan2(\y2-\y1,\x2-\x1)}, \n2 = {veclen(\x2-\x1,\y2-\y1)} in
      node[ellipse, at = ($(\p1)!.5!(\p2)$), rotate=\n1, minimum width=\n2+10pt, minimum height=0.5cm, fill=blue!20, fill opacity=0.4] {};
      \end{scope}

      \begin{scope}[yshift=2.5cm, xshift=8.0cm]
      \node[s, label=below:{$u=\ell_i$}] (A0) at (0,0) {};
      \node[s, label=below:{$v$}] (A1) at (1.5,0) {};
      \node[s, label=below:{$w$}] (A2) at (3,0) {};
      \node[s, label=above:{$\ell_j$}] (A4) at (2.25,1) {};

      \draw (A1)--(A4);
      \draw[-latex] (A0)--(A1) node[midway, above] {$\delta$};
      \draw[-latex] (A1)--(A2) node[midway, above] {$\beta$};

      \path let \p1 = (A1), \p2 = (A0), \n1 = {atan2(\y2-\y1,\x2-\x1)}, \n2 = {veclen(\x2-\x1,\y2-\y1)} in
      node[ellipse, at = ($(\p1)!.5!(\p2)$), rotate=\n1, minimum width=\n2+10pt, minimum height=0.5cm, fill=blue!20, fill opacity=0.4] {};

      \path let \p1 = (A2), \p2 = (A4), \n1 = {atan2(\y2-\y1,\x2-\x1)}, \n2 = {veclen(\x2-\x1,\y2-\y1)} in
      node[ellipse, at = ($(\p1)!.5!(\p2)$), rotate=\n1, minimum width=\n2+10pt, minimum height=0.5cm, fill=blue!20, fill opacity=0.4] {};
      \end{scope}

      \begin{scope}[xshift=4cm]
      \node[s, label=below:{$u$}] (A0) at (0,0) {};
      \node[s, label=below:{$v=\ell_j$}] (A1) at (1.5,0) {};
      \node[s, label=below:{$w$}] (A2) at (3,0) {};
      \node[s, label=above:{$\ell_i$}] (A3) at (0.75,1) {};

      \draw (A0)--(A3);
      \draw[-latex] (A0)--(A1) node[midway, above] {$\beta$};
      \draw[-latex] (A1)--(A2) node[midway, above] {$\delta$};

      \path let \p1 = (A1), \p2 = (A3), \n1 = {atan2(\y2-\y1,\x2-\x1)}, \n2 = {veclen(\x2-\x1,\y2-\y1)} in
      node[ellipse, at = ($(\p1)!.5!(\p2)$), rotate=\n1, minimum width=\n2+10pt, minimum height=0.5cm, fill=blue!20, fill opacity=0.4] {};

      \path let \p1 = (A2), \p2 = (A1), \n1 = {atan2(\y2-\y1,\x2-\x1)}, \n2 = {veclen(\x2-\x1,\y2-\y1)} in
      node[ellipse, at = ($(\p1)!.5!(\p2)$), rotate=\n1, minimum width=\n2+10pt, minimum height=0.5cm, fill=blue!20, fill opacity=0.4] {};
      \end{scope}

      \begin{scope}[xshift=8cm]
      \node[s, label=below:{$u=\ell_i$}] (A0) at (0,0) {};
      \node[s, label=below:{$v=\ell_j$}] (A1) at (1.5,0) {};
      \node[s, label=below:{$w$}] (A2) at (3,0) {};

      \draw[-latex] (A0)--(A1) node[midway, above] {$\delta$};
      \draw[-latex] (A1)--(A2) node[midway, above] {$\delta$};

      \path let \p1 = (A1), \p2 = (A0), \n1 = {atan2(\y2-\y1,\x2-\x1)}, \n2 = {veclen(\x2-\x1,\y2-\y1)} in
      node[ellipse, at = ($(\p1)!.5!(\p2)$), rotate=\n1, minimum width=\n2+10pt, minimum height=0.5cm, fill=blue!20, fill opacity=0.4] {};

      \path let \p1 = (A2), \p2 = (A1), \n1 = {atan2(\y2-\y1,\x2-\x1)}, \n2 = {veclen(\x2-\x1,\y2-\y1)} in
      node[ellipse, at = ($(\p1)!.5!(\p2)$), rotate=\n1, minimum width=\n2+10pt, minimum height=0.5cm, fill=blue!20, fill opacity=0.4] {};
      \end{scope}

      \begin{scope}[yshift=-2.5cm, xshift=2.5cm]
      \node[s, label=below:{$u$}] (A0) at (0,0) {};
      \node[s, label=below:{$v$}] (A1) at (1.5,0) {};
      \node[s, label=below:{$w$}] (A2) at (3,0) {};
      \node[s, label=below:{$y$}] (A3) at (4.5,0) {};
      \node[s, label=above:{$\ell_i$}] (A4) at (0.75,1) {};
      \node[s, label=above:{$\ell_j$}] (A5) at (3.75,1) {};

      \draw (A1)--(A2);
      \draw (A1)--(A4);
      \draw (A2)--(A5);
      \draw[-latex] (A1)--(A0) node[midway, above] {$\beta$};
      \draw[-latex] (A2)--(A3) node[midway, above] {$\beta$};
      \draw[dashed] (A2)--(A4);
      \draw[dashed] (A1)--(A5);

      \path let \p1 = (A0), \p2 = (A4), \n1 = {atan2(\y2-\y1,\x2-\x1)}, \n2 = {veclen(\x2-\x1,\y2-\y1)} in
      node[ellipse, at = ($(\p1)!.5!(\p2)$), rotate=\n1, minimum width=\n2+10pt, minimum height=0.5cm, fill=blue!20, fill opacity=0.4] {};

      \path let \p1 = (A3), \p2 = (A5), \n1 = {atan2(\y2-\y1,\x2-\x1)}, \n2 = {veclen(\x2-\x1,\y2-\y1)} in
      node[ellipse, at = ($(\p1)!.5!(\p2)$), rotate=\n1, minimum width=\n2+10pt, minimum height=0.5cm, fill=blue!20, fill opacity=0.4] {};
      \end{scope}

      \begin{scope}[yshift=-2.5cm, xshift=8cm]
      \node[s, label=below:{$u=\ell_i$}] (A0) at (0,0) {};
      \node[s, label=below:{$v$}] (A1) at (1.5,0) {};
      \node[s, label=below:{$w$}] (A2) at (3,0) {};
      \node[s, label=below:{$y=\ell_j$}] (A3) at (4.5,0) {};

      \draw (A1)--(A2);
      \draw[-latex] (A0)--(A1) node[midway, above] {$\delta$};
      \draw[-latex] (A3)--(A2) node[midway, above] {$\delta$};

      \path let \p1 = (A1), \p2 = (A0), \n1 = {atan2(\y2-\y1,\x2-\x1)}, \n2 = {veclen(\x2-\x1,\y2-\y1)} in
      node[ellipse, at = ($(\p1)!.5!(\p2)$), rotate=\n1, minimum width=\n2+10pt, minimum height=0.5cm, fill=blue!20, fill opacity=0.4] {};

      \path let \p1 = (A2), \p2 = (A3), \n1 = {atan2(\y2-\y1,\x2-\x1)}, \n2 = {veclen(\x2-\x1,\y2-\y1)} in
      node[ellipse, at = ($(\p1)!.5!(\p2)$), rotate=\n1, minimum width=\n2+10pt, minimum height=0.5cm, fill=blue!20, fill opacity=0.4] {};
      \end{scope}

      \begin{scope}[yshift=-5.5cm, xshift=5cm]
      \node[s, label=below:{$u$}] (A0) at (0,0) {};
      \node[s, label=below:{$v$}] (A1) at (1.5,0) {};
      \node[s, label=above:{$w$}] (A2) at (1.5,1.5) {};
      \node[s, label=above:{$y$}] (A3) at (0,1.5) {};
      \node[s, label=below:{$\ell_i=\ell_j$}] (A4) at (0.75,-0.75) {};

      \draw (A0)--(A1);
      \draw (A2)--(A3);
      \draw (A0)--(A4);
      \draw (A1)--(A4);
      \draw[-latex] (A1)--(A2) node[midway, right] {$\beta$};
      \draw[-latex] (A0)--(A3) node[midway, left] {$\beta$};
      \draw[dashed] (A0)--(A2);
      \draw[dashed] (A1)--(A3);

      \path let \p1 = (A2), \p2 = (A4), \n1 = {atan2(\y2-\y1,\x2-\x1)}, \n2 = {veclen(\x2-\x1,\y2-\y1)} in
      node[ellipse, at = ($(\p1)!.5!(\p2)$), rotate=\n1, minimum width=\n2+10pt, minimum height=0.5cm, fill=blue!20, fill opacity=0.4] {};

      \path let \p1 = (A3), \p2 = (A4), \n1 = {atan2(\y2-\y1,\x2-\x1)}, \n2 = {veclen(\x2-\x1,\y2-\y1)} in
      node[ellipse, at = ($(\p1)!.5!(\p2)$), rotate=\n1, minimum width=\n2+10pt, minimum height=0.5cm, fill=blue!20, fill opacity=0.4] {};
      \end{scope}

      \begin{scope}[yshift=-5.5cm, xshift=9cm]
      \node[s, label=below:{$u$}] (A0) at (0,0) {};
      \node[s, label=below:{$v$}] (A1) at (1.5,0) {};
      \node[s, label=above:{$w$}] (A2) at (1.5,1.5) {};
      \node[s, label=above:{$y$}] (A3) at (0,1.5) {};
      \node[s, label=right:{$\ell_i$}] (A4) at (2.5,0.75) {};
      \node[s, label=left:{$\ell_j$}] (A5) at (-1,0.75) {};

      \draw (A0)--(A1);
      \draw (A2)--(A3);
      \draw (A0)--(A5);
      \draw (A1)--(A4);
      \draw[-latex] (A1)--(A2) node[midway, right] {$\beta$};
      \draw[-latex] (A0)--(A3) node[midway, left] {$\beta$};
      \draw[dashed] (A0)--(A4);
      \draw[dashed] (A1)--(A5);

      \path let \p1 = (A2), \p2 = (A4), \n1 = {atan2(\y2-\y1,\x2-\x1)}, \n2 = {veclen(\x2-\x1,\y2-\y1)} in
      node[ellipse, at = ($(\p1)!.5!(\p2)$), rotate=\n1, minimum width=\n2+10pt, minimum height=0.5cm, fill=blue!20, fill opacity=0.4] {};

      \path let \p1 = (A3), \p2 = (A5), \n1 = {atan2(\y2-\y1,\x2-\x1)}, \n2 = {veclen(\x2-\x1,\y2-\y1)} in
      node[ellipse, at = ($(\p1)!.5!(\p2)$), rotate=\n1, minimum width=\n2+10pt, minimum height=0.5cm, fill=blue!20, fill opacity=0.4] {};
      \end{scope}

    \end{tikzpicture}

    \caption{Different paths and cycles considered in the proofs of Lemma~\ref{lem:P3}, Lemma~\ref{lem:P4} and Lemma~\ref{lem:4cycles}. Edges of type~$\alpha$ are illustrated as undirected edges and edges of type~$\beta$ or~$\delta$ as directed edges. The blue ellipses mark SST cuts in~$S$ while the leaders of the cuts are labeled. Dashed lines display edges that are added by the addition presolving if not already present in~$G$.}
    \label{fig:P34Cycles}
\end{figure}

Lemma~\ref{lem:P3} is a first step towards showing an alternating orientation pattern for holes in~$\underline{G}_S$.
We say that  consecutive edges $(e_1,\dots ,e_k)$ have alternating orientations if either oriented edges~$e_i=\{u_i,u_{i+1}\}$ are oriented towards~$u_{i+1}$ for all even~$i\in [k]$ and oriented towards~$u_i$ for all odd~$i\in [k]$ or the other way around.
So in particular, Lemma~\ref{lem:P3} shows that all paths of length three in~$\underline{G}_S$ have alternating orientations.
In the next lemma, we show that this alternating orientation pattern extends to paths of length four.

\begin{lemma}
\label{lem:P4}
Let the induced subgraph of~$\underline{G}_S$ on the nodes~$u$,~$v$,~$w$,~$y$ be a path of length four with edges $\{u,v\},\{v,w\}$ and~$\{w,y\}$.
Then the edges have alternating orientations.
\end{lemma}

\begin{proof}
If~$\{v,w\}$ is not of type~$\alpha$ this is implied by Lemma~\ref{lem:P3}.
So we only have to consider the type sequences $((\beta/\delta)_u,\alpha,(\beta/\delta)_y)$ and $((\beta/\delta)_v,\alpha,(\beta/\delta)_w)$.
If we have $(\delta_u,\alpha,(\beta/\delta)_y)$ (resp.\ $((\beta/\delta)_u,\alpha,\delta_y)$), we can directly see that the edge~$\{u,w\}$ (resp.~$\{v,y\}$) would have been added by an addition operation and therefore contradict the assumed path structure.

Since we consider a path in~$\underline{G}_S$, for the type sequence $(\beta_u,\alpha,\beta_y)$ the path nodes and leaders~$\ell_i$,~$\ell_j$ involved in the~$\beta$ edges have to be pairwise distinct.
See the first graph of the third line of Figure~\ref{fig:P34Cycles} for an illustration of such a graph.
Therefore, there exists a walk $(\{\ell_i,v\},\{v,w\},$ $\{w,\ell_j\})$ in~$G$.
Thus,~$G$ being trivially perfect implies that~$\{v,\ell_j\}\in E$ or~$\{\ell_i,w\}\in E$.
As a consequence of addition presolving, we then have~$\{v,y\}\in E$ or~$\{u,w\}\in E$ which again contradicts the path structure we assumed.

Now consider the type sequence $(\delta_v,\alpha,\delta_w)$.
It is displayed at the end of the third line of Figure~\ref{fig:P34Cycles}.
Without loss of generality, we assume that~$u=\ell_i$ and~$y=\ell_j$ with~$i<j$.
Since~$S$ is stringent, $N(w)\cap \mathfrak{O}'_j=\emptyset$.
This implies that~$v\notin \mathfrak{O}'_j$.
Consequently, since~$v\in O_i$, $v\in \{\ell_1,\dots ,\ell_{j-1}\}$.
Therefore, the graph automorphism associated to the SST cut~$(y,w)$ has to stabilize~$v$ and the existence of the edge~$\{v,w\}$ implies the existence of the edge~$\{v,y\}$ in~$G$ by Lemma~\ref{lem:implied_edge} again contradicting the assumed path structure.

We can apply the same argument for the remaining type sequences, since they all contain the underlying structure of $(\delta_v,\alpha,\delta_w)$ as a substructure.
In the case where we do not directly obtain the edge~$\{v,y\}\in E$ or~$\{u,w\}\in E$, we have three consecutive edges in~$E$ where~$G$ being trivially perfect implies such an edge or the deletion of~$v$ or~$w$ by a deletion operation.
\end{proof}

\begin{proposition}
\label{prop:noOddHole}
If~$S$ is stringent,~$\underline{G}_S$ does not contain any odd holes.
\end{proposition}

\begin{proof}
  Assume there would exist an odd hole~$H$ in~$\underline{G}_S$.
We label the nodes of~$H$ by~$u_0$,~$u_1,\dots ,u_{p-1}$ such that one has $e_i \define \{u_i,u_{(i+1) \bmod p}\}$ being an edge in~$\underline{G}_S$ with $i \in [p-1]_0=\{0,1,\dots ,p-1\}$ where~$p$ is odd.
We claim that the orientation sequence of the edges~$e_0$,~$e_1$,~$e_2$ extends to~$e_3$ which means that if~$e_3$ is oriented, the orientation is uniquely determined in fitting the alternating orientation pattern given by the orientation sequence of the edges~$e_0$,~$e_1$,~$e_2$.
Applying Lemma~\ref{lem:P4} on~$e_0$,~$e_1$,~$e_2$ implies that if~$e_1$ or~$e_2$ is oriented, the alternating orientation pattern extends to $e_3$.
If~$e_1$ and~$e_2$ are both of type~$\alpha$,~$e_0$ as well as~$e_3$ cannot be of type~$\alpha$ because of~$G$ being trivially perfect and~$G$ being a hole.
It is easy to check that~$e_0$ has to be oriented towards~$u_1$ and~$e_3$ towards~$u_3$ because otherwise~$H$ would have chords.
Hence, again the orientation of~$e_3$ matches the orientation pattern on~$e_0$,~$e_1$,~$e_2$ in an alternating way.
Applying this argument on the following three consecutive edges repeatedly leads to a contradiction no later than when the first three edges are reached again since~$H$ has an odd number of edges.
\end{proof}

\begin{lemma}
\label{lem:4cycles}
Let the induced subgraph of~$\underline{G}_S$ on the nodes~$u$,~$v$,~$w$,~$y$ be a cycle without chords with edges $\{u,v\},\{v,w\}, \{w,y\}$ and~$\{y,u\}$.
We denote the cycle by~$C$.
Then the edges have alternating orientations.
This means that all oriented edges are either oriented towards a node in~$\{u,w\}$ or~$\{v,y\}$.
We say the cycle $C$ is oriented towards~$\{u,w\}$ or respectively~$\{v,y\}$.
\end{lemma}

\begin{proof}
By Lemma~\ref{lem:P3}, all neighboring oriented edges have opposite orientations.
So we only have to consider cycles with type sequence $(\alpha ,(\beta /\delta)_w,\alpha ,(\beta /\delta)_y )$ up to symmetry.
If one of the oriented edges is of type~$\delta$, we directly obtain a chord of the cycle by an addition operation.
Therefore, the sequence $(\alpha ,\beta_w,\alpha ,\beta _y )$ is left to consider.

In this case, we have $\{u,v\},\{v,\ell_i\},\{w,y\},\{\ell_j,u\}\in E$ for some~$\ell_i$,~$\ell_j\in L$ with~$w\in O_i$ and~$y\in O_j$.
Clearly, $\ell_i\notin \{v,w\}$ and $\ell_j\notin\{u,y\}$.
If~$\ell_i=y$ or~$\ell_j=w$ we would have chords of type $\alpha$ and if~$\ell_i=u$ or~$\ell_j=v$ we would have chords of type~$\delta$ in~$C$.
Thus, we know that~$\ell_i$,~$\ell_j\notin \{u,v,w,y\}$.
If~$\ell_i=\ell_j$, the chords~$\{u,w\}$ and~$\{v,y\}$ would have been added by addition operations as displayed in the first graph of the fourth line of Figure~\ref{fig:P34Cycles}.
If~$\ell_i \neq \ell_j$, we have the walk $(\{\ell_i,v\},\{v,u\},\{u,\ell_j\})$ in~$G$ implying that~$\{u,\ell_i\}\in E$ or~$\{v,\ell_j\}\in E$ since~$G$ is a \tpgraph resulting in chords~$\{u,w\}$ or~$\{v,y\}$ of type~$\beta$ in~$C$ as depicted in the second graph of the fourth line of Figure~\ref{fig:P34Cycles}.
\end{proof}

\begin{lemma}
\label{lem:orientedP4}
Let~$H$ be an odd antihole in~$\underline{G}_S$.
We label the nodes of~$H$ by~$u_0, \dots, u_{p-1}$, where the indices~$i$ of~$u$ are actually meant as~$i \bmod p$, such that the edges~$\{u_i,u_{i+1}\}$ are not in~$\underline{G}_S$ for all~$i\in [p-1]_0$ where~$p$ is odd.
We denote the induced subgraph of~$H$ on the nodes~$u_i$,~$u_{i+1}$,~$u_{i+2}$,~$u_{i+3}$ by~$P_4^i$.
Then either all edges of type~$\beta$ or~$\delta$ in~$P_4^i$ are oriented towards nodes in~$\{u_i,u_{i+1}\}$ for all~$i\in [p-1]_0$ or all edges of type~$\beta$ or~$\delta$ in~$P_4^i$ are oriented towards nodes in~$\{u_{i+2},u_{i+3}\}$ for all $i\in [p-1]_0$.
In the first case, we say that~$P_4^i$ has a positive orientation for all~$i$ and in the second case we say that~$P_4^i$ has a negative orientation for all~$i$.
\end{lemma}

\begin{proof}
Since~$H$ is an antihole,~$P_4^i$ is a path of length~$4$ from node~$u_{i+1}$ to~$u_{i+2}$.See Figure~\ref{fig:P4i} for an illustration of~$P_4^i$.
Therefore, we can apply Lemma~\ref{lem:P4} and get that~$P_4^i$ is either oriented towards~$\{u_i,u_{i+1}\}$ or oriented towards~$\{u_{i+2},u_{i+3}\}$.
Assume that not all~$P_4^i$'s have the same positive or negative orientation.
Then there exist~$P_4^i$ and~$P_4^{i+1}$ such that~$P_4^i$ has a negative orientation and~$P_4^{i+1}$ has a positive orientation.
In Figure~\ref{fig:P4i} both paths are displayed.
Since~$P_4^i$ and~$P_4^{i+1}$ share the edge~$\{u_{i+1},u_{i+3}\}$ it has to be of type~$\alpha$.
Further, the induced subgraph of~$H$ on the nodes $u_i,u_{i+1},u_{i+3},u_{i+4}$ is a cycle~$C$ without chords.
By Lemma~\ref{lem:4cycles},~$C$ is either oriented towards~$\{u_i,u_{i+1}\}$ or towards~$\{u_{i+3},u_{i+4}\}$.
Hence, since~$P_4^i$ and~$P_4^{i+1}$ have opposite orientations, one of the edges~$\{u_i,u_{i+3}\}$ or~$\{u_{i+1},u_{i+4}\}$ has to be of type~$\alpha$.
Without loss of generality, we assume that the edge~$\{u_i,u_{i+3}\}$ is of type~$\alpha$.
Then, the type sequence of~$P_4^i$ is $(\alpha,\alpha,(\alpha/\beta/\delta))$.
By the same argument as in the proof of Proposition~\ref{prop:noOddHole}, the edge~$\{u_i,u_{i+2}\}$ has to be oriented towards~$u_i$ contradicting that~$P_4^i$ has a negative orientation.
\end{proof}

\begin{figure}[t!]
    \centering
        \begin{tikzpicture}[thick]
    \tikzstyle{s} += [circle,draw=black,inner sep=0pt,minimum size=2mm];
    \def\s{3}
    \def\d{5}
    \def\r{0.2}

        \node[s] (A) at (360/9*1+10:\s) {};
        \node[s] (B) at (360/9*2+10:\s) {};
        \node[s] (C) at (360/9*3+10:\s) {};
        \node[s] (D) at (360/9*4+10:\s) {};
        \node[s] (E) at (360/9*5+10:\s) {};
        \node[s] (F) at (360/9*6+10:\s) {};
        \node[s] (G) at (360/9*7+10:\s) {};
        \node[s] (H) at (360/9*8+10:\s) {};
        \node[s] (I) at (360/9*9+10:\s) {};

        \draw[gray] (A)--(E);
        \draw[gray] (A)--(F);
        \draw[gray] (B)--(F);
        \draw[gray] (B)--(G);
        \draw[gray] (C)--(G);
        \draw[gray] (C)--(H);
        \draw[gray] (D)--(H);
        \draw[gray] (D)--(I);
        \draw[gray] (E)--(I);

        \draw[black] (A)--(C);
        \draw[black] (A)--(D);
        \draw[black] (B)--(D);
        \draw[black] (B)--(E);
        \draw[black] (C)--(E);
        \draw[black] (C)--(F);
        \draw[black] (C)--(I);
        \draw[black] (D)--(F);
        \draw[black] (D)--(G);
        \draw[black] (E)--(G);
        \draw[black] (E)--(H);
        \draw[black] (F)--(H);
        \draw[black] (F)--(I);

        \draw[orange, ultra thick, shorten >=1pt, shorten <=1pt] (H.north) -- (B.south);
        \draw[orange, ultra thick, shorten >=1pt, shorten <=1pt] (A.south) -- (G.north);
        \draw[orange, ultra thick, shorten >=1pt, shorten <=1pt] (G.north) -- (B.south);
        \draw[orange, ultra thick, shorten >=1pt, shorten <=1pt] (A.south) -- (H.north);

        \draw[blue, ultra thick, shorten >=3pt, shorten <=3pt] (H.east) -- (B.east);
        \draw[blue, ultra thick, shorten >=3pt, shorten <=3pt] (A.east) -- (H.east);
        \draw[blue, ultra thick, shorten >=3pt, shorten <=5pt] (I.east) -- (B.east);

        \draw[green, ultra thick, shorten >=3pt, shorten <=3pt] (H.west) -- (A.west);
        \draw[green, ultra thick, shorten >=4pt, shorten <=3pt] (A.west) -- (G.west);
        \draw[green, ultra thick, shorten >=2pt, shorten <=4pt] (G.west) -- (I.west);

        \node[align=center] at (360/9*1+10:\s+0.6) {$u_{i+1}$};
        \node[align=center] at (360/9*2+10:\s+0.6) {$u_i$};
        \node[align=center] at (360/9*8+10:\s+0.6) {$u_{i+3}$};
        \node[align=center] at (360/9*9+10:\s+0.6) {$u_{i+2}$};
        \node[align=center] at (360/9*7+10:\s+0.6) {$u_{i+4}$};
        \node[align=center, blue] at (360/9*0.5+10:\s+0.6) {$P_4^i$};
        \node[align=center, green] at (360/9*8.5+10:\s+0.6) {$P_4^{i+1}$};
        \node[align=center, orange] at (0:\r) {\textbf{$C$}};
\end{tikzpicture}

    \caption{Odd antihole in~$\underline{G}_S$ (exemplary of size~$9$) considered in the proof of Lemma~\ref{lem:orientedP4}.}
    \label{fig:P4i}
\end{figure}

\begin{proposition}
\label{prop:noOddAntihole}
If~$S$ is stringent,~$\underline{G}_S$ does not contain any odd antiholes and is therefore a perfect graph.
\end{proposition}

\begin{proof}
We assume there exists an induced subgraph~$H$ in~$\underline{G}_S$ that is an odd antihole and label the nodes as we did in Lemma~\ref{lem:orientedP4}.
Since an odd antihole of size~$5$ is an odd hole as well, we know by Proposition~\ref{prop:noOddHole} that this cannot exist.
Hence, we can assume that~$H$ consists of at least~$7$ nodes.
As we did in the proof of Lemma~\ref{lem:orientedP4}, we start by considering induced subgraphs~$P_4^i$ of~$H$ on the nodes $u_i,u_{i+1},u_{i+2},u_{i+3}$.
Without loss of generality, using Lemma~\ref{lem:orientedP4}, we assume that~$P_4^{i}$ has a positive orientation for all~$i\in [p-1]_0$.
As we have seen in the proof of Lemma~\ref{lem:orientedP4}, not both edges~$\{u_{i},u_{i+2}\}$ and~$\{u_{i},u_{i+3}\}$ of~$P_4^{i}$ can be of type~$\alpha$.
Further, consider the induced subgraphs~$C_i^j$ of~$H$ on the nodes $u_j,u_{j+1},u_{i},u_{i+1}$ for~$i$, $j\in [p-1]_0,\ \lvert j - i \rvert > 2$ which are cycles without chords.
The edges~$\{u_{i},u_{i+2}\}$ and~$\{u_{i},u_{i+3}\}$ are both part of~$C_{i+2}^{i-1}$.
Therefore, using Lemma~\ref{lem:4cycles}, we get that~$C_{i+2}^{i-1}$ is oriented towards~$\{u_{i-1},u_{i}\}$ for all~$i\in [p-1]_0$.

We claim that there exist~$i$,~$j\in [p-1]_0$ such that~$C_i^j$ is oriented towards~$\{u_{i},u_{i+1}\}$ and~$C_i^{j+1}$ is oriented towards~$\{u_{j+1},u_{j+2}\}$.
The claim is true because~$C_{i}^{i-3}$ is oriented towards~$\{u_{i-3},u_{i-2}\}$ and~$C_{i}^{i+3}$ is oriented towards~$\{u_{i},u_{i+1}\}$.
Since~$C_i^j$ and~$C_i^{j+1}$ share the edges~$\{u_{j+1},u_{i}\}$ and~$\{u_{j+1},u_{i+1}\}$, they both have to be of type~$\alpha$.
Using~$G$ being trivially perfect, we therefore get that the edges~$\{u_{j+2},u_{i}\}$ and~$\{u_{j+2},u_{i+1}\}$ are both oriented towards~$u_{j+2}$, the edge~$\{u_{j},u_{i}\}$ is oriented towards~$u_{i}$ and the edge~$\{u_{j},u_{i+1}\}$ is oriented towards~$u_{i+1}$.
See Figure~\ref{fig:oddAntihole} for an illustration.
Now consider the cycle~$C_{i-1}^j$.
Since it contains the edge~$\{u_j,u_i\}$ which is oriented towards~$u_{i}$, we get that~$C_{i-1}^j$ is oriented towards~$\{u_{i-1},u_{i}\}$.
Similarly, we get that~$C_{i-1}^{j+1}$ is oriented towards~$\{u_{j+1},u_{j+2}\}$.
Hence, we can repeat the argument for~$C_i^j$ and~$C_i^{j+1}$ on~$C_{i-1}^j$ and~$C_{i-1}^{j+1}$ until we reach~$C_{j+3}^j$ for which we get that the edge~$\{u_j,u_{j+4}\}$ is oriented towards~$u_{j+4}$.
But this contradicts that all~$C_{i+3}^i$ are oriented towards~$\{u_i,u_{i+1}\}$.
\end{proof}

\begin{figure}[t!]
    \centering

    \begin{tikzpicture}[thick]
      \tikzstyle{s} += [circle,draw=black,inner sep=0pt,minimum size=2mm];
      \def\s{4}

      \node[s, yshift=2cm, label=below:{$u_{i+1}$}] (A) at (257.5:\s) {};
      \node[s, yshift=2cm, label=below:{$u_{i}$}] (B) at (282.5:\s) {};
      \node[s, yshift=2cm, label=below:{$u_{i-1}$}] (C) at (307.5:\s) {};

      \node[s, label=above:{$u_{j}$}] (D) at (115:\s) {};
      \node[s, label=above:{$u_{j+1}$}] (E) at (90:\s) {};
      \node[s, label=above:{$u_{j+2}$}] (F) at (65:\s) {};
      \node[s, label=above:{$u_{j+3}$}] (G) at (40:\s) {};
      \node[s, label=above:{$u_{j+4}$}] (H) at (15:\s) {};

      \draw[ultra thick, orange, shorten >=0.1cm,shorten <=0.1cm] (E)--(A);
      \draw[ultra thick, orange, shorten >=0.1cm,shorten <=0.1cm] (E)--(B);
      \draw[ultra thick, blue, shorten >=0.1cm,shorten <=0.1cm,-latex] (D)--(A);
      \draw[ultra thick, blue, shorten >=0.1cm,shorten <=0.1cm,-latex] (D)--(B) node [midway, black, fill=white, inner sep=2pt,rotate=15]{$C_i^j$};
      \draw[ultra thick, green, shorten >=0.1cm,shorten <=0.1cm,-latex] (A) -- (F);
      \draw[ultra thick, green, shorten >=0.1cm,shorten <=0.1cm,-latex] (B) -- (F);

      \draw[ultra thick, orange, shorten >=0.1cm,shorten <=0.1cm] (E)--(C);
      \draw[ultra thick, blue, shorten >=0.1cm,shorten <=0.1cm,-latex] (D)--(C);
      \draw[ultra thick, green, shorten >=0.1cm,shorten <=0.1cm,-latex] (C) -- (F);

      \draw[ultra thick, black, shorten >=0.1cm,shorten <=0.1cm] (E)--(H);
      \draw[ultra thick, blue, shorten >=0.1cm,shorten <=0.1cm,-latex] (D)--(H);
      \draw[ultra thick, black, shorten >=0.1cm,shorten <=0.1cm] (E)--(G);
      \draw[ultra thick, black, shorten >=0.1cm,shorten <=0.1cm] (D)--(G);

      \draw[loosely dotted, ultra thick, shorten >=0.7cm,shorten <=0.7cm ] (C) to [bend right] (H);

      \node [black, fill=white, inner sep=2pt, rotate=-20] at (65:\s-0.8) {$C_{j+3}^{j}$};
      \node [black, fill=white, inner sep=2pt, rotate=-20] at ($(A)!0.5!(F)$) {$C_i^{j+1}$};

    \end{tikzpicture}

    \caption{Parts of an odd antihole in~$\underline{G}_S$ considered in the proof of~\ref{prop:noOddAntihole}. The orange edges are edges of type~$\alpha$ and the blue and green edges are oriented edges with their orientation indicated through directed arcs.}
    \label{fig:oddAntihole}
\end{figure}

\section{Computational Results}
\label{sec:ComputationalResults}

In this section we report on computational experiments with SST cuts and
presolving that in particular investigate the impact of the order in which
leaders are selected.

Our implementation is based on the framework SCIP~\cite{BestuzhevaEtAl23}
and the branch-and-cut solver for stable set problems BACS
(\url{https://github.com/dopt-TUDa/bacs}). One main component of the
implementation is a clique separator, which first applies a greedy
algorithm and, if unsuccessful, a combinatorial algorithm for the maximal
weight clique problem implemented in SCIP.
So in particular, BACS does not use the edge formulation of the problem
but instead separates clique inequalities in the branch-and-bound tree throughout
the solving process. Moreover, the implementation
features a presolver that automatically detects connected components and solves them individually.
All other specialized presolving/propagation methods of BACS are turned off.
SST presolving is performed on the original graph before the optimization process starts.
After computing symmetry for SST presolving, BACS first checks whether there are symmetric components and possibly merges them.
This means only one symmetric component is kept in the graph as a representative with adapted objective.
More precisely, the variables of all except one of the symmetric components are set to~$0$ while the objective value of the remaining components variables are updated by the number of symmetric components times their initial objective values.
To be able to translate the solution of the transformed problem to a solution of the initial problem, we store the correspondence of nodes and their representatives.
Additionally, isolated nodes of nonnegative weight get directly fixed to 1.
SST presolving is then iteratively run until no further modifications of the graph are made.
Instead of computing full stabilizers for the leaders and followers, we filter the generators and therefore only operate on a subgroup of the actual symmetry group of the current graph.
SST clique cuts are separated by iterating through all
leader orbits and calling the above mentioned combinatorial algorithm for
finding maximal weight cliques from SCIP.
Since SCIP presolving might modify the graph after the SST presolving, we recompute symmetry if additionally SST cuts or SST clique cuts are separated.
If leaders and followers have been computed for SST cuts, we reuse them to separate SST clique cuts.

To detect automorphisms of graphs, we apply bliss~\cite{bliss} and sassy~\cite{andSS23}.
The corresponding running time
is usually very small---the maximal time for one exceptional instance was
2.50 seconds, see Table~\ref{tab:instance_statistics} in the appendix. To
highlight the effect on the dual bound, we initialized the runs with a
cutoff using the best primal value. All primal heuristics are turned off.

\paragraph{Computational setup} For our computations we use a developer version of SCIP
10.0.3~\cite{SCIPOptSuite10} (githash: a396028) and CPLEX 22.1.2.0 as LP
solver. The experiments were run on a Linux cluster with 3.5 GHz Intel Xeon
E5-1620 Quad-Core CPUs, having 32~GB main memory and 10~MB cache each. All
computations were run single-threaded and with a time limit of one
hour. Detailed results for each instance can be found in an online
supplement at \url{https://www2.mathematik.tu-darmstadt.de/~pfetsch/SST-cuts.html}.

\paragraph{Instances}

We selected instances containing symmetry from various sources, see Table~\ref{tab:instances}, to generate a testset which we call \emph{Mixed-testset}.
Many of the instances were designed for the
maximum clique problem (as indicated by column ``clique'' in
Table~\ref{tab:instances}). We use the complemented graphs for these
instances.
To obtain results less sensitive to the particular node labeling of the instances, we added two copies of each instance with differently permuted node labels.
For the
selection of presented instances, we then ran BACS with all different settings used in each of our tests, on all three permuted versions of the instances listed in this table (plus some additional ones).
We removed all instances which could be solved in less than 1 second for all settings.
This way the aggregated results become less homogenized and therefore more informative.
All the selected instances make up the Mixed-testset containing 526 instances.
Note that the testset does not necessarily contain each instance in all three permuted versions, because of variances in running times among different instance versions.
We excluded the instance \texttt{15-22-10-10-compl-weigh} from the Mixed-testset, since this weighted graph with over 37 million edges exceeded the available memory.

Further, we generated a testset consisting of~266 Johnson graphs, which we call \emph{Johnson-testset}.
For this testset, we considered all Johnson graphs with at least 100 and at most \num{10000} nodes and number of edges between \num{1000} and \num{1000000}.
Again, we additionally considered two differently permuted versions of the generated graphs and selected the graphs for the testset exactly like we did for the Mixed-testset.

For the third testset, we considered graphs whose nodes correspond to binary code words of length~$n$ and whose edges exist between all nodes whose Hamming distance is at most $d$.
We denote these graphs by \emph{BinaryCode\_$n$\_$d$} and call the corresponding testset \emph{BinaryCode-testset}.
Following the exact same approach we used for the Johnson-testset to create this testset, we obtained~88 instances in the BinaryCode-testset.

The Johnson- and BinaryCode-testset are available on github (\url{https://github.com/dopt-TUDa/SymGraphs}).

\begin{table}[tb]
  \caption{Sources and selection of instances for Mixed-testset.}
  \label{tab:instances}
  \begin{tabularx}{\textwidth}{@{}l@{\;\;}r@{\;\;}r@{\;\;}r@{\;\;}r@{\;\;}r@{\;\;}>{\raggedright\arraybackslash}X@{}}\toprule
    name & ref. & clique & weighted &total & chosen & notes\\\midrule
    aim\_A    & \cite{CliSAT}    & yes &  no &  44 &  1 & constraint satisfaction instances\\
    Coloring  &                  &  no & yes & 105 &  4 & pricing problems for
                                                     solving graph coloring
                                                     problems supplied by Rolf van der Hulst\\
    Dimacs    & \cite{DIMACS}    & yes & no  &  66 & 12 & standard max.\ clique benchmark set\\
    ECC       & \cite{ECC}       &  no & yes &  15 & 15 & error correcting codes\\
    ECC-C     &                  & yes &  no &  15 & 14 & complemented unweighted versions of ECC\\
    ECC-CW    &                  & yes & yes &  15 & 14 & complemented weighted versions of ECC\\
    ehi\_A    & \cite{CliSAT}    & yes &  no &  12 & 10 & constraint satisfaction instances\\
    ehi\_B    & \cite{CliSAT}    & yes &  no &  13 &  8 & constraint satisfaction instances\\
    Kidney    & \cite{Kidney}    & yes &  no &  20 & 20 & kidney exchange instances from~\cite{MaxClique}\\
    Kidney-W  &                  & yes & yes &  20 & 20 & kidney exchange instances from~\cite{MaxClique}\\
    Monotone  & \cite{CliSAT}    & yes &  no &   3 &  3 & monotone matrices\\
    OEIS      & \cite{OEIS}      &  no &  no &  32 & 17 & Online Encyclopedia of Integer Sequences\\
    REF-W     & \cite{MaxClique} & yes & yes & 100 & 47 & Research Excellence Framework \\
    UAI       & \cite{UAI}       &  no & yes &  79 &  6 & UAI Competition 2014 (some $\infty$ weights) \\
    VC        & \cite{CliSAT}    & yes &  no &  55 &  1 & based on vertex cover instances from~\cite{PACS19}\\
    \bottomrule
  \end{tabularx}
\end{table}

\paragraph{Settings}
In our experiments, we compare various variants of combining SST (clique) cuts and presolving.
More precisely, we are considering the settings displayed in Table~\ref{tab:settings}.
Column ``SST presolving'' indicates whether SST presolving is applied and the column ``add edges'' gives the additional information whether the edges found in SST presolving are added to the graph or not.
``SCIP symmetry'' indicates whether SCIP symmetry handling is applied.
The last two columns display if SST (clique) cuts are separated.
Finally, the ``x'' in the settings name is a placeholder for the concrete rule how to choose leaders used for the presolving and/or the SST (clique) cuts.
We considered three different orbitrules:
\begin{enumerate}[label=-]
\item min: while choosing the order of the leaders,
  pick the next leader as the first element from a smallest nontrivial orbit,
\item max: same as ``min'' but using largest orbits,
\item str: choose stringent leaders (Section~\ref{sec:triviallyPerfect}) from a largest nontrivial orbit.
\end{enumerate}

\begin{table}
  \caption{Description of the experimental settings.}
  \label{tab:settings}
  \begin{tabularx}{\textwidth}{@{}>{\raggedright\arraybackslash}Xrrrrr@{}}\toprule
    Setting & SST presolving & add edges & SCIP symmetry & SST cuts & SST clique cuts \\ \midrule
    base             & no  & no  & no  & no  & no  \\
    base-sym         & no  & no  & yes & no  & no  \\
    SST-pre-x        & yes & yes & no  & no  & no  \\
    SST-pre-x-ne     & yes & no  & no  & no  & no  \\
    SST-pre-x-sym    & yes & yes & yes & no  & no  \\
    SST-pre-x-ne-sym & yes & no  & yes & no  & no  \\
    SST-C            & no  & no  & no  & yes & no  \\
    SST-CC           & no  & no  & no  & no  & yes \\
    SST-CCC          & no  & no  & no  & yes & yes \\
    SST-pre-x-C      & yes & yes & no  & yes & no  \\
    SST-pre-x-CC     & yes & yes & no  & no  & yes \\
    SST-pre-x-CCC    & yes & yes & no  & yes & yes \\
    \bottomrule
  \end{tabularx}
\end{table}

\paragraph{Evaluation of experiments}
Table~\ref{tab:all_SSTtestset} presents a summary of our results for the Mixed-testset.
For each setting, Table~\ref{tab:all_SSTtestset} gives the number of
instances solved to optimality (from 526) (column ``\#opt''), the shifted
geometric mean\footnote{The shifted geometric mean of
  values~$t_1,\dots,t_n$ is defined
  as~$\big(\prod_{i=1}^n (t_i + s)\big)^{1/n} -s$, where the shift~$s$ is
  100 for the number of nodes and 1 for time. Unsolved instances contribute with 3600 seconds.}  of the total time in
seconds and the number of nodes in the branch-and-bound tree. Then for SST
presolving, the average number of fixed nodes, number of added edges, and
the SST presolving time in seconds are presented. Finally, we list the
average of the total number of leaders $\card{L}$ and followers
$\card{S(L)}$.
Note that $\card{L}$ and $\card{S(L)}$ represent the sum of all computed
leaders and followers, respectively, in variants that perform SST presolving.
For variants where only cuts are separated, symmetry is computed only once and reported in Table~\ref{tab:all_SSTtestset}.
For the three variants that combine SST presolving and the separation of SST (clique) cuts, the numbers of leaders and followers computed in the SST presolving and the ones recomputed for the cut separation are added up.

The results of Table~\ref{tab:all_SSTtestset} show that using SST presolving or cuts
significantly improves on the \textit{base} settings. However, using \textit{base-sym} is
clearly the fastest. Using the ``min'' orbit rule is significantly worse that
using any of the other tested rules. The other rules behave quite
similarly. Here not generating edges always performs worse. Additionally, using
the symmetry handling of SCIP does not help to improve performance
and is even counterproductive. This is likely, because performing SST presolving removes a lot of
symmetry, such that the remaining part is not enough to potentially improve
performance too much (see also~\cite{JaegerPfetsch2025}). The much better performance of the SCIP symmetry
handling methods in comparison to the SST based methods seems to come from
the ability to dynamically handle symmetry within the tree.
Interestingly, just using SST cuts (\textit{SST-C}) already provides a big
improvement.
Comparing the number of leaders when using SST presolving with the stringent orbitrule (\textit{SST-pre-str}) and the settings that combine SST presolving and cuts, we see that after SST presolving not much symmetry is left.
Therefore, the number of separated SST (clique) cuts is small, resulting in a modest effect on the solution process.

\begin{table}
  \caption{Comparison of different SST variants for the Mixed-testset.}
  \label{tab:all_SSTtestset}
  \begin{tabular*}{\textwidth}{@{}l@{\;\extracolsep{\fill}}rrrrrrrr@{}}\toprule
    & & & & \multicolumn{3}{c}{SST presolving} &\\\cmidrule(lr){5-7}
    Setting            & \#opt & time & \#nodes & \#fixed & \#edges & time & $\card{L}$ & $\card{S(L)}$\\
    \midrule
    base                 & \num{ 432} & \num{   19.58} & \num{   448.8} & \num{   0.0} & \num{     0.0} & \num{  0.00} & \num{   0.0} & \num{   0.0}\\
    base-sym             & \num{ 467} & \num{   15.28} & \num{   250.7} & \num{   0.0} & \num{     0.0} & \num{  0.00} & \num{   0.0} & \num{   0.0}\\
    SST-pre-min          & \num{ 435} & \num{   21.49} & \num{   445.7} & \num{  57.2} & \num{   235.0} & \num{  2.67} & \num{  85.2} & \num{ 421.1}\\
    SST-pre-min-ne       & \num{ 435} & \num{   22.24} & \num{   466.4} & \num{  57.2} & \num{     0.0} & \num{  2.66} & \num{  82.6} & \num{ 390.5}\\
    SST-pre-min-sym      & \num{ 435} & \num{   22.71} & \num{   445.9} & \num{  57.2} & \num{   235.0} & \num{  2.66} & \num{  85.2} & \num{ 421.1}\\
    SST-pre-min-ne-sym   & \num{ 440} & \num{   21.10} & \num{   399.7} & \num{  57.2} & \num{     0.0} & \num{  2.64} & \num{  82.6} & \num{ 390.5}\\
    SST-pre-max          & \num{ 449} & \num{   17.71} & \num{   358.9} & \num{ 117.2} & \num{  8806.7} & \num{  2.59} & \num{  88.5} & \num{ 602.1}\\
    SST-pre-max-ne       & \num{ 448} & \num{   19.24} & \num{   413.6} & \num{ 113.7} & \num{     0.0} & \num{  2.55} & \num{  83.3} & \num{ 561.5}\\
    SST-pre-max-sym      & \num{ 449} & \num{   18.81} & \num{   358.7} & \num{ 117.2} & \num{  8806.7} & \num{  2.58} & \num{  88.5} & \num{ 602.1}\\
    SST-pre-max-ne-sym   & \num{ 449} & \num{   19.72} & \num{   389.8} & \num{ 113.7} & \num{     0.0} & \num{  2.55} & \num{  83.3} & \num{ 561.5}\\
    SST-pre-str          & \num{ 450} & \num{   17.74} & \num{   356.9} & \num{ 117.3} & \num{  8794.6} & \num{  2.62} & \num{  88.6} & \num{ 603.6}\\
    SST-pre-str-ne       & \num{ 449} & \num{   19.19} & \num{   410.5} & \num{ 113.8} & \num{     0.0} & \num{  2.57} & \num{  83.4} & \num{ 562.2}\\
    SST-pre-str-sym      & \num{ 449} & \num{   18.76} & \num{   354.0} & \num{ 117.3} & \num{  8794.6} & \num{  2.59} & \num{  88.6} & \num{ 603.6}\\
    SST-pre-str-ne-sym   & \num{ 449} & \num{   19.69} & \num{   388.1} & \num{ 113.8} & \num{     0.0} & \num{  2.58} & \num{  83.4} & \num{ 562.2}\\
    SST-C                & \num{ 448} & \num{   17.60} & \num{   343.6} & \num{   0.0} & \num{     0.0} & \num{  0.00} & \num{  68.4} & \num{ 398.8}\\
    SST-CC               & \num{ 445} & \num{   17.71} & \num{   355.3} & \num{   0.0} & \num{     0.0} & \num{  0.00} & \num{  53.4} & \num{ 200.4}\\
    SST-CCC              & \num{ 447} & \num{   17.96} & \num{   341.0} & \num{   0.0} & \num{     0.0} & \num{  0.00} & \num{  68.4} & \num{ 398.8}\\
    SST-pre-C            & \num{ 449} & \num{   17.74} & \num{   354.9} & \num{ 117.3} & \num{  8794.6} & \num{  2.59} & \num{ 103.4} & \num{ 763.4}\\
    SST-pre-CC           & \num{ 449} & \num{   17.75} & \num{   356.9} & \num{ 117.3} & \num{  8794.6} & \num{  2.60} & \num{  89.6} & \num{ 607.2}\\
    SST-pre-CCC          & \num{ 450} & \num{   17.79} & \num{   354.9} & \num{ 117.3} & \num{  8794.6} & \num{  2.61} & \num{ 103.4} & \num{ 763.4}\\
    \bottomrule
  \end{tabular*}
\end{table}

\begin{table}
  \caption{Comparison of different SST variants for the Johnson-testset.}
  \label{tab:all_Johnson}
  \begin{tabular*}{\textwidth}{@{}l@{\;\extracolsep{\fill}}rrrrrrrr@{}}\toprule
    & & & & \multicolumn{3}{c}{SST presolving} &\\\cmidrule(lr){5-7}
    Setting            & \#opt & time & \#nodes & \#fixed & \#edges & time & $\card{L}$ & $\card{S(L)}$\\
    \midrule
    base                 & \num{ 124} & \num{  397.93} & \num{  1936.7} & \num{   0.0} & \num{     0.0} & \num{  0.00} & \num{   0.0} & \num{   0.0}\\
    base-sym             & \num{ 152} & \num{  135.83} & \num{   772.0} & \num{   0.0} & \num{     0.0} & \num{  0.00} & \num{   0.0} & \num{   0.0}\\
    SST-pre-min          & \num{ 146} & \num{  123.08} & \num{  1117.4} & \num{ 234.6} & \num{  2814.3} & \num{  0.37} & \num{  19.3} & \num{1785.8}\\
    SST-pre-min-ne       & \num{ 141} & \num{  142.13} & \num{  1250.9} & \num{ 234.8} & \num{     0.0} & \num{  0.45} & \num{  19.3} & \num{1799.3}\\
    SST-pre-min-sym      & \num{ 146} & \num{  123.72} & \num{  1130.2} & \num{ 234.6} & \num{  2814.3} & \num{  0.37} & \num{  19.3} & \num{1785.8}\\
    SST-pre-min-ne-sym   & \num{ 176} & \num{   57.19} & \num{   514.7} & \num{ 234.8} & \num{     0.0} & \num{  0.44} & \num{  19.3} & \num{1799.3}\\
    SST-pre-max          & \num{ 170} & \num{   63.36} & \num{   738.6} & \num{ 327.7} & \num{ 38247.1} & \num{  0.38} & \num{  25.2} & \num{2764.0}\\
    SST-pre-max-ne       & \num{ 167} & \num{   71.37} & \num{   807.6} & \num{ 340.1} & \num{     0.0} & \num{  0.38} & \num{   8.8} & \num{2688.5}\\
    SST-pre-max-sym      & \num{ 170} & \num{   63.68} & \num{   737.7} & \num{ 327.7} & \num{ 38247.1} & \num{  0.37} & \num{  25.2} & \num{2764.0}\\
    SST-pre-max-ne-sym   & \num{ 167} & \num{   67.19} & \num{   757.1} & \num{ 340.1} & \num{     0.0} & \num{  0.38} & \num{   8.8} & \num{2688.5}\\
    SST-pre-str          & \num{ 170} & \num{   64.42} & \num{   755.9} & \num{ 327.5} & \num{ 38241.1} & \num{  0.38} & \num{  25.2} & \num{2771.5}\\
    SST-pre-str-ne       & \num{ 167} & \num{   70.37} & \num{   801.6} & \num{ 340.9} & \num{     0.0} & \num{  0.39} & \num{   9.5} & \num{2691.9}\\
    SST-pre-str-sym      & \num{ 170} & \num{   64.66} & \num{   755.9} & \num{ 327.5} & \num{ 38241.1} & \num{  0.38} & \num{  25.2} & \num{2771.5}\\
    SST-pre-str-ne-sym   & \num{ 167} & \num{   65.68} & \num{   746.0} & \num{ 340.9} & \num{     0.0} & \num{  0.39} & \num{   9.5} & \num{2691.9}\\
    SST-C                & \num{ 171} & \num{  103.27} & \num{   830.4} & \num{   0.0} & \num{     0.0} & \num{  0.00} & \num{   4.7} & \num{2558.6}\\
    SST-CC               & \num{ 160} & \num{  176.96} & \num{   857.6} & \num{   0.0} & \num{     0.0} & \num{  0.00} & \num{   4.0} & \num{1755.8}\\
    SST-CCC              & \num{ 165} & \num{  124.67} & \num{   825.6} & \num{   0.0} & \num{     0.0} & \num{  0.00} & \num{   4.7} & \num{2558.6}\\
    SST-pre-C            & \num{ 170} & \num{   64.22} & \num{   753.6} & \num{ 327.5} & \num{ 38241.1} & \num{  0.38} & \num{  27.5} & \num{2779.6}\\
    SST-pre-CC           & \num{ 170} & \num{   64.24} & \num{   757.8} & \num{ 327.5} & \num{ 38241.1} & \num{  0.38} & \num{  25.8} & \num{2772.9}\\
    SST-pre-CCC          & \num{ 170} & \num{   64.20} & \num{   753.1} & \num{ 327.5} & \num{ 38241.1} & \num{  0.38} & \num{  27.5} & \num{2779.6}\\
    \bottomrule
  \end{tabular*}
\end{table}

\begin{table}
  \caption{Comparison of different SST variants for the BinaryCode-testset.}
  \label{tab:all_BinaryCode}
  \begin{tabular*}{\textwidth}{@{}l@{\;\extracolsep{\fill}}rrrrrrrr@{}}\toprule
    & & & & \multicolumn{3}{c}{SST presolving} &\\\cmidrule(lr){5-7}
    Setting            & \#opt & time & \#nodes & \#fixed & \#edges & time & $\card{L}$ & $\card{S(L)}$\\
    \midrule
    base                 & \num{  49} & \num{  246.17} & \num{  1070.3} & \num{   0.0} & \num{     0.0} & \num{  0.00} & \num{   0.0} & \num{   0.0}\\
    base-sym             & \num{  64} & \num{  122.34} & \num{   346.9} & \num{   0.0} & \num{     0.0} & \num{  0.00} & \num{   0.0} & \num{   0.0}\\
    SST-pre-min          & \num{  61} & \num{   31.91} & \num{   419.2} & \num{ 414.5} & \num{   134.9} & \num{  1.38} & \num{  92.3} & \num{34943.2}\\
    SST-pre-min-ne       & \num{  61} & \num{   31.85} & \num{   425.8} & \num{ 414.5} & \num{     0.0} & \num{  1.39} & \num{  92.4} & \num{34943.6}\\
    SST-pre-min-sym      & \num{  61} & \num{   32.00} & \num{   419.7} & \num{ 414.5} & \num{   134.9} & \num{  1.40} & \num{  92.3} & \num{34943.2}\\
    SST-pre-min-ne-sym   & \num{  61} & \num{   31.91} & \num{   425.7} & \num{ 414.5} & \num{     0.0} & \num{  1.41} & \num{  92.4} & \num{34943.6}\\
    SST-pre-max          & \num{  65} & \num{   19.80} & \num{   252.5} & \num{ 447.7} & \num{ 27531.2} & \num{  1.37} & \num{  89.5} & \num{35505.0}\\
    SST-pre-max-ne       & \num{  64} & \num{   22.68} & \num{   306.0} & \num{ 444.4} & \num{     0.0} & \num{  1.36} & \num{  86.5} & \num{35508.4}\\
    SST-pre-max-sym      & \num{  65} & \num{   19.84} & \num{   252.6} & \num{ 447.7} & \num{ 27531.2} & \num{  1.37} & \num{  89.5} & \num{35505.0}\\
    SST-pre-max-ne-sym   & \num{  64} & \num{   22.70} & \num{   303.6} & \num{ 444.4} & \num{     0.0} & \num{  1.35} & \num{  86.5} & \num{35508.4}\\
    SST-pre-str          & \num{  65} & \num{   19.94} & \num{   254.0} & \num{ 447.5} & \num{ 27516.6} & \num{  1.41} & \num{  89.6} & \num{35505.2}\\
    SST-pre-str-ne       & \num{  64} & \num{   22.62} & \num{   302.7} & \num{ 444.4} & \num{     0.0} & \num{  1.39} & \num{  86.5} & \num{35508.5}\\
    SST-pre-str-sym      & \num{  65} & \num{   19.97} & \num{   253.9} & \num{ 447.5} & \num{ 27516.6} & \num{  1.41} & \num{  89.6} & \num{35505.2}\\
    SST-pre-str-ne-sym   & \num{  64} & \num{   22.70} & \num{   299.6} & \num{ 444.4} & \num{     0.0} & \num{  1.40} & \num{  86.5} & \num{35508.5}\\
    SST-C                & \num{  65} & \num{   65.92} & \num{   269.5} & \num{   0.0} & \num{     0.0} & \num{  0.00} & \num{  81.8} & \num{35361.0}\\
    SST-CC               & \num{  58} & \num{  157.41} & \num{   279.6} & \num{   0.0} & \num{     0.0} & \num{  0.00} & \num{   2.9} & \num{1207.2}\\
    SST-CCC              & \num{  61} & \num{  117.67} & \num{   252.2} & \num{   0.0} & \num{     0.0} & \num{  0.00} & \num{  81.8} & \num{35361.0}\\
    SST-pre-C            & \num{  65} & \num{   19.93} & \num{   253.9} & \num{ 447.5} & \num{ 27516.6} & \num{  1.41} & \num{  89.7} & \num{35505.3}\\
    SST-pre-CC           & \num{  65} & \num{   19.95} & \num{   254.0} & \num{ 447.5} & \num{ 27516.6} & \num{  1.43} & \num{  89.7} & \num{35506.2}\\
    SST-pre-CCC          & \num{  65} & \num{   19.94} & \num{   253.8} & \num{ 447.5} & \num{ 27516.6} & \num{  1.41} & \num{  89.7} & \num{35505.3}\\
    \bottomrule
  \end{tabular*}
\end{table}

The picture changes when considering the Johnson and BinaryCode
testsets. Here, the SST based methods are clearly superior to \textit{base} by
around a factor of 6 and also to \textit{base-sym} by a factor of about
2 (see Table~\ref{tab:all_Johnson} and Table~\ref{tab:all_BinaryCode}). Interestingly, \textit{SST-pre-min-ne-sym} performs best on the Johnson-testset. But we suspect
that this favorable performance highly depends on the particluar structure of Johnson graphs.
Nonetheless, we cannot explain why on Johnson graphs \textit{SST-pre-min-ne-sym} achieves better performance than the other orbitrules in the same setting, while it is strongly outperformed by the other orbitrules in the other three settings.

For the three settings \textit{SST-pre-min}, \textit{SST-pre-max}, and \textit{SST-pre-str} that
perform full SST presolving while not using SCIP symmetry handling, we observe a significant number of operations:
between \num{57.2} ($\approx$ \SI{1.7}{\percent}) and \num{117.3} ($\approx$
\SI{3.4}{\percent}) variables are removed (fixed to 0) and between \num{235.0}
($\approx$ \SI{0.02}{\percent}) and \num{8806.7} ($\approx$
\SI{0.73}{\percent}) edges are added, on average on the Mixed-testset.
On the Johnson-testset and BinaryCode-testset the number of operations is
considerably higher with up to \SI{30.9}{\percent} of the variables being
fixed and the number of edges being increased by up to
\SI{21.6}{\percent}.
Johnson graphs and binary code instances apparently having a structure favorable for SST presolving seems correlated with the SST presolving being very effective on the corresponding testsets.
With one exception which we discussed earlier (\textit{SST-pre-min-ne-sym} on the Johnson-testset), no significant improvement but sometimes worsening of performance can be seen when no edges are added in the SST presolving or SCIP symmetry handling is additionally used.
This emphasizes the usage of both SST presolving operations as an alternative for rather than an addition to SCIP symmetry handling.
The overall fastest of these three
presolving options is \textit{SST-pre-max} (closely followed by \textit{SST-pre-str}) with a
speed-up of about \SI{10}{\percent} in relation to \textit{base} on the Mixed-testset.
It also solved 17 more instances (\textit{SST-pre-str} solved 18 more instances).
The other presolving variant \textit{SST-pre-min}
solves more instances than the default, but has a worse running time.
Generally, with respect to both measures, it performs worse than \textit{SST-pre-str} and
\textit{SST-pre-max}.
All these observations can be confirmed on the other two testsets (except \textit{SST-pre-min-ne-sym} on the Johnson-testset) and are even present to a much greater extent.

To investigate the influence of the SST presolving and the SST (clique) cuts on the LP relaxations, we ran our computations again, but this time only processed the root node.
Additionally, we turned off the propagator that causes the saparate consideration of connected components.
Comparing the gaps between the dual bound and the primal bound we obtain very similar relative performance relations as for the complete computations among the different ``min'', ``max'' and ``str'' settings on the three testsets.
Anyhow, the application of SST (clique) cuts and especially of the SST presolving seems to strengthen the LP relaxations on all testsets since the average gaps for \textit{base} as well as \textit{base-sym} are larger than for all the other settings on all testsets.
Further, the unexpected superior performance for the \textit{SST-pre-min-ne-sym} setting does not seem to be caused by a beneficial interplay with the LP relaxations.
See the online supplement for more details.

Overall these results nicely support and complement the theoretical
results: SST presolving is easy to use and a very valuable tool. For
general stable set instances, it is quite effective, but not as effective
as established symmetry handling methods. This is likely, because the symmetry groups only move
a small number of variables and SST presolving is a static method. For
instances like Johnson graphs and binary code instances, the usage of SST
presolving is faster by a factor of at least about 2 and 6,
respectively. This highlights the dramatic benefits of SST cuts on these
instances.

The selection of the leaders has significant impact both theoretically as
well as practically. Exploiting graph structure as done for SST clique cuts
helps for the polyhedral results and also slightly speeds up the solution.

\section{Conclusions and Outlook}

Concerning our leading question from the introduction, our theoretical
results show that with respect to a computational complexity and polyhedral
perspective, there is a---maybe surprising---dependency on the kind of SST
cuts (stringent vs.\ others).  One open question is thus whether $P(G,S)$
can be separated in polynomial time if $G$ is a perfect graph and $S$ are
stringent. Moreover, an alternative way to prove
Theorem~\ref{thm:SSTclique} would be to show that the graph $G'$
corresponding to the extended formulation in
Proposition~\ref{prop:extendedFormulation} is perfect for stringent SST
cuts. We have not succeeded in this direction, but it is an open question
for which classes of perfect graphs and SST cuts, $G'$ remains perfect.

\bibliographystyle{plain}
\bibliography{sst_cuts_preprint_V2}

@Misc{DIMACS,
  Title =	 {2nd {DIMACS} challenge ``{NP}-hard problems: Maximum clique, graph coloring, and
                  satisfiability''},
  note =	 {{Instances} available at
                  \url{http://archive.dimacs.rutgers.edu/pub/challenge/graph/benchmarks/clique/}},
  year =	 {1992},
  key =		 {DIMACS}
}

@article{claus_danielson_fundamental_2021,
	title = {Fundamental {Domains} for {Symmetric} {Optimization} {Construction} and {Search}},
	volume = {31},
	journal = {SIAM J. on Optimization},
	author = {{Claus Danielson}},
	year = {2021},
	pages = {1827--1849},
}

@Misc{ECC,
  key =		 {ECC},
  title =	 {Error Correcting Codes Instances},
  year =	 {2001},
  note =	 {Available at
                  \url{https://github.com/jamestrimble/max-weight-clique-instances/tree/master/error-correcting-codes}}
}

@Misc{Kidney,
  key =		 {Kindney},
  title =	 {Kidney Exchange Instances},
  year =	 {2016},
  note =	 {Available at
                  \url{https://github.com/jamestrimble/max-weight-clique-instances/tree/master/kidney-exchange}}
}

@Misc{OEIS,
  key =		 {OEIS},
  author =	 {N. J. A. Sloane},
  title =	 {Challenge Problems: Independent Sets in Graphs},
  note =	 {Available at \url{https://oeis.org/A265032/a265032.html}}
}

@Misc{CliSAT,
  key =		 {CliSAT},
  title =	 {{CliSAT} -- an efficient state-of-the-art exact algorithm for the Maximum Clique
                  Problem},
  note =	 {Available at \url{https://github.com/psanse/CliSAT}}
}

@Misc{MaxClique,
  key =		 {MaxClique},
  author =	 {James Trimble},
  title =	 {Maximum weight clique instances},
  year =	 {2017},
  note =	 {Available at
                  \url{https://github.com/jamestrimble/max-weight-clique-instances/tree/master}}
}

@Misc{PACS19,
  key =		 {PACE19},
  author =	 {Max Bannach and Sebastian Berndt and Holger Dell and Bart M. P. Jansen and Philipp
                  Kindermann and Andr\'e Nichterlein and Christian Schulz and Soeren Terziadis },
  title =	 {PACE Challenge 2019 -- Parameterized Algorithms and Computational Experiments},
  year =	 {2019},
  note =	 {Available at \url{https://pacechallenge.org/2019/vc/vc_exact/}}
}

@Misc{UAI,
  key =		 {UAI},
  author =	 {James Trimble},
  title =	 {{UAI} Competition 2014 Instances},
  note =	 {Available at
                  \url{https://github.com/jamestrimble/max-weight-clique-instances/tree/master/UAI},
                  original source: \url{https://github.com/dechterlab/uai-competitions}}
}

@InCollection{Friedman2007,
  Title =	 {Fundamental Domains for Integer Programs with Symmetries},
  Author =	 {Friedman, Eric J.},
  Booktitle =	 {Combinatorial Optimization and Applications},
  Publisher =	 {Springer},
  Year =	 {2007},
  Editor =	 {Andreas Dress and Yinfeng Xu and Binhai Zhu},
  Pages =	 {146--153},
  Series =	 {LNCS},
  Volume =	 {4616},
  Doi =		 {10.1007/978-3-540-73556-4_17},
  noUrl =	 {http://dx.doi.org/10.1007/978-3-540-73556-4_17}
}

@Article{HojnyPfetsch2019,
  author =	 {Christopher Hojny and Marc E. Pfetsch},
  title =	 {Polytopes associated with symmetry handling},
  journal =	 {Mathematical Programming},
  year =	 {2019},
  pages =	 {197--240},
  volume =	 {175},
  issn =	 {1436-4646},
  doi =		 {10.1007/s10107-018-1239-7},
  nourl =	 {https://doi.org/10.1007/s10107-018-1239-7}
}

@article{KaibelPfetsch2008,
  year =	 {2008},
  issn =	 {0025-5610},
  journal =	 {Mathematical Programming},
  volume =	 {114},
  number =	 {1},
  doi =		 {10.1007/s10107-006-0081-5},
  title =	 {Packing and partitioning orbitopes},
  nourl =	 {http://dx.doi.org/10.1007/s10107-006-0081-5},
  publisher =	 {Springer},
  author =	 {Volker Kaibel and Marc E. Pfetsch},
  pages =	 {1--36}
}

@article{Liberti2012a,
  year =	 {2012},
  issn =	 {0025-5610},
  journal =	 {Mathematical Programming},
  volume =	 {131},
  number =	 {1-2},
  doi =		 {10.1007/s10107-010-0351-0},
  title =	 {Reformulations in mathematical programming: automatic symmetry detection and
                  exploitation},
  nourl =	 {http://dx.doi.org/10.1007/s10107-010-0351-0},
  author =	 {Liberti, Leo},
  pages =	 {273--304}
}

@Article{LibO14,
  author =	 {Leo Liberti and Jim Ostrowski},
  title =	 {Stabilizer-based Symmetry Breaking constraints for Mathematical Programs},
  journal =	 {Journal of Global Optimization},
  year =	 {2014},
  volume =	 {60},
  pages =	 {183--194}
}

@article{Mar02,
  author =	 {Margot, Fran{\c{c}}ois},
  year =	 {2002},
  journal =	 {Mathematical Programming},
  volume =	 {94},
  number =	 {1},
  doi =		 {10.1007/s10107-002-0358-2},
  title =	 {Pruning by isomorphism in branch-and-cut},
  pages =	 {71--90}
}

@article{Mar03a,
  author =	 {Margot, Fran{\c{c}}ois},
  year =	 {2003},
  journal =	 {Mathematical Programming},
  volume =	 {98},
  number =	 {1--3},
  doi =		 {10.1007/s10107-003-0394-6},
  title =	 {Exploiting orbits in symmetric {ILP}},
  pages =	 {3--21}
}

@Article{NemhauserTrotter1975,
  author =	 {George L Nemhauser and Leslie E Trotter},
  title =	 {Vertex packings: {Structural} properties and algorithms},
  journal =	 {Mathematical Programming},
  year =	 {1975},
  volume =	 {8},
  pages =	 {232--248}
}

@Article{Padberg1973,
  Title =	 {On the facial structure of set packing polyhedra},
  Author =	 {Manfred Padberg},
  Journal =	 {Mathematical Programming},
  Year =	 {1973},
  Pages =	 {199--215},
  Volume =	 {5}
}

@InProceedings{Sal18,
  author =	 {Salvagnin, Domenico},
  editor =	 {van Hoeve, Willem-Jan},
  title =	 {Symmetry Breaking Inequalities from the {S}chreier-{S}ims Table},
  booktitle =	 {Integration of Constraint Programming, Artificial Intelligence, and Operations
                  Research},
  year =	 {2018},
  publisher =	 {Springer},
  address =	 {Cham},
  pages =	 {521--529}
}

@Book{Seress2003,
  author =	 {\'Akos Seress},
  title =	 {Permutation Group Algorithms},
  publisher =	 {Cambridge University Press},
  year =	 {2003}
}

@article{verschae_geometry_2023,
  title =	 {On the geometry of symmetry breaking inequalities},
  volume =	 {197},
  number =	 {2},
  pages =	 {693--719},
  journal =	 {Mathematical Programming},
  author =	 {Verschae, José and Villagra, Matías and von Niederhäusern, Léonard},
  doi =		 {10.1007/s10107-022-01819-2},
  year =	 {2023}
}

@incollection{FonluptUhry1982,
  title =	 {Transformations which Preserve Perfectness and H-Perfectness of Graphs},
  editor =	 {Achim Bachem and Martin Gr\"otschel and Bernhard Korte},
  series =	 {North-Holland Mathematics Studies},
  publisher =	 {North-Holland},
  volume =	 {66},
  pages =	 {83--95},
  year =	 {1982},
  booktitle =	 {Bonn Workshop on Combinatorial Optimization},
  issn =	 {0304-0208},
  doi =		 {10.1016/S0304-0208(08)72445-9},
  author =	 {J. Fonlupt and J.P. Uhry}
}

@Book{GrotschelLovaszSchrijver1993,
  author =	 {Martin Gr{\"o}tschel and L{\'a}szl{\'o} Lov{\'a}sz and Alexander Schrijver},
  title =	 {Geometric Algorithms and Combinatorial Optimization},
  publisher =	 {Springer},
  year =	 {1993}
}

@incollection{Karp72,
  author =	 {R. M. Karp},
  title =	 {Reducibility among combinatorial problems},
  booktitle =	 {Complexity of Computer Computations},
  editor =	 {R. Miller and J. Thatcher},
  pages =	 {85--103},
  publisher =	 {Plenum Press},
  year =	 {1972}
}

@article{Tro75,
  title =	 {A class of facet producing graphs for vertex packing polyhedra},
  journal =	 {Discrete Mathematics},
  volume =	 {12},
  number =	 {4},
  pages =	 {373--388},
  year =	 {1975},
  doi =		 {10.1016/0012-365X(75)90077-1},
  author =	 {L.E. Trotter}
}

@article{BestuzhevaEtAl23,
  author =	 {Bestuzheva, Ksenia and Besan\c{c}on, Mathieu and Chen, Wei-Kun and Chmiela,
                  Antonia and Donkiewicz, Tim and van Doornmalen, Jasper and Eifler, Leon and Gaul,
                  Oliver and Gamrath, Gerald and Gleixner, Ambros and Gottwald, Leona and Graczyk,
                  Christoph and Halbig, Katrin and Hoen, Alexander and Hojny, Christopher and van
                  der Hulst, Rolf and Koch, Thorsten and L\"{u}bbecke, Marco and Maher, Stephen
                  J. and Matter, Frederic and M\"{u}hmer, Erik and M\"{u}ller, Benjamin and Pfetsch,
                  Marc E. and Rehfeldt, Daniel and Schlein, Steffan and Schl\"{o}sser, Franziska and
                  Serrano, Felipe and Shinano, Yuji and Sofranac, Boro and Turner, Mark and
                  Vigerske, Stefan and Wegscheider, Fabian and Wellner, Philipp and Weninger, Dieter
                  and Witzig, Jakob},
  title =	 {Enabling Research Through the {SCIP Optimization Suite} 8.0},
  year =	 {2023},
  volume =	 {49},
  number =	 {2},
  note =	 {Article 22},
  doi =		 {10.1145/3585516},
  journal =	 {ACM Trans. Math. Softw.}
}

@phdthesis{Borndoerfer1998,
  author =	 {Ralf Bornd{\"o}rfer},
  title =	 {Aspects of Set Packing, Partitioning, and Covering},
  school =	 {TU Berlin},
  url =		 {http://opus.kobv.de/zib/volltexte/1998/1028/},
  year =	 {1998}
}

@article{GroetschelLovaszSchrijver1984,
  year =	 {1984},
  journal =	 {Annals of Discrete Mathematics},
  volume =	 {21},
  title =	 {Polynomial algorithms for perfect graphs},
  author =	 {Martin Gr\"otschel and L\'asl\'o Lov\'asz and Alexander Schrijver},
  pages =	 {325--356}
}

@article{Hojny2020a,
  title =	 {Polynomial size IP formulations of knapsack may require exponentially large
                  coefficients},
  journal =	 {Operations Research Letters},
  volume =	 {48},
  number =	 {5},
  pages =	 {612--618},
  year =	 {2020},
  issn =	 {0167-6377},
  doi =		 {https://doi.org/10.1016/j.orl.2020.07.013},
  url =		 {https://www.sciencedirect.com/science/article/pii/S0167637720301103},
  author =	 {Christopher Hojny}
}

@techreport{SCIPOptSuite10,
  author = {Christopher Hojny and Mathieu Besan{\c{c}}on and Ksenia Bestuzheva and Sander Borst and Antonia Chmiela and Jo{\~{a}}o Dion{\'{i}}sio and Leon Eifler and Mohammed Ghannam and Ambros Gleixner and Adrian G\"{o}{\ss} and Alexander Hoen and Rolf van der Hulst and Dominik Kamp and Thorsten Koch and Kevin Kofler and Jurgen Lentz and Stephen J. Maher and Gioni Mexi and Erik M\"{u}hmer and Marc E. Pfetsch and Sebastian Pokutta and Felipe Serrano and Yuji Shinano and Mark Turner and Stefan Vigerske and Matthias Walter and Dieter Weninger and Liding Xu},
  title = {{The SCIP Optimization Suite 10.0}},
  type = {Technical Report},
  institution = {Optimization Online},
  month = {November},
  year = {2025},
  url = {https://optimization-online.org/2025/11/the-scip-optimization-suite-10-0/}
}

@article{BendottiEtAl2021,
  Title =	 {Orbitopal fixing for the full (sub-)orbitope and application to the Unit
                  Commitment Problem},
  Author =	 {Pascale Bendotti and Pierre Fouilhoux and C\'ecile Rottner},
  journal =	 {Mathematical Programming},
  volume =	 {186},
  year =	 {2021},
  pages =	 {337--372},
  doi =		 {10.1007/s10107-019-01457-1}
}

@article{KaibelEtAl2011,
  title =	 {Orbitopal fixing},
  journal =	 {Discrete Optimization},
  volume =	 {8},
  number =	 {4},
  pages =	 {595--610},
  year =	 {2011},
  issn =	 {1572-5286},
  doi =		 {10.1016/j.disopt.2011.07.001},
  nourl =	 {http://www.sciencedirect.com/science/article/pii/S1572528611000430},
  author =	 {Volker Kaibel and Matthias Peinhardt and Marc E. Pfetsch}
}

@InCollection{KaibelLoos2011,
  Title =	 {Finding Descriptions of Polytopes via Extended Formulations and Liftings},
  Author =	 {Volker Kaibel and Andreas Loos},
  Booktitle =	 {Progress in Combinatorial Optimization},
  Publisher =	 {Wiley},
  Year =	 {2011},
  Editor =	 {A. Ridha Mahjoub}
}

@article{vDoornmaalenHojny2022,
  author =	 {van Doornmalen, Jasper and Hojny, Christopher},
  title =	 {Efficient Propagation Techniques for Handling Cyclic Symmetries in Binary
                  Programs},
  journal =	 {INFORMS Journal on Computing},
  volume =	 {36},
  number =	 {3},
  pages =	 {868-883},
  year =	 {2024},
  doi =		 {10.1287/ijoc.2022.0060}
}

@PhdThesis{Ostrowski2008,
  author =	 {James Ostrowski},
  title =	 {Symmetry in Integer Programming},
  school =	 {Lehigh University},
  year =	 {2008}
}

@article{OstrowskiEtAl2011,
  year =	 {2011},
  issn =	 {0025-5610},
  journal =	 {Mathematical Programming},
  volume =	 {126},
  number =	 {1},
  doi =		 {10.1007/s10107-009-0273-x},
  title =	 {Orbital branching},
  nourl =	 {http://dx.doi.org/10.1007/s10107-009-0273-x},
  publisher =	 {Springer-Verlag},
  author =	 {Ostrowski, James and Linderoth, Jeff and Rossi, Fabrizio and Smriglio, Stefano},
  pages =	 {147--178}
}

@article{OstrowskiAnjosVannelli2015,
  year =	 {2015},
  journal =	 {Mathematical Programming},
  volume =	 {150},
  number =	 {1},
  doi =		 {10.1007/s10107-014-0812-y},
  title =	 {Modified orbital branching for structured symmetry with an application to unit
                  commitment},
  nourl =	 {http://dx.doi.org/10.1007/s10107-014-0812-y},
  publisher =	 {Springer Berlin Heidelberg},
  author =	 {James Ostrowski and Miguel F. Anjos and Anthony Vannelli},
  pages =	 {99--129}
}

@article{LinderothEtAl2021,
  author =	 {Linderoth, Jeff and N\'{u}\~{n}ez Ares, Jos\'{e} and Ostrowski, James and Rossi,
                  Fabrizio and Smriglio, Stefano},
  title =	 {Orbital Conflict: Cutting Planes for Symmetric Integer Programs},
  journal =	 {INFORMS Journal on Optimization},
  volume =	 {3},
  number =	 {2},
  pages =	 {139--153},
  year =	 {2021},
  doi =		 {10.1287/ijoo.2019.0044},
  URL =		 {https://doi.org/10.1287/ijoo.2019.0044},
  eprint =	 {https://doi.org/10.1287/ijoo.2019.0044}
}

@Misc{bliss,
  key =		 {bliss},
  author =	 {T. Junttila and P. Kaski},
  title =	 {bliss: A Tool for Computing Automorphism Groups and Canonical Labelings of Graphs},
  howpublished = {\url{https://users.aalto.fi/~tjunttil/bliss/}},
  year =	 {2012}
}

@InCollection{Sew96,
  author =	 {E. C. Sewell},
  title =	 {An improved algorithm for exact graph coloring},
  booktitle =	 {Cliques, coloring, and satisfiability. Second DIMACS implementation
                  challenge. Proceedings of a workshop held at DIMACS, 1993},
  pages =	 {359--373},
  publisher =	 {AMS, DIMACS},
  year =	 {1996},
  editor =	 {David S. Johnson and Michael Trick},
  volume =	 {26},
  series =	 {Ser. Discrete Math. Theor. Comput. Sci.}
}

@incollection{Johnson1982,
  author = {Johnson, Ellis L. and Padberg, Manfred W.},
  editor = {Achim Bachem and Martin Grötschel and Bemhard Korte},
  title = {Degree-two Inequalities, Clique Facets, and Biperfect Graphs},
  booktitle = {Bonn Workshop on Combinatorial Optimization},
  year = {1982},
  publisher = {North-Holland},
  pages = {169--187},
  doi = {https://doi.org/10.1016/S0304-0208(08)72450-2},
  series = {North-Holland Mathematics Studies},
  volume = {66},
  issn = {0304-0208},
  url = {https://www.sciencedirect.com/science/article/pii/S0304020808724502}
}

@article{Tamura1997,
  title = {The generalized stable set problem for perfect bidirected graphs},
  author = {Akihisa Tamura},
year = {1997},
doi = {10.15807/jorsj.40.401},
volume = {40},
pages = {401--414},
journal = {Journal of the Operations Research Society of Japan},
issn = {0453-4514},
publisher = {Operations Research Society of Japan},
number = {3}
}

@unpublished{Ikebe1996,
  author       = {Yoshiko T. Ikebe and Akihisa Tamura},
  title        = {Perfect bidirected graphs},
  note         = {Report CSIM96-2, Department of Computer Science and Information Mathematics, University of Electro-Communications, Tokyo},
  year         = {1996}
}

@phdthesis{vDoornmalen2024,
  title =	 "Propagation Algorithms for Handling Symmetries of Mathematical Programs",
  author =	 "van Doornmalen, Jasper",
  year =	 "2024",
  isbn =	 "978-90-386-6041-7",
  school =	 "Eindhoven University of Technology"
}

@InProceedings{andSS23,
  author =	 {Anders, Markus and Schweitzer, Pascal and Stie{\ss}, Julian},
  title =	 {Engineering a Preprocessor for Symmetry Detection},
  booktitle =	 {21st International Symposium on Experimental Algorithms (SEA 2023)},
  pages =	 {1:1--1:21},
  series =	 {Leibniz International Proceedings in Informatics (LIPIcs)},
  year =	 {2023},
  volume =	 {265},
  editor =	 {Georgiadis, Loukas},
  publisher =	 {Schloss Dagstuhl -- Leibniz-Zentrum f{\"u}r Informatik},
  address =	 {Dagstuhl, Germany},
  doi =		 {10.4230/LIPIcs.SEA.2023.1}
}

@inProceedings{JaegerPfetsch2025,
    author = {J\"{a}ger, Annika and Pfetsch, Marc E.},
    title = {Structure-Preserving Symmetry Presolving for Mixed-Binary Linear Problems},
    year = {2026},
    isbn = {978-3-032-28690-1},
    publisher = {Springer-Verlag},
    address = {Berlin, Heidelberg},
    url = {https://doi.org/10.1007/978-3-032-28691-8_30},
    doi = {10.1007/978-3-032-28691-8_30},
    booktitle = {Integer Programming and Combinatorial Optimization: 27th International Conference, IPCO 2026, Padua, Italy, June 17–19, 2026, Proceedings},
    pages = {457--473},
    numpages = {17},
    location = {Padua, Italy}
}

\clearpage
\appendix

\section{Additional Details for the Computational Experiments}
\label{sec:AdditionalComputationalResults}

In this appendix, we provide detailed numerical results for each of the instances in the Mixed-testset
used in our experiments, see Table~\ref{tab:instance_statistics}.
Each row contains information about one instance that appears in the Mixed-testset in at least one version out of the original and the two permuted versions.
If it is contained in the testset in several versions, the symtime shown is the maximum over all symtimes.
The number of generators is exemplary for one version and all the other information given is static among the three different versions of one instance.

\smallskip

\begin{scriptsize}
  \tablecaption{Instance statistics: source, number of nodes (``$\card{V}$'') and edges (``$\card{E}$'') in the graph,
      the number of generators (``\# gen'') and size (``$\card{\group}$'') of the automorphism group $\group$, the time for symmetry computation (``symtime'') in seconds, and the number of permutations of the instance contained in the Mixed-testset (``\# perm'').}
  \tablehead{
    \toprule
    name &
    source &
    $\card{V}$ &
    $\card{E}$ &
    \# gen &
    $\card{\group}$ &
    symtime &
    \# perm
    \\ 
    \midrule}
  \tabletail{
    \midrule
    continued on next page \\
    \bottomrule
    }
  \tablelasttail{\bottomrule}
  \label{tab:instance_statistics}
  \begin{supertabular*}{\textwidth}{@{}l@{\;\extracolsep{\fill}}lrrrrrr@{}}
    aim-50-1-6-unsat-2                  &    aim\_A & \num{     536} & \num{   11863} & \num{     1} & \num{     2} & \num{  0.00} & \num{1}\\
    1-Insertions\_5                     &  Coloring & \num{     202} & \num{    1227} & \num{     2} & \num{    10} & \num{  0.00} & \num{3}\\
    DSJR500.1c                          &  Coloring & \num{     259} & \num{   32261} & \num{     2} & \num{     6} & \num{  0.00} & \num{3}\\
    myciel7                             &  Coloring & \num{     191} & \num{    2360} & \num{     2} & \num{    10} & \num{  0.00} & \num{3}\\
    queen10\_10                         &  Coloring & \num{     100} & \num{    1470} & \num{     2} & \num{     8} & \num{  0.00} & \num{3}\\
    c-fat200-5                          &    Dimacs & \num{     200} & \num{   11427} & \num{   193} & \num{3e+212} & \num{  0.08} & \num{3}\\
    c-fat500-1                          &    Dimacs & \num{     500} & \num{  120291} & \num{   420} & \num{6e+245} & \num{  1.90} & \num{3}\\
    c-fat500-10                         &    Dimacs & \num{     500} & \num{   78123} & \num{   492} & $\infty$ & \num{  0.00} & \num{3}\\
    c-fat500-2                          &    Dimacs & \num{     500} & \num{  115611} & \num{   460} & $\infty$ & \num{  0.01} & \num{3}\\
    c-fat500-5                          &    Dimacs & \num{     500} & \num{  101559} & \num{   484} & $\infty$ & \num{  1.73} & \num{3}\\
    hamming10-4                         &    Dimacs & \num{    1024} & \num{   89600} & \num{    10} & \num{ 4e+09} & \num{  0.01} & \num{3}\\
    keller4                             &    Dimacs & \num{     171} & \num{    5100} & \num{     4} & \num{ 4e+02} & \num{  0.00} & \num{3}\\
    keller5                             &    Dimacs & \num{     776} & \num{   74710} & \num{     5} & \num{ 4e+03} & \num{  0.00} & \num{3}\\
    keller6                             &    Dimacs & \num{    3361} & \num{ 1026582} & \num{     6} & \num{ 5e+04} & \num{  0.03} & \num{3}\\
    MANN\_a27                           &    Dimacs & \num{     378} & \num{     702} & \num{     9} & \num{ 3e+05} & \num{  0.00} & \num{3}\\
    MANN\_a45                           &    Dimacs & \num{    1035} & \num{    1980} & \num{     5} & \num{ 4e+02} & \num{  0.00} & \num{3}\\
    MANN\_a81                           &    Dimacs & \num{    3321} & \num{    6480} & \num{    14} & \num{ 2e+09} & \num{  0.00} & \num{3}\\
    01-11-4-4                           &       ECC & \num{     150} & \num{    9365} & \num{     8} & \num{ 2e+03} & \num{  0.00} & \num{3}\\
    02-12-4-6                           &       ECC & \num{     230} & \num{   22453} & \num{     4} & \num{ 4e+02} & \num{  0.00} & \num{3}\\
    03-14-4-7                           &       ECC & \num{     223} & \num{   19900} & \num{     3} & \num{    10} & \num{  0.00} & \num{3}\\
    04-14-6-6                           &       ECC & \num{     807} & \num{  186477} & \num{     6} & \num{ 3e+03} & \num{  0.01} & \num{3}\\
    05-16-4-5                           &       ECC & \num{     156} & \num{    9730} & \num{     3} & \num{     8} & \num{  0.00} & \num{3}\\
    06-16-8-8                           &       ECC & \num{    2246} & \num{  349235} & \num{   118} & \num{ 3e+37} & \num{  0.01} & \num{3}\\
    07-17-4-4                           &       ECC & \num{     132} & \num{    5892} & \num{     2} & \num{     8} & \num{  0.00} & \num{2}\\
    08-17-6-6                           &       ECC & \num{     558} & \num{   95583} & \num{     6} & \num{ 3e+03} & \num{  0.01} & \num{3}\\
    09-19-4-6                           &       ECC & \num{     263} & \num{   29080} & \num{     4} & \num{    70} & \num{  0.00} & \num{3}\\
    10-19-8-8                           &       ECC & \num{    2124} & \num{  595392} & \num{     5} & \num{ 4e+02} & \num{  0.02} & \num{3}\\
    11-20-6-5                           &       ECC & \num{    1302} & \num{  344541} & \num{     3} & \num{    80} & \num{  0.01} & \num{3}\\
    12-20-6-6                           &       ECC & \num{    1490} & \num{  428359} & \num{     3} & \num{     8} & \num{  0.01} & \num{3}\\
    13-20-8-10                          &       ECC & \num{    2510} & \num{  590958} & \num{     3} & \num{    20} & \num{  0.02} & \num{3}\\
    14-21-10-9                          &       ECC & \num{    5098} & \num{ 2867431} & \num{   327} & \num{8e+101} & \num{  0.11} & \num{3}\\
    15-22-10-10                         &       ECC & \num{    8914} & \num{ 2694426} & \num{     5} & \num{ 2e+02} & \num{  0.09} & \num{3}\\
    01-11-4-4.cmpl                      &     ECC-C & \num{     150} & \num{    1810} & \num{     5} & \num{ 2e+03} & \num{  0.00} & \num{3}\\
    02-12-4-6.cmpl                      &     ECC-C & \num{     230} & \num{    3882} & \num{     4} & \num{ 4e+02} & \num{  0.00} & \num{3}\\
    03-14-4-7.cmpl                      &     ECC-C & \num{     223} & \num{    4853} & \num{     2} & \num{    10} & \num{  0.00} & \num{3}\\
    04-14-6-6.cmpl                      &     ECC-C & \num{     807} & \num{  138744} & \num{     5} & \num{ 3e+03} & \num{  0.01} & \num{3}\\
    05-16-4-5.cmpl                      &     ECC-C & \num{     156} & \num{    2360} & \num{     3} & \num{     8} & \num{  0.00} & \num{3}\\
    06-16-8-8.cmpl                      &     ECC-C & \num{    2246} & \num{ 2171900} & \num{   117} & \num{ 3e+37} & \num{  0.05} & \num{3}\\
    08-17-6-6.cmpl                      &     ECC-C & \num{     558} & \num{   59820} & \num{     7} & \num{ 3e+03} & \num{  0.01} & \num{3}\\
    09-19-4-6.cmpl                      &     ECC-C & \num{     263} & \num{    5373} & \num{     4} & \num{    70} & \num{  0.00} & \num{3}\\
    10-19-8-8.cmpl                      &     ECC-C & \num{    2124} & \num{ 1659234} & \num{     6} & \num{ 4e+02} & \num{  0.04} & \num{3}\\
    11-20-6-5.cmpl                      &     ECC-C & \num{    1302} & \num{  502410} & \num{     4} & \num{    80} & \num{  0.02} & \num{3}\\
    12-20-6-6.cmpl                      &     ECC-C & \num{    1490} & \num{  680946} & \num{     2} & \num{     8} & \num{  0.02} & \num{3}\\
    13-20-8-10.cmpl                     &     ECC-C & \num{    2510} & \num{ 2557837} & \num{     3} & \num{    20} & \num{  0.07} & \num{3}\\
    14-21-10-9.cmpl                     &     ECC-C & \num{    5098} & \num{10124822} & \num{   327} & \num{8e+101} & \num{  0.11} & \num{3}\\
    15-22-10-10.cmpl                    &     ECC-C & \num{    8914} & \num{37030815} & \num{     5} & \num{ 2e+02} & \num{  1.23} & \num{3}\\
    01-11-4-4-compl-weigh               &    ECC-CW & \num{     150} & \num{    1810} & \num{     8} & \num{ 2e+03} & \num{  0.00} & \num{3}\\
    02-12-4-6-compl-weigh               &    ECC-CW & \num{     230} & \num{    3882} & \num{     4} & \num{ 4e+02} & \num{  0.00} & \num{3}\\
    03-14-4-7-compl-weigh               &    ECC-CW & \num{     223} & \num{    4853} & \num{     3} & \num{    10} & \num{  0.00} & \num{3}\\
    04-14-6-6-compl-weigh               &    ECC-CW & \num{     807} & \num{  138744} & \num{     5} & \num{ 3e+03} & \num{  0.01} & \num{3}\\
    05-16-4-5-compl-weigh               &    ECC-CW & \num{     156} & \num{    2360} & \num{     3} & \num{     8} & \num{  0.00} & \num{3}\\
    06-16-8-8-compl-weigh               &    ECC-CW & \num{    2246} & \num{ 2171900} & \num{   118} & \num{ 3e+37} & \num{  0.05} & \num{3}\\
    07-17-4-4-compl-weigh               &    ECC-CW & \num{     132} & \num{    2754} & \num{     3} & \num{     8} & \num{  0.00} & \num{2}\\
    08-17-6-6-compl-weigh               &    ECC-CW & \num{     558} & \num{   59820} & \num{     6} & \num{ 3e+03} & \num{  0.01} & \num{3}\\
    09-19-4-6-compl-weigh               &    ECC-CW & \num{     263} & \num{    5373} & \num{     4} & \num{    70} & \num{  0.00} & \num{3}\\
    10-19-8-8-compl-weigh               &    ECC-CW & \num{    2124} & \num{ 1659234} & \num{     5} & \num{ 4e+02} & \num{  0.04} & \num{3}\\
    11-20-6-5-compl-weigh               &    ECC-CW & \num{    1302} & \num{  502410} & \num{     3} & \num{    80} & \num{  0.01} & \num{3}\\
    12-20-6-6-compl-weigh               &    ECC-CW & \num{    1490} & \num{  680946} & \num{     3} & \num{     8} & \num{  0.01} & \num{3}\\
    13-20-8-10-compl-weigh              &    ECC-CW & \num{    2510} & \num{ 2557837} & \num{     3} & \num{    20} & \num{  0.07} & \num{3}\\
    14-21-10-9-compl-weigh              &    ECC-CW & \num{    5098} & \num{10124822} & \num{   327} & \num{8e+101} & \num{  0.11} & \num{3}\\
    ehi-85-297-00                       &    ehi\_A & \num{    2079} & \num{  108240} & \num{    12} & \num{ 4e+03} & \num{  0.01} & \num{3}\\
    ehi-85-297-12                       &    ehi\_A & \num{    2079} & \num{  108302} & \num{    24} & \num{ 2e+07} & \num{  0.01} & \num{3}\\
    ehi-85-297-23                       &    ehi\_A & \num{    2079} & \num{  108643} & \num{     6} & \num{    60} & \num{  0.01} & \num{3}\\
    ehi-85-297-28                       &    ehi\_A & \num{    2079} & \num{  108331} & \num{    24} & \num{ 2e+07} & \num{  0.01} & \num{3}\\
    ehi-85-297-36                       &    ehi\_A & \num{    2079} & \num{  108481} & \num{    26} & \num{ 5e+08} & \num{  0.01} & \num{3}\\
    ehi-85-297-44                       &    ehi\_A & \num{    2079} & \num{  108730} & \num{    12} & \num{ 4e+03} & \num{  0.00} & \num{3}\\
    ehi-85-297-52                       &    ehi\_A & \num{    2079} & \num{  108346} & \num{    12} & \num{ 4e+03} & \num{  0.00} & \num{3}\\
    ehi-85-297-60                       &    ehi\_A & \num{    2079} & \num{  108541} & \num{    12} & \num{ 4e+03} & \num{  0.01} & \num{3}\\
    ehi-85-297-76                       &    ehi\_A & \num{    2079} & \num{  108401} & \num{    12} & \num{ 4e+03} & \num{  0.00} & \num{3}\\
    ehi-85-297-92                       &    ehi\_A & \num{    2079} & \num{  108448} & \num{     6} & \num{    60} & \num{  0.00} & \num{3}\\
    ehi-90-315-00                       &    ehi\_B & \num{    2205} & \num{  114973} & \num{    12} & \num{ 4e+03} & \num{  0.01} & \num{3}\\
    ehi-90-315-08                       &    ehi\_B & \num{    2205} & \num{  115260} & \num{    12} & \num{ 4e+03} & \num{  0.01} & \num{3}\\
    ehi-90-315-40                       &    ehi\_B & \num{    2205} & \num{  115449} & \num{     6} & \num{    60} & \num{  0.01} & \num{3}\\
    ehi-90-315-48                       &    ehi\_B & \num{    2205} & \num{  115451} & \num{     6} & \num{    60} & \num{  0.00} & \num{3}\\
    ehi-90-315-56                       &    ehi\_B & \num{    2205} & \num{  115388} & \num{     6} & \num{    60} & \num{  0.01} & \num{3}\\
    ehi-90-315-60                       &    ehi\_B & \num{    2205} & \num{  115303} & \num{    19} & \num{ 5e+05} & \num{  0.01} & \num{3}\\
    ehi-90-315-76                       &    ehi\_B & \num{    2205} & \num{  115188} & \num{     6} & \num{    60} & \num{  0.00} & \num{3}\\
    ehi-90-315-92                       &    ehi\_B & \num{    2205} & \num{  115587} & \num{     6} & \num{    60} & \num{  0.01} & \num{3}\\
    101                                 &    Kidney & \num{    4741} & \num{ 1943309} & \num{   102} & \num{ 5e+30} & \num{  0.06} & \num{3}\\
    102                                 &    Kidney & \num{    3717} & \num{ 1187313} & \num{    82} & \num{ 7e+24} & \num{  0.03} & \num{3}\\
    103                                 &    Kidney & \num{    4673} & \num{ 1884498} & \num{    73} & \num{ 9e+21} & \num{  0.06} & \num{3}\\
    104                                 &    Kidney & \num{    4846} & \num{ 2035984} & \num{   123} & \num{ 1e+37} & \num{  0.05} & \num{3}\\
    105                                 &    Kidney & \num{    4663} & \num{ 1985826} & \num{   246} & \num{ 1e+74} & \num{  0.05} & \num{3}\\
    106                                 &    Kidney & \num{    3790} & \num{ 1318808} & \num{   141} & \num{ 1e+44} & \num{  0.04} & \num{3}\\
    107                                 &    Kidney & \num{    5207} & \num{ 2392844} & \num{   276} & \num{ 3e+84} & \num{  0.07} & \num{3}\\
    108                                 &    Kidney & \num{    5529} & \num{ 2709071} & \num{   258} & \num{ 5e+77} & \num{  0.08} & \num{3}\\
    109                                 &    Kidney & \num{    4490} & \num{ 1701998} & \num{   167} & \num{ 2e+50} & \num{  0.05} & \num{3}\\
    110                                 &    Kidney & \num{    4802} & \num{ 2023807} & \num{    64} & \num{ 3e+19} & \num{  0.06} & \num{3}\\
    111                                 &    Kidney & \num{    8953} & \num{ 7106080} & \num{   390} & \num{1e+124} & \num{  0.20} & \num{3}\\
    112                                 &    Kidney & \num{    8288} & \num{ 6213021} & \num{   363} & \num{2e+111} & \num{  0.18} & \num{3}\\
    113                                 &    Kidney & \num{    6870} & \num{ 3877305} & \num{   125} & \num{ 4e+37} & \num{  0.11} & \num{3}\\
    114                                 &    Kidney & \num{    8169} & \num{ 6337819} & \num{   659} & \num{7e+229} & \num{  0.17} & \num{3}\\
    115                                 &    Kidney & \num{    4934} & \num{ 2666778} & \num{   292} & \num{ 4e+95} & \num{  0.07} & \num{3}\\
    116                                 &    Kidney & \num{    5451} & \num{ 2665013} & \num{   238} & \num{ 2e+80} & \num{  0.07} & \num{3}\\
    117                                 &    Kidney & \num{    4979} & \num{ 2415559} & \num{   201} & \num{ 1e+63} & \num{  0.07} & \num{3}\\
    118                                 &    Kidney & \num{    7592} & \num{ 4753874} & \num{   250} & \num{ 2e+75} & \num{  0.14} & \num{3}\\
    119                                 &    Kidney & \num{    5218} & \num{ 2335980} & \num{    89} & \num{ 6e+26} & \num{  0.07} & \num{3}\\
    120                                 &    Kidney & \num{    8072} & \num{ 5214597} & \num{   260} & \num{ 4e+78} & \num{  0.15} & \num{3}\\
    101-compl-weigh                     &  Kidney-W & \num{    4741} & \num{ 1943309} & \num{   101} & \num{ 3e+30} & \num{  0.06} & \num{3}\\
    102-compl-weigh                     &  Kidney-W & \num{    3717} & \num{ 1187313} & \num{    51} & \num{ 2e+15} & \num{  0.04} & \num{3}\\
    103-compl-weigh                     &  Kidney-W & \num{    4673} & \num{ 1884498} & \num{    73} & \num{ 9e+21} & \num{  0.07} & \num{3}\\
    104-compl-weigh                     &  Kidney-W & \num{    4846} & \num{ 2035984} & \num{   123} & \num{ 1e+37} & \num{  0.07} & \num{3}\\
    105-compl-weigh                     &  Kidney-W & \num{    4663} & \num{ 1985826} & \num{   236} & \num{ 1e+71} & \num{  0.05} & \num{3}\\
    106-compl-weigh                     &  Kidney-W & \num{    3790} & \num{ 1318808} & \num{    87} & \num{ 2e+26} & \num{  0.04} & \num{3}\\
    107-compl-weigh                     &  Kidney-W & \num{    5207} & \num{ 2392844} & \num{   265} & \num{ 6e+79} & \num{  0.06} & \num{3}\\
    108-compl-weigh                     &  Kidney-W & \num{    5529} & \num{ 2709071} & \num{   247} & \num{ 2e+74} & \num{  0.09} & \num{3}\\
    109-compl-weigh                     &  Kidney-W & \num{    4490} & \num{ 1701998} & \num{   165} & \num{ 5e+49} & \num{  0.05} & \num{3}\\
    110-compl-weigh                     &  Kidney-W & \num{    4802} & \num{ 2023807} & \num{    61} & \num{ 2e+18} & \num{  0.07} & \num{3}\\
    111-compl-weigh                     &  Kidney-W & \num{    8953} & \num{ 7106080} & \num{   358} & \num{5e+110} & \num{  0.20} & \num{3}\\
    112-compl-weigh                     &  Kidney-W & \num{    8288} & \num{ 6213021} & \num{   347} & \num{3e+104} & \num{  0.18} & \num{3}\\
    113-compl-weigh                     &  Kidney-W & \num{    6870} & \num{ 3877305} & \num{   113} & \num{ 1e+34} & \num{  0.12} & \num{3}\\
    114-compl-weigh                     &  Kidney-W & \num{    8169} & \num{ 6337819} & \num{   582} & \num{3e+193} & \num{  0.18} & \num{3}\\
    115-compl-weigh                     &  Kidney-W & \num{    4934} & \num{ 2666778} & \num{   264} & \num{ 2e+86} & \num{  0.09} & \num{3}\\
    116-compl-weigh                     &  Kidney-W & \num{    5451} & \num{ 2665013} & \num{   185} & \num{ 5e+61} & \num{  0.08} & \num{3}\\
    117-compl-weigh                     &  Kidney-W & \num{    4979} & \num{ 2415559} & \num{   184} & \num{ 6e+56} & \num{  0.07} & \num{3}\\
    118-compl-weigh                     &  Kidney-W & \num{    7592} & \num{ 4753874} & \num{   243} & \num{ 1e+73} & \num{  0.14} & \num{3}\\
    119-compl-weigh                     &  Kidney-W & \num{    5218} & \num{ 2335980} & \num{    80} & \num{ 1e+24} & \num{  0.07} & \num{3}\\
    120-compl-weigh                     &  Kidney-W & \num{    8072} & \num{ 5214597} & \num{   250} & \num{ 2e+75} & \num{  0.15} & \num{3}\\
    monoton-7                           &  Monotone & \num{     343} & \num{   12348} & \num{     3} & \num{    10} & \num{  0.00} & \num{3}\\
    monoton-8                           &  Monotone & \num{     512} & \num{   24192} & \num{     3} & \num{    10} & \num{  0.00} & \num{3}\\
    monoton-9                           &  Monotone & \num{     729} & \num{   43740} & \num{     3} & \num{    10} & \num{  0.00} & \num{3}\\
    a265032\_1dc.1024                   &      OEIS & \num{    1024} & \num{   24063} & \num{     2} & \num{     4} & \num{  0.00} & \num{3}\\
    a265032\_1dc.2048                   &      OEIS & \num{    2048} & \num{   58367} & \num{     2} & \num{     4} & \num{  0.00} & \num{3}\\
    a265032\_1dc.512                    &      OEIS & \num{     512} & \num{    9727} & \num{     2} & \num{     4} & \num{  0.00} & \num{3}\\
    a265032\_1et.1024                   &      OEIS & \num{    1024} & \num{    9600} & \num{    24} & \num{ 2e+13} & \num{  0.00} & \num{3}\\
    a265032\_1et.2048                   &      OEIS & \num{    2048} & \num{   22528} & \num{    25} & \num{ 8e+14} & \num{  0.00} & \num{3}\\
    a265032\_1et.512                    &      OEIS & \num{     512} & \num{    4032} & \num{    20} & \num{ 2e+11} & \num{  0.00} & \num{3}\\
    a265032\_1tc.1024                   &      OEIS & \num{    1024} & \num{    7936} & \num{    11} & \num{ 2e+04} & \num{  0.00} & \num{3}\\
    a265032\_1tc.2048                   &      OEIS & \num{    2048} & \num{   18944} & \num{    11} & \num{ 3e+04} & \num{  0.00} & \num{3}\\
    a265032\_1tc.512                    &      OEIS & \num{     512} & \num{    3264} & \num{     9} & \num{ 4e+03} & \num{  0.00} & \num{3}\\
    a265032\_1zc.1024                   &      OEIS & \num{    1024} & \num{   16640} & \num{    10} & \num{ 7e+06} & \num{  0.00} & \num{3}\\
    a265032\_1zc.2048                   &      OEIS & \num{    2048} & \num{   39424} & \num{    11} & \num{ 8e+07} & \num{  0.00} & \num{3}\\
    a265032\_1zc.256                    &      OEIS & \num{     256} & \num{    2816} & \num{     5} & \num{ 8e+04} & \num{  0.00} & \num{1}\\
    a265032\_1zc.4096                   &      OEIS & \num{    4096} & \num{   92160} & \num{    12} & \num{ 1e+09} & \num{  0.00} & \num{3}\\
    a265032\_1zc.512                    &      OEIS & \num{     512} & \num{    6912} & \num{     9} & \num{ 7e+05} & \num{  0.00} & \num{3}\\
    a265032\_2dc.1024                   &      OEIS & \num{    1024} & \num{  169162} & \num{     2} & \num{     4} & \num{  0.01} & \num{3}\\
    a265032\_2dc.2048                   &      OEIS & \num{    2048} & \num{  504451} & \num{     2} & \num{     4} & \num{  0.02} & \num{3}\\
    a265032\_2dc.512                    &      OEIS & \num{     512} & \num{   54895} & \num{     2} & \num{     4} & \num{  0.00} & \num{3}\\
    ref-20-30-00-compl-weigh            &     REF-W & \num{     300} & \num{   21980} & \num{     1} & \num{     2} & \num{  0.00} & \num{1}\\
    ref-20-30-01-compl-weigh            &     REF-W & \num{     300} & \num{   20791} & \num{     1} & \num{     2} & \num{  0.01} & \num{3}\\
    ref-20-30-02-compl-weigh            &     REF-W & \num{     300} & \num{   21672} & \num{     1} & \num{     2} & \num{  0.01} & \num{1}\\
    ref-20-30-03-compl-weigh            &     REF-W & \num{     300} & \num{   20759} & \num{     1} & \num{     2} & \num{  0.01} & \num{2}\\
    ref-20-30-05-compl-weigh            &     REF-W & \num{     300} & \num{   21351} & \num{     2} & \num{     4} & \num{  0.00} & \num{1}\\
    ref-20-30-06-compl-weigh            &     REF-W & \num{     300} & \num{   22234} & \num{     2} & \num{     4} & \num{  0.00} & \num{1}\\
    ref-20-30-07-compl-weigh            &     REF-W & \num{     300} & \num{   21904} & \num{     2} & \num{     4} & \num{  0.00} & \num{3}\\
    ref-20-30-08-compl-weigh            &     REF-W & \num{     300} & \num{   22820} & \num{     2} & \num{     4} & \num{  0.00} & \num{2}\\
    ref-20-40-02-compl-weigh            &     REF-W & \num{     300} & \num{   16759} & \num{     1} & \num{     2} & \num{  0.00} & \num{3}\\
    ref-20-50-01-compl-weigh            &     REF-W & \num{     300} & \num{   14276} & \num{     4} & \num{    20} & \num{  0.00} & \num{1}\\
    ref-20-50-06-compl-weigh            &     REF-W & \num{     300} & \num{   14805} & \num{     4} & \num{    20} & \num{  0.00} & \num{2}\\
    ref-25-25-00-compl-weigh            &     REF-W & \num{     375} & \num{   38311} & \num{     8} & \num{ 3e+02} & \num{  0.00} & \num{1}\\
    ref-25-25-01-compl-weigh            &     REF-W & \num{     375} & \num{   37695} & \num{    10} & \num{ 1e+03} & \num{  0.00} & \num{2}\\
    ref-25-25-03-compl-weigh            &     REF-W & \num{     375} & \num{   38642} & \num{     4} & \num{    20} & \num{  0.00} & \num{1}\\
    ref-25-25-04-compl-weigh            &     REF-W & \num{     375} & \num{   37159} & \num{     1} & \num{     2} & \num{  0.00} & \num{1}\\
    ref-25-25-05-compl-weigh            &     REF-W & \num{     375} & \num{   41484} & \num{     4} & \num{    20} & \num{  0.01} & \num{1}\\
    ref-25-25-07-compl-weigh            &     REF-W & \num{     375} & \num{   38464} & \num{     8} & \num{ 3e+02} & \num{  0.01} & \num{3}\\
    ref-25-25-09-compl-weigh            &     REF-W & \num{     375} & \num{   39210} & \num{     8} & \num{ 3e+02} & \num{  0.00} & \num{2}\\
    ref-25-30-00-compl-weigh            &     REF-W & \num{     375} & \num{   33508} & \num{     3} & \num{     8} & \num{  0.00} & \num{1}\\
    ref-25-30-01-compl-weigh            &     REF-W & \num{     375} & \num{   33475} & \num{     3} & \num{     8} & \num{  0.00} & \num{1}\\
    ref-25-30-02-compl-weigh            &     REF-W & \num{     375} & \num{   32584} & \num{     4} & \num{    20} & \num{  0.00} & \num{2}\\
    ref-25-30-03-compl-weigh            &     REF-W & \num{     375} & \num{   33416} & \num{     2} & \num{     4} & \num{  0.00} & \num{1}\\
    ref-25-30-04-compl-weigh            &     REF-W & \num{     375} & \num{   31898} & \num{     1} & \num{     2} & \num{  0.00} & \num{1}\\
    ref-25-30-07-compl-weigh            &     REF-W & \num{     375} & \num{   32967} & \num{     2} & \num{     4} & \num{  0.00} & \num{1}\\
    ref-25-30-08-compl-weigh            &     REF-W & \num{     375} & \num{   33444} & \num{     4} & \num{    20} & \num{  0.00} & \num{1}\\
    ref-25-30-09-compl-weigh            &     REF-W & \num{     375} & \num{   32687} & \num{     7} & \num{ 1e+02} & \num{  0.00} & \num{2}\\
    ref-25-35-00-compl-weigh            &     REF-W & \num{     375} & \num{   26794} & \num{     2} & \num{     4} & \num{  0.00} & \num{3}\\
    ref-25-35-01-compl-weigh            &     REF-W & \num{     375} & \num{   28579} & \num{     1} & \num{     2} & \num{  0.00} & \num{3}\\
    ref-25-35-02-compl-weigh            &     REF-W & \num{     375} & \num{   30076} & \num{     1} & \num{     2} & \num{  0.00} & \num{3}\\
    ref-25-35-03-compl-weigh            &     REF-W & \num{     375} & \num{   28611} & \num{     1} & \num{     2} & \num{  0.00} & \num{3}\\
    ref-25-35-04-compl-weigh            &     REF-W & \num{     375} & \num{   29161} & \num{     2} & \num{     4} & \num{  0.00} & \num{3}\\
    ref-25-35-07-compl-weigh            &     REF-W & \num{     375} & \num{   29135} & \num{     2} & \num{     4} & \num{  0.00} & \num{3}\\
    ref-25-35-08-compl-weigh            &     REF-W & \num{     375} & \num{   30297} & \num{     7} & \num{ 1e+02} & \num{  0.00} & \num{3}\\
    ref-25-40-05-compl-weigh            &     REF-W & \num{     375} & \num{   25834} & \num{     2} & \num{     4} & \num{  0.00} & \num{3}\\
    ref-25-40-06-compl-weigh            &     REF-W & \num{     375} & \num{   26071} & \num{     1} & \num{     2} & \num{  0.00} & \num{3}\\
    ref-30-30-00-compl-weigh            &     REF-W & \num{     450} & \num{   48085} & \num{     2} & \num{     4} & \num{  0.00} & \num{1}\\
    ref-30-30-01-compl-weigh            &     REF-W & \num{     450} & \num{   44844} & \num{     4} & \num{    20} & \num{  0.00} & \num{3}\\
    ref-30-30-02-compl-weigh            &     REF-W & \num{     450} & \num{   48996} & \num{     2} & \num{     4} & \num{  0.00} & \num{2}\\
    ref-30-30-03-compl-weigh            &     REF-W & \num{     450} & \num{   48186} & \num{     1} & \num{     2} & \num{  0.01} & \num{3}\\
    ref-30-30-04-compl-weigh            &     REF-W & \num{     450} & \num{   46130} & \num{     4} & \num{    20} & \num{  0.00} & \num{2}\\
    ref-30-30-05-compl-weigh            &     REF-W & \num{     450} & \num{   47447} & \num{     2} & \num{     4} & \num{  0.00} & \num{2}\\
    ref-30-30-06-compl-weigh            &     REF-W & \num{     450} & \num{   45528} & \num{     4} & \num{    20} & \num{  0.00} & \num{2}\\
    ref-30-30-07-compl-weigh            &     REF-W & \num{     450} & \num{   47063} & \num{     5} & \num{    30} & \num{  0.01} & \num{3}\\
    ref-30-30-08-compl-weigh            &     REF-W & \num{     450} & \num{   47452} & \num{     1} & \num{     2} & \num{  0.00} & \num{3}\\
    ref-30-30-09-compl-weigh            &     REF-W & \num{     450} & \num{   47274} & \num{     5} & \num{    30} & \num{  0.00} & \num{2}\\
    ref-30-40-04-compl-weigh            &     REF-W & \num{     450} & \num{   39187} & \num{     1} & \num{     2} & \num{  0.01} & \num{3}\\
    ref-30-40-06-compl-weigh            &     REF-W & \num{     450} & \num{   39999} & \num{     2} & \num{     4} & \num{  0.00} & \num{3}\\
    Alchemy\_11                         &       UAI & \num{    2360} & \num{    3630} & \num{     0} & \num{     1} & \num{  0.00} & \num{3}\\
    relational\_1                       &       UAI & \num{   62750} & \num{   62500} & \num{ 62499} & $\infty$ & \num{  2.50} & \num{3}\\
    relational\_2                       &       UAI & \num{   50750} & \num{   50500} & \num{ 24749} & $\infty$ & \num{  0.02} & \num{3}\\
    relational\_3                       &       UAI & \num{    9700} & \num{   18600} & \num{    99} & \num{9e+157} & \num{  0.06} & \num{3}\\
    relational\_4                       &       UAI & \num{   40200} & \num{   50000} & \num{  9999} & $\infty$ & \num{  0.01} & \num{3}\\
    relational\_5                       &       UAI & \num{   60000} & \num{   50000} & \num{ 49999} & $\infty$ & \num{  1.85} & \num{3}\\
    vc-exact\_038                       &        VC & \num{     786} & \num{   14024} & \num{    22} & \num{ 6e+08} & \num{  0.00} & \num{3}\\
    \midrule
    averages: & & \num{  3499.9} & \num{1213558.9} & \num{ 911.7} &  & \num{  0.06} & \num{526}\\
    \multicolumn{2}{@{}l}{(192 instances, 526 permutations)} & & & & & & \\
  \end{supertabular*}
\end{scriptsize}

\end{document}